\documentclass[11pt]{article}
\usepackage{amsmath}
\usepackage{amsfonts}
\usepackage{amssymb}
\usepackage{amsthm}
\usepackage{graphicx}
\usepackage{etoolbox,xstring,mfirstuc,textcase}
\usepackage{empheq}
\usepackage{indentfirst}
\usepackage{cite} 
\usepackage{mathrsfs}
\usepackage{cases}
\usepackage{graphics}
\usepackage{xcolor}
\usepackage{bm}
\usepackage{times}

\usepackage{hyperref}
\hypersetup{colorlinks=true, linkcolor=blue, citecolor=purple, urlcolor=blue}

\usepackage[left=2.5cm, right=2.5cm, top=2cm, bottom=2cm]{geometry}

\newtheorem{theorem}{\textbf{Theorem}}[section]
\newtheorem{lemma}{\textbf{Lemma}}[section]
\newtheorem{proposition}{\textbf{Proposition}}[section]
\newtheorem{corollary}{\textbf{Corollary}}[section]
\newtheorem{remark}{\textbf{Remark}}[section]
\newtheorem{definition}{\textbf{Definition}}[section]
\allowdisplaybreaks[4] 
\def\be{\begin{equation}}
\def\ee{\end{equation}}
\def\bea{\begin{eqnarray}}
\def\eea{\end{eqnarray}}
\def\bt{\begin{theorem}}
\def\et{\end{theorem}}
\def\bl{\begin{lemma}}
\def\el{\end{lemma}}
\def\br{\begin{remark}}
\def\er{\end{remark}}
\def\bp{\begin{proposition}}
\def\ep{\end{proposition}}
\def\bc{\begin{corollary}}
\def\ec{\end{corollary}}
\def\bd{\begin{definition}}
\def\ed{\end{definition}}
\begin{document}
	
\title{Global Weak Solutions to a Navier--Stokes--Cahn--Hilliard System with Chemotaxis and Mass Transport: \\
Cross Diffusion versus Logistic Degradation}

\author{
Andrea Giorgini
\footnote{Dipartimento di Matematica, Politecnico di Milano, Milano 20133, Italy. E-mail: andrea.giorgini@polimi.it},
\ \
Jingning He
\footnote{Corresponding author. School of Mathematics, Hangzhou Normal University, Hangzhou 311121, P.R. China. Email: hejingning@hznu.edu.cn},
\ \
Hao Wu
\footnote{School of Mathematical Sciences, Fudan University, Shanghai 200433, P.R. China. Email: haowufd@fudan.edu.cn}
}

\date{\today}
\maketitle


\begin{abstract}
\noindent
We analyze a diffuse interface model that describes the dynamics of incompressible two-phase flows influenced by interactions with a soluble chemical substance, encompassing the chemotaxis effect, mass transport, and reactions. In the resulting coupled evolutionary system, the macroscopic fluid velocity field $\bm{v}$ satisfies a Navier--Stokes system driven by a capillary force, the phase field variable $\varphi$ is governed by a convective Cahn--Hilliard equation incorporating a mass source and a singular potential (e.g., the Flory--Huggins type), and the chemical concentration $\sigma$ obeys an advection-reaction-diffusion equation with logistic degradation, exhibiting a cross-diffusion structure akin to the Keller--Segel model for chemotaxis. Under general structural assumptions, we establish the existence of global weak solutions to the initial boundary value problem within a bounded smooth domain $\Omega\subset \mathbb{R}^d$, $d=2,3$. The proof hinges on a novel semi-Galerkin scheme for a suitably regularized system, featuring a non-standard approximation of the singular potential. 
Moreover, with more restrictive assumptions on coefficients and data, we establish regularity properties and uniqueness of global weak solutions in the two-dimensional case.
Our analysis contributes to a further understanding of phase separation processes under the interplay of fluid dynamics and chemotaxis, in particular, the influence of cross diffusion and logistic degradation.
\end{abstract}

\section{Introduction}
\setcounter{equation}{0}
\noindent

The diffuse interface method has attracted significant attention as an efficient and versatile methodology for moving interface problems arising from Materials Science \cite{CH}, fluid dynamics \cite{AMW,LT98}, and mathematical biology, e.g., the tumor growth process \cite{GLSS,OHP}. Within this framework, the classical hypersurface description of free interfaces in the so-called sharp interface models is replaced by a thin layer that permits microscopic mixing of macroscopically distinct components. The diffuse interface approach avoids explicit tracking of free interfaces in both mathematical formulations and numerical computations, and enables handling large deformations and topological changes of the interfaces in a natural way \cite{AMW,DF20}.

In this work, we consider the following Navier--Stokes--Cahn--Hilliard system for viscous incompressible two-phase flows incorporating the effects of chemotaxis, mass transport and reaction:
\begin{subequations}
\label{maineq}
\begin{alignat}{3}
&\partial_t  \bm{ v}+\bm{ v} \cdot \nabla  \bm {v}-\mathrm{div} \big(  2\eta(\varphi) D\bm{v} \big)+\nabla p=(\mu+\chi \sigma)\nabla \varphi,&&\quad \textrm{in}\ \Omega \times (0,T),
\label{main.1} \\
&\mathrm{div}\, \bm{v}=0,&&\quad \textrm{in}\ \Omega \times (0,T),\ \label{main.11}\\
&\partial_t\varphi+\bm{v}\cdot\nabla\varphi =\mathrm{div}\big(m(\varphi)\nabla \mu\big) +S(\varphi,\sigma),&&\quad \textrm{in}\ \Omega \times (0,T),
\label{main.2} \\
&\mu=- \varepsilon \Delta \varphi + \frac{1}{\varepsilon} \Psi'(\varphi) - \chi \sigma, &&\quad \textrm{in}\ \Omega \times (0,T),
\label{main.3}\\
&\partial_t\sigma+\bm{v}\cdot\nabla\sigma-\mathrm{div}\big[ n(\varphi)\sigma\nabla (\ln \sigma +\chi (1-\varphi))\big]
=R(\varphi, \sigma), &&\quad \textrm{in}\ \Omega \times (0,T),
\label{main.4}
\end{alignat}
\end{subequations}
subject to the boundary conditions
\begin{alignat}{3}
&\bm{v}=\mathbf{0},\quad\partial_{\boldsymbol{n}} \varphi=m(\varphi)\nabla   \mu\cdot\bm{n}= \big[n(\varphi)\sigma\nabla (\ln \sigma +\chi (1-\varphi))\big]\cdot\bm{n}=0,\quad \textrm{on}\ \partial\Omega\times(0,T),
\label{boundary}
\end{alignat}
as well as the initial conditions
\begin{alignat}{3}
&\bm{v}|_{t=0}=\bm{v}_{0},\quad \varphi|_{t=0}=\varphi_{0},\quad \sigma|_{t=0}=\sigma_{0}, \qquad &\textrm{in}&\ \Omega.
\label{ini0}
\end{alignat}
Here, we assume that $\Omega \subset\mathbb{R}^d$ (with spatial dimension $d=2,3$) is a bounded domain with smooth boundary $\partial\Omega$, and $T>0$ is a given final time of arbitrary magnitude. The notation $\partial_{\boldsymbol{n}}f$  denotes the normal derivative of a function $f$ on the boundary
with outer unit normal $\bm{n}$. The scalar function $\varphi: \Omega\times (0,T)\to [-1,1]$ is the so-called order parameter (or phase field) that describes the difference in volume fractions of the binary mixture. The level sets $\{\varphi=1\}$ and  $\{\varphi=-1\}$ represent regions occupied by the pure phases of fluid 1 and fluid 2, respectively, while the free interface corresponds to a narrow transition layer of thickness scaling as $\varepsilon\in (0,1)$, in which $\{-1<\varphi<1\}$. The fluid velocity $\bm{v}: \Omega\times (0,T)\to \mathbb{R}^d$ is taken as the volume-averaged velocity with $D\bm{v}=\frac{1}{2}(\nabla\bm{ v}+(\nabla\bm{ v})^\mathrm{T})$ being its symmetrized gradient, and $p:\Omega\times (0,T)\to \mathbb{R}$ is the (modified) pressure. The unknown variable $\sigma:\Omega\times (0,T)\to \mathbb{R}$ denotes the concentration of an unspecific chemical substance and $\mu:\Omega\times (0,T)\to \mathbb{R}$ stands for the chemical potential of the phase separation process.
In the coupled system \eqref{maineq},
equations \eqref{main.1} and \eqref{main.11} represent the momentum balance for the viscous incompressible two-phase flow, equations \eqref{main.2} and \eqref{main.3} constitute a convective Cahn--Hilliard system for the order parameter $\varphi$, and equation \eqref{main.4} yields an advection-diffusion-reaction equation for the chemical concentration $\sigma$.

Throughout this paper, we assume that the density of the binary fluid mixture as well as densities of the individual constituents is constant (all set to one), while the fluid viscosities are allowed to be unmatched. Let $\eta_1$, $\eta_2>0$ be viscosities of the two homogeneous components, respectively. The averaged viscosity is modeled by a concentration dependent function $\eta=\eta(\varphi)$, for instance,
\be
\eta(r)=\eta_1\frac{1+r}{2}+\eta_2\frac{1-r}{2}, \quad \forall\, r\in[-1,1].
\label{vis}
\ee
In \eqref{main.3}, the nonlinear function $\Psi'$ denotes the derivative of a potential $\Psi$ that has a double-well structure, with two minima and a local unstable maximum in between. The physically significant example is the logarithmic type (also called the Flory--Huggins potential):
 \be
\Psi (r)=\frac{\theta}{2}[(1-r)\ln(1-r)+(1+r)\ln(1+r)]+\frac{\theta_{c}}{2}(1-r^2),\quad \forall\, r\in[-1,1],
\label{pot}
 \ee
 with $0<\theta<\theta_{c}$ (see, e.g., \cite{CH,Mi19}). It is referred to as a singular potential since its derivative $\Psi'$ blows up at the pure phases $\pm 1$. In the literature, the singular potential $\Psi$ is often approximated by a fourth-order polynomial
  \be
\Psi(r)=\frac{1}{4}(1-r^2)^2,\quad r\in\mathbb{R}, \label{regular}
 \ee
or some more general polynomial functions \cite{Mi19}. The functions $m(\varphi)$ and $n(\varphi)$ in \eqref{main.2} and \eqref{main.4} are nonnegative mobility functions related to the phase field and the chemical concentration, respectively.

The coupled system \eqref{maineq} is thermodynamical consistent due to its variational structure, that is, solutions to \eqref{maineq} (under the boundary conditions \eqref{boundary}) (formally) satisfy the basic energy law:
\begin{align}
& \frac{\mathrm{d}}{\mathrm{d}t} \int_{\Omega} \Big( \frac{1}{2}|\bm{v}|^2 +\frac{\varepsilon}{2}|\nabla \varphi|^2
+ \frac{1}{\varepsilon}\Psi(\varphi)
+\sigma(\ln\sigma-1 )+\chi\sigma(1-\varphi) \Big)\, \mathrm{d}x
\notag \\
& \qquad +\int_{\Omega} \Big( 2\eta(\varphi)|D\bm{v}|^2 + m(\varphi)|\nabla \mu|^2+ n(\varphi)|\sigma^{\frac12}\nabla(\sigma+\chi(1-\varphi))|^2\Big)\, \mathrm{d}x
\nonumber\\
&\quad  =\int_\Omega \big[S(\varphi, \sigma)\mu+R(\varphi, \sigma)(\ln \sigma+\chi(1-\varphi))\big] \,\mathrm{d}x.
\label{BEL}
\end{align}
As seen in \eqref{BEL}, the total energy density is the sum of the kinetic energy density $\frac12 |\bm{v}|^2$ for the macroscopic fluid, the Ginzburg--Landau energy density $\frac{\varepsilon}{2}|\nabla \varphi|^2
+  \frac{1}{\varepsilon} \Psi(\varphi)$ for mixing of the binary mixture, and the chemical free energy density given by
\be
 N(\varphi, \sigma)=\sigma(\ln \sigma-1)+\chi \sigma(1-\varphi).\label{energy density}
\ee
In the context of tumor growth modeling (see, e.g., \cite{GLSS, GL17e}), the constant coefficient $\chi$ is related to certain key transport mechanisms, such as chemotaxis and active transport (when $\chi>0$). The source term $S(\varphi, \sigma)$ in \eqref{main.2} corresponds to biological mechanisms like proliferation, apoptosis of cells (see \cite{Mi19} for further discussions on biologically relevant mass source terms). Besides, the function $R(\varphi,\sigma)$ in \eqref{main.4} models certain reaction for the chemical substance. The presence of mass source/reaction terms not only affects the mass dynamics for the mixture and the chemical substance, but also contributes to changes in energy, such that the total energy of \eqref{maineq} may increase with over time.

Our system \eqref{maineq} contains the well-known Navier--Stokes-Cahn--Hilliard system as a subsystem
(i.e., the ``Model H'' \cite{HH, Gur}), which has been extensively analyzed in the literature, see \cite{A2009,B,GG2010,GMT,MT,ZWH} and the references cited therein.
In the Navier--Stokes part \eqref{main.1}--\eqref{main.11},  the coupling structure is reflected in terms of the capillary force $(\mu+\chi \sigma)\nabla \varphi$ (only depending on $\varphi$ in view of \eqref{main.3}) and the viscous stress tensor with a concentration dependent viscosity $\eta(\varphi)$. Besides, the macroscopic velocity $\bm{v}$ influences the dynamics of $\varphi$ and $\sigma$ via
the advection terms $\bm{v}\cdot \nabla \varphi$ and $\bm{v}\cdot \nabla \sigma$. Concerning the interaction between the phase field $\varphi$ and the chemical concentration $\sigma$, we reformulate \eqref{main.2} and \eqref{main.4} as follows
\begin{align*}
&\partial_t \varphi+ \bm{v} \cdot \nabla \varphi +\mathrm{div}\,\bm{q}_\varphi=S(\varphi,\sigma),\\
& \partial_t \sigma+ \bm{v} \cdot \nabla \sigma + \mathrm{div}\,\bm{q}_\sigma=R(\varphi, \sigma),
\end{align*}
 with mass fluxes
\begin{align*}
&\bm{q}_\varphi:= -m(\varphi)\nabla \mu= -m(\varphi)\nabla( - \varepsilon\Delta \varphi + \varepsilon^{-1}\Psi'(\varphi)+N_{\varphi}),
\quad
\bm{q}_\sigma:=-n(\varphi)\sigma\nabla N_{\sigma},
\end{align*}
where $N_{\varphi}=-\chi\sigma $ and $N_{\sigma}=\ln \sigma +\chi(1-\varphi)$ denote the variational derivatives of $N(\varphi, \sigma)$ with respect to $\varphi$ and $\sigma$, respectively.
The term $\chi \nabla \sigma$ in $\bm{q}_\varphi$ represents the chemotactic response of the mixture to the chemical substance, while the other term $\chi\sigma\nabla (1-\varphi)$ in $\bm{q}_\sigma$ propels the chemical
substance via the concentration gradient of the mixture, which differs from the conventional diffusion mechanism.

In comparison with the general thermodynamically consistent diffuse interface model derived in \cite{LW} for a mixture of two viscous incompressible fluids interacting with a chemical substance (see also \cite{Sitka} in the context of tumor growth modeling), a distinct feature of equation \eqref{main.4} is the nonlinear cross-diffusion term $\chi\mathrm{div}(n(\varphi)\sigma \nabla\varphi)$. In \cite{LW} and related works for tumor growth, this is usually replaced by a simpler one with linear mass transport (see e.g., \cite{EG19jde,GL17e,GLSS}):
\be
\partial_t\sigma-\Delta\sigma+ \chi \Delta \varphi =R(\varphi, \sigma),
\label{diffu}
\ee
where for simplicity, we have taken a constant mobility coefficient $n(\varphi)=1$ and neglected the advection term $\bm{v} \cdot \nabla \sigma$. Since $\sigma$ plays the role of a concentration, its nonnegativity is an expected property. However, the term $\chi\Delta\varphi$ in \eqref{diffu} has no definite sign, preventing the application of the  minimum principle to guarantee that $\sigma$ will stay non-negative during its evolution. This unsatisfactory issue was avoided in \cite{EL21} by considering a quasi-steady-state equation $\Delta\sigma-\sigma h(\varphi)=0$ for $\sigma$. On the other hand, inspired by the Keller--Segel system for chemotaxis (see e.g., \cite{BBTW,KS}):
\begin{equation}
\left\{ \begin{array}{l}
\partial_t u=\mathrm{div}(\gamma (v)\nabla u-u\chi (v)\nabla v),\\
\tau \partial_t v=\Delta v + u -v,
\end{array} \right.
\label{KS-a}
\end{equation}
the authors of \cite{RSS} proposed an alternative of \eqref{diffu} such that
\be
\partial _t\sigma - \Delta \sigma + \chi  \mathrm{div} (\sigma\nabla \varphi) = R(\varphi,\sigma).
\label{signew}
\ee
With a properly chosen reaction term $R(\varphi,\sigma)$, they were able to show that $\sigma$ satisfies the required non-negativity preserving property, provided that its initial value $\sigma_0$ is non-negative. See also \cite{AS24} for a similar consideration in a multi-species tumor growth model with chemotaxis and angiogenesis.
From a modeling perspective, equation \eqref{signew} is more natural for the mass transfer process of tumor growth, since the mass flux  $\chi\sigma\nabla \varphi$ depends on both the chemical substance concentration $\sigma$ and the gradient of mixture distribution $\nabla \varphi$. On the other hand, a new mathematical structure emerges in \eqref{maineq} with the choice of \eqref{signew}: the equation \eqref{main.4} for the chemical concentration $\sigma$ is similar to the equation for the cell density $u$ in the Keller--Segel system, while the fourth-order Cahn--Hilliard equation \eqref{main.2}--\eqref{main.3} for the order parameter $\varphi$ replaces the linear, second-order reaction-diffusion equation for the chemical signal concentration $v$ in \eqref{KS-a}. This requires a rather different treatment for $\varphi$ comparing with that for $v$, which introduces challenges in the mathematical analysis.

To the best of our knowledge, there is no analytic result concerning the initial boundary value problem \eqref{maineq}--\eqref{ini0} in the literature. A fluid-free version has been analyzed in \cite{RSS}, that is, the Navier--Stokes equations for the fluid velocity $\bm{v}$ were neglected.
Assuming a logistic type reaction term like $R(\varphi, \sigma)=\beta(\varphi)\big(\kappa_0\sigma-\kappa_{\infty}\sigma^p\big)$ for some $p\in(1,2]$, the authors proved the existence of global weak solutions in two and three dimensional cases when $\chi>0$. Further regularity properties were established under more restrictive assumptions on structural coefficients and data. The proof therein is mainly based on \emph{a priori} estimates and compactness methods. At the level of weak solutions, one crucial issue is about the coercivity of the energy functional (more precisely, \eqref{energy density}), which can be handled with the uniform boundedness of $\varphi$, hence on the choice of a singular configuration potential like \eqref{pot}. The authors of \cite{RSS} outlined a possible regularization scheme in which both the singular potential $\Psi$ and the chemical concentration $\sigma$ are properly truncated. They were able to control the crossing term $\chi \sigma(1-\varphi)$ in \eqref{energy density} and obtained uniform estimates for approximate solutions, which allowed them to pass to the limit. In particular, the logistic growth with respect to $\sigma$ helped the derivation of adequate \emph{a priori} estimates. We also mention the recent work \cite{Sch24} on tumor growth processes influenced by chemotaxis and mass transport, in which a Cahn--Hilliard--Brinkman system with a singular potential was analyzed. Using a different approximating scheme, the author demonstrated the existence of global weak solutions in two and three dimensions, provided that the chemotactic sensitivity function exhibits a \emph{nonlinear} dependence on the chemical concentration like $\frac{\sigma}{1+\sigma^{q-1}}$ for some $q\in (1,2]$ (the so-called ``degenerate sensitivity'' characterized by a slower growth for large values of $\sigma$), but without incorporating the logistic degradation.

Our aim is to establish the existence of global weak solutions to the initial boundary value problem \eqref{maineq}--\eqref{ini0} in both two and three dimensions, under the choice of a singular potential like \eqref{pot} and a reaction term $R(\varphi,\sigma)$ involving logistic degradation. As shown in Theorem \ref{main}, problem \eqref{maineq}--\eqref{ini0} admits a global weak solution $(\bm{v},\varphi, \sigma)$ on $[0,T]$, provided that the initial total energy is finite, that is, $ \bm{v}_{0} \in \bm{L}^2(\Omega)$ with $\mathrm{div}\,\bm{v}_0=0$, $\varphi_{0}\in H^1(\Omega)$ with $\|  \varphi_{0} \|_{L^{\infty}(\Omega)} \le 1$, $|\overline{\varphi}_{0}|<1$, and $\sigma_{0} \ln \sigma_{0} \in L^{1}(\Omega)$ with $\sigma_{0} \geq 0$ almost everywhere in $\Omega$. 
When the spatial dimension is two, further conclusions can be achieved: the global weak solution satisfies better regularity properties under a more regular initial datum $\sigma_0\in L^2(\Omega)$ (see Theorem \ref{main}); moreover, assuming constant mobility functions and a constant mass source, the aforementioned global weak solution is unique (see Theorem \ref{main-2}).

Let us highlight the main novelties of our study.
Firstly, we have proposed a new approximation of the singular potential $\Psi$. The existence of global weak solutions to problem \eqref{maineq}--\eqref{ini0} is primarily based on the energy balance \eqref{BEL}. To handle the Cahn--Hilliard equation \eqref{main.2}--\eqref{main.3} with a singular potential like \eqref{pot}, a conventional strategy is to smooth out the singular term $\Psi(\varphi)$ and subsequently apply some local existence theorem, for instance, an appropriate Faedo--Galerkin scheme. This can be accomplished through a polynomial-type regularization (see e.g., \cite{MT,GGW}), or a Yosida approximation (see e.g., \cite{B2010}). Nonetheless, the uniform boundedness of the approximate solution for $\varphi$ cannot be guaranteed, which is essential to manage the crossing term $\chi \sigma(1-\varphi)$. Furthermore, neither the growth of the polynomial regularization nor the Yosida approximation is sufficient to handle the product $\sigma\varphi$ in conjunction with the mixing entropy $\sigma\ln \sigma$. In light of the generalized Young's inequality (see Lemma \ref{You}), we find that the aforementioned crossing term can be controlled by the mixing entropy for $\sigma$ combined with a proper exponential function of $\varphi$ (depending on the strength of chemotactic effect characterized by $|\chi|$). This insight motives us to introduce a new approximate potential $\Psi_\epsilon(\varphi)$ that is regular on $\mathbb{R}$ and exhibits exponential growth for large values of $\varphi$.
This enables us to maintain the uniform coercivity of the free energy at the level of approximate solutions, without involving any truncation for $\sigma$ as was done in \cite{RSS}. We note that the degenerate sensitivity considered in \cite{Sch24} leads to an additional contribution $\frac{1}{q(q-1)}\sigma^q$ in the chemical free energy density $N(\varphi, \sigma)$. This gives a stronger coercivity on $\sigma$ so that a sufficiently higher order polynomial-type regularization of $\Psi(\varphi)$ can be adopted.

Secondly, we have devised a new semi-Galerkin approximation scheme that incorporates a $p$-Laplace regularization. Due to the nonlinear interaction with the Navier--Stokes system, it is unclear whether the approximating scheme outlined in \cite{RSS} can be extended to our system \eqref{maineq}. Besides, it remains uncertain how to apply the approximation schemes presented in \cite{Lan2016,W2016,W2019} for the Keller--Segel--Navier--Stokes system with or without a logistic source, owing to the current coupling between the Cahn--Hilliard and Keller--Segel equations as well as the singular character of the potential $\Psi$.
A feasible approach for the Navier--Stokes--Cahn--Hilliard system with a singular potential and chemotactic effects is the so-called semi-Galerkin scheme. For instance, in \cite{H}, the author proved the existence of global weak solutions to a system analogous to \eqref{maineq} in a simpler scenario with linear mass transport akin to \eqref{diffu}. Roughly speaking, the author performed a Galerkin approximation solely for the Navier--Stokes system, but solved the coupled system for $(\varphi, \sigma)$ independently. This approach ensures that the approximate solution for $\varphi$ takes its values in $[-1,1]$ and facilitates the completion of a fixed point argument.
Unfortunately, the semi-Galerkin scheme used in \cite{H} is not applicable here, due to the cross-diffusion structure inherent in \eqref{main.4} (cf. \cite{RSS}). Indeed, it yields weaker energy for the chemical substance (that is, $\sigma(\ln \sigma -1)$ versus $\frac{1}{2}\sigma^2$ as in \cite{GLSS,LW}) and weaker energy dissipation in the basic energy law \eqref{BEL} (that is, $n(\varphi)|\sigma^{\frac12}\nabla(\sigma+\chi(1-\varphi))|^2$ versus
$n(\varphi)|\nabla(\sigma+\chi(1-\varphi))|^2$ as in \cite{GLSS,LW}).
To overcome these difficulties, we propose a new semi-Galerkin scheme. This involves performing a Faedo--Galerkin approximation to the Navier--Stokes--Cahn--Hilliard system \eqref{main.1}--\eqref{main.3} with the aforementioned regularization of the singular potential $\Psi$, while solving the advection-diffusion-reaction equation \eqref{main.4} for $\sigma$ independently. In particular, the approximate solution of $\sigma$ is solved in the classical sense such that it is strictly positive in $\Omega\times(0,T)$ thanks to the strong maximum principle. This key property allows us to use $\ln \sigma$ as a test function in the derivation of uniform estimates for approximate solutions.

It is worth mentioning that a (regular) potential function with exponential growth for the Cahn--Hilliard equation is critical in two dimensions and super-critical in three dimensions in view of the Sobolev embedding theorem for $H^1(\Omega)$. To overcome this difficulty and the lack of maximum principle, we include an additional $p$-Laplace regularization (with $p=4$ ) in the chemical potential $\mu$, which yields an $L^\infty$-bound for the approximate solution of $\varphi$ (not necessarily in the physical interval $[-1,1]$), since $W^{1,4}(\Omega)\hookrightarrow L^\infty(\Omega)$ for $d=2,3$. This method has been applied in \cite{BG24} within a periodic setting, where a diffuse interface model describing the complex rheology and interfacial dynamics during phase separation in a polar liquid-crystalline emulsion was analyzed. Taking advantage of recent works \cite{CM14,CM19} on the $p$-Laplace system, we successfully extend this regularization technique to our problem \eqref{maineq}--\eqref{ini0} in a general bounded smooth domain, which may have its independent interest. Comparing with the regularization via a bi-Laplacian as described in \cite{Sch24}, the $p$-Laplace regularization appears to be more manageable.

In this study, we illustrate the approximating scheme in full detail and rigorously justify the associated uniform estimates. The presented method has the potential for application to scenarios when the Navier--Stokes system \eqref{main.1}--\eqref{main.11} for the fluid velocity is substituted with either a Brinkman's system \cite{EG19jde,CG}, or a Darcy's system \cite{GLSS,GGW}. It is well-known that logistic degradation exerts a regularizing effect on chemotaxis models \cite{W2010,TW,L2021}, and this phenomenon has also been observed in chemotaxis-fluid systems \cite{TW15,W2019,Lan2016}. Our analysis explicitly indicates the role of the logistic term in the highly nonlinear system \eqref{maineq} within the framework of global weak solutions. It will be interesting to investigate the problem of global regularity and boundedness of weak solutions as well as their long-time behavior. These issues will be addressed in future works.

The remaining part of this paper is organized as follows. In Section \ref{pm}, we introduce the functional settings and state the main results (i.e., Theorems \ref{main}, \ref{main-2}).
In Section \ref{ws}, we propose a semi-Galerkin scheme for a suitably regularized system, and then demonstrate its solvability. Section \ref{proof-main} is devoted to the proof of Theorem \ref{main}. We first derive uniform estimates for approximate solutions and then construct global weak solutions by weak compactness methods. In Section \ref{sec:uni}, we prove Theorem \ref{main-2} on the uniqueness of global weak solutions in the two-dimensional case. In the Appendix, we provide some details of the analysis for the semi-Galerkin scheme.

\section{Main Results}
\setcounter{equation}{0}
\label{pm}

\subsection{Preliminaries}

Let $X$ be a real Banach space. We denote its norm by $\|\cdot\|_X$, its dual space by $X'$, and the duality pairing by $\langle \cdot,\cdot\rangle_{X',X}$. The bold letter $\bm{X}$ denotes the generic space of vectors or matrices, with each component belonging to $X$. Given a measurable set $J\subset \mathbb{R}$, $L^q(J;X)$ with $q\in [1,+\infty]$ denotes the space of Bochner measurable $q$-integrable/essentially bounded functions with values in the Banach space $X$. If $J=(a, b)$ is an interval, we write for simplicity $L^q(a,b;X)$. For $q\in [1,+\infty]$, $W^{1,q}(J;X)$ denotes the space of functions $f$ such that $f\in L^q(J;X)$ with $\partial_t f\in L^q(J;X)$, where $\partial_t $ means the vector-valued distributional derivative of $f$. When $q=2$, we set $H^1(J;X):=W^{1,2}(J;X)$. Besides, we denote by $C(J;X)$ (or $C_w(J;X)$)  the space of functions that are strong (or weak) continuous from $J$ to $X$.

We assume that $\Omega \subset\mathbb{R}^d$ ($d=2,3$) is a bounded domain with sufficiently smooth boundary $\partial\Omega$ such that $\partial\Omega$ is a $(d-1)$-dimensional compact and connected hypersurface without boundary.
For any $q \in [1,+\infty]$, $L^{q}(\Omega)$ and $L^{q}(\partial\Omega)$  denote the Lebesgue spaces on $\Omega$ and $\partial\Omega$, respectively.
For $s\geq 0$ and $q\in [1,\infty)$,
we denote by $H^{s,q}(\Omega )$ the Bessel-potential spaces and by $%
W^{s,q}(\Omega )$ the Slobodeckij spaces. It holds that $H^{s,2}(\Omega
)=W^{s,2}(\Omega )$ for all $s$, but for $q\neq 2$ the identity $%
H^{s,q}(\Omega )=W^{s,q}(\Omega )$ is only true if $s\in \mathbb{N}$.
For $s\in \mathbb{N}$, $H^{s,q}(\Omega )$ and $W^{s,q}(\Omega )$
coincide with the usual Sobolev spaces. The corresponding function spaces
over the boundary $\partial\Omega$ are defined via local charts.
If $q=2$ and $s\in \mathbb{Z}^+$, we shall use the standard notation $%
H^{s}(\Omega ):=H^{s,2}(\Omega )=W^{s,2}(\Omega )$ for functions defined in $\Omega$, and $H^{s}(\partial\Omega ):=H^{s,2}(\partial\Omega )=W^{s,2}(\partial\Omega )$ for functions defined on $\partial \Omega$.
For simplicity, the norm and inner product in the basic space $L^{2}(\Omega)$ (as well as $\bm{L}^{2}(\Omega)$) are denoted by $\|\cdot\|$ and $(\cdot,\cdot)$, respectively.

For every $f\in (H^1(\Omega))'$, we define its generalized mean over $\Omega$ by
$\overline{f}=\frac{1}{|\Omega|} \langle f,1\rangle_{(H^1(\Omega))',H^1(\Omega)}$; if $f\in L^1(\Omega)$, then it holds $\overline{f}= \frac{1}{|\Omega|} \int_\Omega f \,\mathrm{d}x$.
Define the linear subspaces
\begin{align*}
&L^2_{0}(\Omega):=\big\{f\in L^2(\Omega)\ |\ \overline{f} =0\big\},
\qquad V_0:= H^1(\Omega)\cap L_0^2(\Omega),
\\
& V_0^{-1}:= \big\{ f \in (H^1(\Omega))'\ |\  \overline{f}=0 \big\}\subset (V_0)'.
\end{align*}
In view of the homogeneous Neumann boundary condition, we also introduce the space
$$
H^2_{N}(\Omega):= \big\{f\in H^2(\Omega)\ |\  \partial_{\bm{n}}f=0 \ \textrm{on}\  \partial \Omega\big\}.
$$
Let $\mathcal{A}_N\in \mathcal{L}(H^1(\Omega),(H^1(\Omega))')$ be the realization of minus Laplacian $-\Delta$ subject to the homogeneous Neumann boundary condition such that
\begin{equation}\nonumber
	\langle \mathcal{A}_N u,v\rangle_{(H^1)',H^1} := \int_\Omega \nabla u\cdot \nabla v \, \mathrm{d}x,\quad \forall\, u,v\in H^1(\Omega).
\end{equation}
The restriction of $\mathcal{A}_N$ from $V_0$ onto $V_0^{-1}$ is an isomorphism. Besides, $\mathcal{A}_N$ is positively defined on $V_0$ and self-adjoint. We denote the inverse map by $\mathcal{N} =\mathcal{A}_N^{-1}: V_0^{-1} \to V_0$. For every $f\in V_0^{-1}$, $u= \mathcal{N} f \in V_0$ is the unique weak solution of the Neumann problem
$$
\begin{cases}
	-\Delta u=f, \quad \text{in} \ \Omega,\\
	\partial_{\bm{n}} u=0, \quad \ \  \text{on}\ \partial \Omega.
\end{cases}
$$
For every $f\in V_0^{-1}$, we set $\|f\|_{V_0^{-1}}=\|\nabla \mathcal{N} f\|$.
It is well-known that $f \to \|f\|_{V_0^{-1}}$ and $
f \to\big(\|f-\overline{f}\|_{V_0^{-1}}^2+|\overline{f}|^2\big)^\frac12$
are norms on $V_0^{-1}$ and $(H^1(\Omega))'$,
respectively, which are equivalent to the standard ones (see e.g., \cite{MZ04}). From the Poincar\'{e}--Wirtinger inequality:
\begin{equation}
	\notag
	\|f-\overline{f}\|\leq C \|\nabla f\|,\qquad \forall\,
	f\in H^1(\Omega),
\end{equation}
where $C>0$ only depends on $\Omega$, we find that $f\to \|\nabla f\|$ and  $f\to \big(\|\nabla f\|^2+|\overline{f}|^2\big)^\frac12$  are norms on $V_0$ and $H^1(\Omega)$, respectively, which are equivalent to the standard ones.
Moreover, we report the following standard Hilbert interpolation inequality
\begin{align*}
	&\|f\|  \leq \|f\|_{V_0^{-1}}^{\frac12} \| \nabla f\|^{\frac12},
	 \qquad  \forall\, f \in V_0.
\end{align*}
For later use, we also consider the elliptic operator $\mathcal{A}_1 := -\Delta+I $ subject to the  homogeneous Neumann boundary condition, which is an unbounded operator in $L^2(\Omega)$ with domain $D(\mathcal{A}_1) =H^2_N(\Omega)$. It is well-known that
$\mathcal{A}_1$ is a positive and self-adjoint operator in $L^2(\Omega)$ with a compact inverse denoted by $\mathcal{N}_1:=\mathcal{A}^{-1}_1$. The spectral theory allows us to define the powers $\mathcal{A}_1^s$ for $s\in \mathbb{R}$.

Let us now introduce the Hilbert spaces of solenoidal vector-valued functions. As in \cite{S}, we denote by $\bm{L}^2_{0,\sigma}(\Omega)$, $\bm{H}^1_{0,\sigma}(\Omega) $ the closure of   $C_{0,\sigma}^{\infty}(\Omega;\mathbb{R}^d)=\big\{\bm{f}\in C_0^{\infty}(\Omega;\mathbb{R}^d):\, \mathrm{div}\bm{f}=0\big\}$ in  $\bm{L}^2(\Omega)$ and $\bm{H}^1(\Omega)$,
respectively\footnote{The subscript $\sigma$ is a conventional notation for spaces of divergence-free functions. It should not be related to the solution $\sigma$.}.
For simplicity, we use $(\cdot,\cdot)$ and $\|\cdot\|$ for the inner product and norm in $\bm{L}^2_{0,\sigma}(\Omega)$. For any function $\bm{f} \in \bm{L}^2(\Omega)$, the Helmholtz--Weyl decomposition holds (see  \cite[Chapter \uppercase\expandafter{\romannumeral3}]{G}):
\be
\bm{f}=\bm{f}_{0}+\nabla z,\quad\text{where}\ \ \bm{f}_{0} \in \bm{L}^2_{0,\sigma}(\Omega),\ z \in H^1(\Omega).\nonumber
\ee
Define the Leray projection $\bm{P}:\bm{L}^2(\Omega)\to \bm{L}^2_{0,\sigma}(\Omega)$ such that $\bm{P}(\bm{f})=\bm{f}_{0}$.
It holds $\|\bm{P}(\bm{f})\|\leq \|\bm{f}\|$ for all $\bm{f}\in \bm{L}^2(\Omega)$.
The space $\bm{H}^1_{0,\sigma}(\Omega)$ is equipped with the inner product $(\bm{u},\bm{v})_{\bm{H}^1_{0,\sigma}}:=(\nabla \bm{u},\nabla \bm{v})$ and the norm $\|\bm{u}\|_{\bm{H}^1_{0,\sigma}}=\|\nabla \bm{u}\|$.
Owing to Korn's inequality
$$
\|\nabla \bm{u}\|\leq \sqrt{2}\|D\bm{u}\|
\leq \sqrt{2}\|\nabla \bm{u}\|,\quad \forall\, \bm{u}\in \bm{H}^1_{0,\sigma}(\Omega),
$$
$\|D\bm{u}\|$ is an equivalent norm for $\bm{H}^1_{0,\sigma}(\Omega)$. Next, we introduce the Stokes operator $\bm{S}=\bm{P}(-\Delta)$ with domain $D(\bm{S})= \bm{H}^1_{0,\sigma}(\Omega)\cap\bm{H}^2(\Omega)$.
The Hilbert space $D(\bm{S})$ is equipped with the inner product $(\bm{S}\bm{u},\bm{S}\bm{v})$ and the norm $\|\bm{S}\bm{u}\|$  (see e.g., \cite[Chapter III]{S}). For any $\bm{u}\in D(\bm{S})$ and $\bm{\zeta} \in \bm{H}^1_{0,\sigma}(\Omega)$, it holds
$(\bm{S}\bm{u},\bm{\zeta})=(\nabla \bm{u},\nabla\bm{\zeta})$.
The operator $\bm{S}$ is a canonical isomorphism from $\bm{H}^1_{0,\sigma}(\Omega)$ to $(\bm{H}^1_{0,\sigma}(\Omega))'$ and we denote its inverse by $\bm{S}^{-1}$. For any $\bm{f}\in (\bm{H}^1_{0,\sigma}(\Omega))'$, there is a unique $\bm{u}=\bm{S}^{-1}\bm{f}\in\bm{H}^1_{0,\sigma}(\Omega)$ such that
$(\nabla\bm{S}^{-1}\bm{f},\nabla \bm{\zeta}) =\langle\bm{f},\bm{\zeta}\rangle_{(\bm{H}^1_{0,\sigma})', \bm{H}^1_{0,\sigma}}$  for all $\bm{\zeta} \in \bm{H}^1_{0,\sigma}(\Omega)$.
Hence, $\|\nabla\bm{S}^{-1}\bm{f}\| =\langle\bm{f},\bm{S}^{-1}\bm{f}\rangle_{(\bm{H}^1_{0,\sigma})',\bm{H}^1_{0,\sigma}}^{\frac{1}{2}}$ yields an equivalent norm on the dual space $(\bm{H}^1_{0,\sigma}(\Omega))'$. 


Finally, we recall the following generalized Young's inequality (see, e.g., \cite{MZ04}):
\begin{lemma}\label{You}
Let
\be
f(a):=\mathrm{e}^a-a-1, \quad g(b):=(b+1) \ln (b+1)-b.
\label{general young}
\ee
Then it holds
$$
a b \le f(a)+g(b), \quad \forall\, a, b \ge 0.
$$
\end{lemma}

In the subsequent sections, the symbols $C$, $C_i$ stand for generic positive constants that may even change within the same line. Specific dependence of these constants in terms of the data will be pointed out if necessary.

\subsection{Statement of main result}
\noindent
We make the following hypotheses, which will be kept for the remainder of this paper.
\begin{itemize}
\item[(H1)] The singular potential $\Psi$ belongs to the class of functions $C([-1,1])\cap C^{2}(-1,1)$. It can be written into the following form
\begin{equation}
\Psi(r)=\Psi_{0}(r)-\frac{\theta_{0}}{2}r^2,\nonumber
\end{equation}
such that
\begin{equation}
\lim_{r\to \pm 1} \Psi_{0}'(r)=\pm \infty ,\quad \text{and}\ \  \Psi_{0}''(r)\ge \theta,\quad \forall\, r\in (-1,1),\nonumber
\end{equation}
where $\theta$ is a strictly positive constant and $\theta_0\in \mathbb{R}$. In addition, there exists $\epsilon_0\in(0,1)$ such that $\Psi_{0}''$ is nondecreasing in $[1-\epsilon_0,1)$ and nonincreasing in $(-1,-1+\epsilon_0]$. We make the extension $\Psi_{0}(r)=+\infty$ for any $r\notin[-1,1]$. Without loss of generality, we also set $\Psi_0(0)=\Psi_0'(0)=0$.
\item[(H2)] The viscosity function $\eta\in C^{1}(\mathbb{R})$ is globally Lipschitz continuous in $\mathbb{R}$ and
\be
\eta_{*} \leq \eta(r)\leq \eta^*,\quad \forall\, r \in \mathbb{R},\nonumber
\ee
where $\eta_{*}<\eta^*$ are given positive constants.
\item[(H3)] The mobility function $m\in C^{1}(\mathbb{R})$ is globally Lipschitz continuous in $\mathbb{R}$ and
\be
m_{*} \leq m(r)\leq m^*,\quad \forall\, r \in \mathbb{R},\nonumber
\ee
where $m_{*}<m^*$ are given positive constants. Moreover, we set the mobility function $n(r)\equiv 1$ for all $r\in \mathbb{R}$.
\item[(H4)] The mass source term $S$ satisfies
$$
S(\varphi, \sigma)= -\alpha \varphi+ h(\varphi, \sigma),\quad \forall\, (\varphi,\sigma)\in \mathbb{R}\times \mathbb{R},
$$
where $\alpha$ is a given positive constant. The function $h$ is  uniformly bounded and global Lipschitz continuous with respect to its variables. Moreover, the following compatibility condition holds
$$
\alpha > \|h\|_{L^\infty(\mathbb{R}\times \mathbb{R})} =: h^*\geq 0.
$$
\item[(H5)] The reaction term $R$ satisfies
$$
R(\varphi, \sigma)= \beta(\varphi)\sigma -\kappa  \sigma^2, \quad \forall\, (\varphi,\sigma)\in \mathbb{R}\times \mathbb{R},
$$
where $\kappa$ is a given positive constant.
The function $\beta \in C^{1}(\mathbb{R})$ satisfies
\be
\begin{aligned}
&|\beta(r)| \leq b^*, &&\quad \forall\, r \in \mathbb{R},
\notag\\
&\beta(r)=0,&&\quad \forall\, r\in (-\infty,-2]\cup [2, +\infty),
\end{aligned}
\ee
where $b^*$ is a given positive constant.
%
\item[(H6)] The coefficients $\varepsilon$,  $\chi$ are prescribed constants such that
\be
\varepsilon=1,\quad \chi \in \mathbb{R}.
\nonumber
\ee
\end{itemize}
\begin{remark}
The physically relevant logarithmic potential \eqref{pot} fulfills the assumption (H1). Besides, as noted in \cite[Remark 2.1]{GMT}, one can easily extend the linear viscosity function \eqref{vis} to $\mathbb{R}$ in such a way to comply (H2). Since the singular potential guarantees that $\varphi\in [-1,1]$, the value of $\eta$ outside of $[-1,1]$ is not important and can be appropriately chosen as in (H2). Similarly, the only significant values of $m$, $n$, $\beta$ will be those for $r\in[-1,1]$.
We extend them outside $[-1,1]$ due to the Galerkin approximation for $\varphi$. To focus on the difficulty from the cross-diffusion structure and  the regularizing effect due to logistic degradation,  we keep a simpler form of the $\sigma$-equation by assuming $n(\varphi)\equiv 1$ and $\kappa$ a positive constant.
\end{remark}

\begin{remark}
Concerning the mass source term $S$, only the behavior of $h$ in the physical reference set $[-1,1]\times[0, +\infty)$ is significant. If $h$ is a constant belonging to $(-\alpha,\alpha)$, $S$ corresponds to the classical Oono's interaction \cite{GG, MT}. The case $S\equiv 0$ (i.e., $\alpha=0$ and $h\equiv 0$) yields the conservation of mass, is admissible too, and in fact simpler to handle.
The magnitude of the chemotaxis sensitivity $\chi$ will play a role in part of the results. The case $\chi=0$ is somewhat trivial since the effects of chemotaxis and mass transport are switched off. Comparing with \cite{RSS}, in which $\chi$ was assumed to be a strictly positive constant, here we are able to deal with any $\chi\in \mathbb{R}$. Since we do not consider the so-called sharp interface limit as $\varepsilon\to 0^+$, the magnitude of $\varepsilon$ has no influence on the subsequent analysis. For the sake of simplicity, we set $\varepsilon=1$.
\end{remark}

The assumptions (H1)--(H6) allow us to rewrite problem \eqref{maineq}--\eqref{ini0} in the following form:
\begin{subequations}
\label{sysNew}
\begin{alignat}{3}
&\partial_t  \bm{ v}+\bm{ v} \cdot \nabla  \bm {v}-\mathrm{div} \big(  2\eta(\varphi) D\bm{v} \big)+\nabla p=(\mu+\chi \sigma)\nabla \varphi,
\label{main.1new} \\
&\mathrm{div}\, \bm{v}=0,
\label{main.11new}\\
&\partial_t\varphi+\bm{v}\cdot\nabla\varphi =\mathrm{div}\big(m(\varphi)\nabla \mu\big) -\alpha \varphi+ h(\varphi, \sigma),
\label{main.2new} \\
&\mu=- \Delta \varphi + \Psi'(\varphi) - \chi \sigma,
\label{main.3new}\\
&\partial_t\sigma+\bm{v}\cdot\nabla\sigma -\Delta \sigma + \chi \mathrm{div}
(\sigma\nabla  \varphi)
=\beta(\varphi) \sigma-\kappa \sigma^2,
\label{main.4new}
\end{alignat}
\end{subequations}
in $\Omega \times (0,T)$, subject to the boundary conditions
\begin{alignat}{3}
&\bm{v}=\mathbf{0},\quad\partial_{\boldsymbol{n}} \varphi = \partial_{\boldsymbol{n}} \mu = \partial_{\boldsymbol{n}}\sigma =0,\quad \textrm{on}\ \partial\Omega\times(0,T),
\label{boundarynew}
\end{alignat}
and the initial conditions
\begin{alignat}{3}
&\bm{v}|_{t=0}=\bm{v}_{0},\quad \varphi|_{t=0}=\varphi_{0},\quad   \sigma|_{t=0}=\sigma_{0}, \qquad &\textrm{in}&\ \Omega.
\label{ini0new}
\end{alignat}

Next, we introduce the definition of weak solutions.
\bd \label{maind}
Let $d =2,3$. Suppose that (H1)--(H6) are satisfied and $T \in (0,+\infty)$ is a given final time. For any given initial data $(\bm{v}_0,\varphi_0,\sigma_0)$ satisfying
$\bm {v}_{0} \in \bm{L}^2_{0,\sigma}(\Omega)$, $\varphi_{0}\in H^1(\Omega)$, $\|  \varphi_{0} \|_{L^{\infty}(\Omega)} \le 1$,
$|\overline{\varphi_{0}}|<1$, $\sigma_{0} \ln \sigma_{0} \in L^{1}(\Omega)$ and $\sigma_{0} \geq 0$ almost everywhere in $\Omega$,
a quadruple $(\bm{v},\varphi,\mu,\sigma)$ is called a weak solution to problem \eqref{sysNew}--\eqref{ini0new} on $[0,T]$, if it satisfies the regularity properties
\begin{align}
&\bm{v} \in L^{\infty}(0,T;\bm{L}^2_{0,\sigma}(\Omega)) \cap L^{2}(0,T;\bm{H}^1_{0,\sigma}(\Omega))\cap  W^{1,\frac{4}{3}}(0,T;(\bm{H}^1_{0,\sigma}(\Omega))'),
\notag \\
&\varphi \in L^{\infty}(0,T;H^1(\Omega))\cap L^{2}(0,T;H^2_{N}(\Omega)) \cap H^{1}(0,T;(H^1(\Omega))'),
\notag \\
&\varphi\in L^{\infty}(\Omega\times (0,T))\ \textrm{with}\ \ |\varphi(x,t)|<1\ \ \textrm{a.e.\ in}\ \Omega\times(0,T),
\notag\\
&\mu \in   L^{2}(0,T;H^1(\Omega)),\quad \Psi'(\varphi)\in L^2(0,T;L^2(\Omega)),
\notag \\
&\sigma \in L^{\infty}(0, T ; L^{1}(\Omega))\cap L^2(0,T;L^2(\Omega))\cap L^{\frac{4}{3}}(0, T ; W^{1,\frac{4}{3}}(\Omega)),
\notag\\
&\sigma(x, t) \geq 0 \quad \text {a.e. in}\ \Omega\times(0,T),
\notag
\end{align}
and the following equalities hold
\begin{subequations}
\begin{alignat}{3}
& \langle\partial_t  \bm{ v},\bm{\zeta} \rangle_{(\bm{H}^1_{0,\sigma}(\Omega))',\bm{H}^1_{0,\sigma}(\Omega)}
{\color{black}{- (\bm{v} \otimes\bm{ v},D\bm{ \zeta})}}
+\big(2\eta(\varphi) D\bm{v},D\bm{ \zeta}\big)
={\color{black}{\langle (\mu+\chi\sigma) \nabla \varphi, \bm {\zeta}\rangle_{\bm{L}^{\frac{6}{5}}(\Omega),\bm{L}^{6}(\Omega)}}},
\\
& \langle \partial_t \varphi,\xi \rangle_{(H^1(\Omega))',H^1(\Omega)}
+ ({\bm{v} \cdot \nabla \varphi},\xi)+ \big(m(\varphi)\nabla \mu,\nabla \xi\big)
=
\big(-\alpha \varphi+ h(\varphi, \sigma),\xi\big),\
\end{alignat}
\end{subequations}
almost everywhere in $(0,T)$ for any  $\bm {\zeta} \in \bm{H}^{1}_{0,\sigma}(\Omega)$, $\xi\in H^1(\Omega)$,
 \begin{align}
 & \mu= - \Delta \varphi +\Psi'(\varphi) -\chi \sigma,\qquad \textrm{a.e. in }\Omega\times(0,T),
 \label{test4.d}
 \end{align}
 and
\begin{align}
&-\int_0^T\!\int_\Omega\sigma\partial_t w\,\mathrm{d}x\mathrm{d}t
-\int_\Omega \sigma_0w(\cdot,0)\,\mathrm{d}x
-\int_0^T\!\int_\Omega \sigma\bm{v}\cdot\nabla w \,\mathrm{d}x\mathrm{d}t
-\int_0^T\!\int_\Omega \sigma \Delta w \,\mathrm{d}x \mathrm{d}t
\notag\\
&\quad =\chi\int_0^T\!\int_\Omega \sigma\nabla \varphi\cdot\nabla w \,\mathrm{d}x\mathrm{d}t
+\int_0^T\!\int_\Omega \big(\beta(\varphi)\sigma-\kappa\sigma^2\big) w \,\mathrm{d}x\mathrm{d}t,
  \label{test2.b}
\end{align}
 for any $w\in C^1([0,T];H^2_N(\Omega))$ with $w(T)=0$.
 Moreover,
\begin{align}
&\bm{v}|_{t=0}=\bm{v}_{0},\quad \varphi|_{t=0}=\varphi_{0}, \quad\textrm{a.e. in } \Omega.	
\notag
\end{align}
\ed
\begin{remark}\rm
Thanks to the regularity of $(\bm{v}, \varphi)$, their initial values are attained in the following sense (see, e.g., \cite{B}):
$\bm{v} \in C_w([0,T];\bm{L}^2_{0,\sigma}(\Omega))$, $\varphi\in C_w([0,T]; H^1(\Omega))$.
\end{remark}
	
We are now in a position to state the main results of this paper.
\bt(Existence of global weak solutions).
\label{main}
Suppose that the assumptions (H1)--(H6) are satisfied and $T \in (0,+\infty)$ is a given final time of arbitrary magnitude.

(1) Let $d =2,3$. For any given initial data $(\bm{v}_0,\varphi_0,\sigma_0)$ satisfying
\begin{align*}
&\bm {v}_{0} \in \bm{L}^2_{0,\sigma}(\Omega),\ \ \varphi_{0}\in H^1(\Omega),\ \ \|  \varphi_{0} \|_{L^{\infty}(\Omega)} \le 1,\ \
|\overline{\varphi_{0}}|<1,\\
 &\sigma_{0} \geq 0 \quad \textrm{a.e. in } \Omega,\ \quad \sigma_{0} \ln \sigma_{0} \in L^{1}(\Omega),
\end{align*}
problem \eqref{sysNew}--\eqref{ini0new} admits a global weak solution $(\bm{v},\varphi,\mu,\sigma)$ on $[0,T]$ in the sense of Definition \ref{maind}.

 (2) Let $d=2$. Assume in addition that $\sigma_0\in L^2(\Omega)$. We have
 \begin{align*}
 &\bm{v} \in  H^{1}(0,T;(\bm{H}^1_{0,\sigma}(\Omega))'),
 \\
 &\varphi\in L^4(0,T;H^2_N(\Omega))\cap L^2(0,T;W^{2,q}(\Omega)),\quad 
 \Psi'(\varphi)\in L^2(0,T;L^q(\Omega)),
 \\
 &\sigma\in C([0,T];L^2(\Omega))\cap L^2(0,T;H^1(\Omega))\cap H^1(0,T;(H^1(\Omega))')
 \end{align*}
 for any $q\in [2,+\infty)$. Furthermore,
 \begin{align}
 & \langle \partial_t \sigma,\xi\rangle_{(H^1(\Omega))',H^1(\Omega)}
 - (\bm{v} \sigma,\nabla\xi) + (\nabla \sigma,\nabla \xi) -\chi(\sigma\nabla\varphi,\nabla\xi)
= \big(\beta(\varphi)\sigma
 -\kappa \sigma^2,\xi\big)
 \notag
 \end{align}
holds almost everywhere in $(0,T)$ and for any $\xi\in H^1(\Omega)$. Besides, $\sigma|_{t=0}=\sigma_{0}$ almost everywhere in $\Omega$.
\et

\begin{remark}
Due to the low regularity of the chemical concentration $\sigma$ as well as the nonconstant mobility $m(\varphi)$, uniqueness of weak solutions remains an open question even if the spatial dimension is two (cf. \cite{CGGG}).
\end{remark}

Below we present a first-step result on the uniqueness of weak solutions in two dimensions, provided that the mobility $m(\varphi)$ is a positive constant, the mass source $h(\varphi, \sigma)$ is a constant and $\sigma_0\in L^2(\Omega)$. This extends the uniqueness result in \cite{RSS} for ``strong'' solutions of a fluid-free version associated with system \eqref{maineq}, as well as the uniqueness result in \cite{H}, where $\sigma$ satisfies a linear mass transport as in \eqref{diffu}.

\begin{theorem}[Uniqueness of global weak solution in two dimensions]
\label{main-2}
Let the assumptions in Theorem \ref{main}-(2) be satisfied. Moreover,  assume that
$$m(\varphi)\equiv \widehat{m},\quad h(\varphi,\sigma)\equiv\widehat{h},$$
where $\widehat{m}>0$, $\widehat{h}\in \mathbb{R}$ are given constants (recalling (H4), we now have $h^*=|\widehat{h}|$). Then, the global weak solution $(\bm{v},\varphi,\mu,\sigma)$ to problem \eqref{sysNew}--\eqref{ini0new} is unique.
\end{theorem}



\section{A Regularized System and its Semi-Galerkin Scheme}\label{ws}
\setcounter{equation}{0}

In this section, we introduce an appropriate regularized system for problem \eqref{sysNew}--\eqref{ini0new}, with a non-standard approximation for the singular potential $\Psi$ and a $p$-Laplace regularization in the chemical potential $\mu$. Then we solve the regularized system by a suitable semi-Galerkin scheme, specifically designed to handle the nonlinear coupling structure of the system. In summary, we proceed as follows: first, given a finite dimensional Galerkin ansatz $(\bm{u}^k,\psi^k)$ for the fluid velocity and the phase variable, we solve the advection-diffusion-reaction equation for the chemical concentration $\sigma^k$; next, with the sufficiently smooth solution $\sigma^k$, we solve the Faedo--Galerkin approximation of the regularized Navier--Stokes--Cahn--Hilliard system and obtain a finite dimensional solution $(\bm{v}^k,\varphi^k)$; finally, we prove the existence of a fixed point for the nonlinear mapping $\Phi: \Phi(\bm{u}^k,\psi^k)= (\bm{v}^k, \varphi^k)$ by Schauder's fixed point theorem, which yields a local solution $(\bm{v}^k, \varphi^k, \sigma^k)$ of the semi-Galerkin scheme.

\subsection{The regularized system}
\textbf{Regularization of the singular potential $\Psi$.} Let $\epsilon_0\in (0,1)$ be the constant given in (H1). There exists $\epsilon_1\in (0,\epsilon_0)$ such that
$$
\Psi_0'(-1+\epsilon)\leq -1,\quad \Psi_0'(1-\epsilon)\geq 1, \quad \forall \, \epsilon \in (0,\epsilon_1).
$$
Given $\epsilon \in (0,\epsilon_1)$, we introduce the following approximation of the singular nonlinearity $\Psi'_0$:
\begin{equation}
\Psi'_ {0,\epsilon}(r)=\left\{
\begin{aligned}
&\Psi_0'(-1+\epsilon) - \Psi_0''(-1+\epsilon) \left(\frac{4|\chi|+3}{4(|\chi|+1)}+\epsilon \right) &&\\
&\quad  - \Psi_0''(-1+\epsilon) \mathrm{e}^{-4(|\chi|+1)r-8(|\chi|+1)-\ln 4(|\chi|+1)},&&\quad r\le-2,\\
&\Psi_0'(-1+\epsilon) + \Psi_0''(-1+\epsilon)(r+1-\epsilon),&&\quad -2<r<-1+\epsilon,\\
&\Psi_0'(r),&&\quad \ |r|\leq 1-\epsilon,\\\\[2pt]
&\Psi_0'(1-\epsilon) + \Psi_0''(1-\epsilon)(r-1+\epsilon),&&\quad  1-\epsilon<r<2,\\
&\Psi_0'(1-\epsilon) + \Psi_0''(1-\epsilon) \left(\frac{4|\chi|+3}{4(|\chi|+1)}+\epsilon \right)
&&\\
&\quad + \Psi_0''(1-\epsilon) \mathrm{e}^{4(|\chi|+1)r-8(|\chi|+1)-\ln 4(|\chi|+1)},&&\quad\, r\geq 2.
\end{aligned}
\right.\notag
\end{equation}
Set
\begin{align*}
\Psi_{0,\epsilon}(r)=\int_0^{r}\Psi'_{0,\epsilon}(s) \, \mathrm{d}s
\qquad \text{and} \qquad
\Psi_\epsilon(r)=\Psi_{0,\epsilon}(r)-\frac{\theta_0}{2}r^2.
\end{align*}
It is straightforward to check that $\Psi_{0,\epsilon}'\in C^{1}(\mathbb{R})$ and thus $\Psi_{0,\epsilon}\in C^{2}(\mathbb{R})$. According to (H1), it holds
\begin{align*}
&\Psi_ {0,\epsilon}''(r)\ge \theta,\qquad r\Psi_{0,\epsilon}'(r) \geq 0,\qquad \Psi_ {0,\epsilon}(r)\geq -L,\qquad \forall\, r\in \mathbb{R},
\end{align*}
where $L>0$ is a constant independent of $\epsilon$ and $r$.
Similar to \cite{GGW}, we also find
$$
\Psi_{0,\epsilon}(r)\leq \Psi_{0}(r),\quad \forall\, r\in [-1,1].
$$
\begin{remark}
The approximate potential function $\Psi_{0,\epsilon}(r)$ depends on the chemical sensitivity $\chi$. For any given $\epsilon\in (0,\epsilon_1)$, $\Psi_{0,\epsilon}(r)$ provides a sufficiently fast growth for large values of $r$, while for $r\in [-1,1]$, it provides sufficient coercivity for small $\epsilon$. The former will be helpful to derive necessary estimates for approximate solutions, since the Galerkin ansatz of $\varphi$ does not have a uniform $L^\infty$-bound. The latter enables us to recover the physical constraint $\varphi\in [-1,1]$ after passing to the limit as $\epsilon\to 0^+$.
\end{remark}

\textbf{Approximation of the initial data.}
For any given $\varphi_{0}\in H^1(\Omega)$ with $\|  \varphi_{0} \|_{L^{\infty}(\Omega)} \le 1$ and $|\overline{\varphi_{0}}|<1$,
we define $\varphi_{0,\gamma}$ as the unique solution to the Neumann problem
\begin{equation*}
\begin{cases}
 \varphi_{0,\gamma} - \gamma \Delta \varphi_{0,\gamma}= (1-\gamma )\varphi_{0}, \, \ \quad \text{in}\ \Omega,
 \\
  \partial_{\bm{n}}  \varphi_{0,\gamma} =0,  \qquad  \qquad \qquad \qquad \ \ \text{on}\ \partial\Omega,
\end{cases}
\end{equation*}
where $\gamma\in \left( 0, \frac12 \right]$.
The classical elliptic theory yields
\begin{align}
\varphi_{0,\gamma}\in H^2_N(\Omega)\cap H^3(\Omega),\quad \text{with}\quad
\overline{\varphi_{0,\gamma}}= (1-\gamma)\overline{\varphi_{0}}\in \big[-|\overline{\varphi_{0}}|,\,|\overline{\varphi_{0}}|\big].
\label{Lvp1k}
\end{align}
Thanks to the maximum principle, we also find
\be
 \|  \varphi_{0,\gamma} \|_{L^\infty(\Omega)} \le 1-\gamma.
 \label{Lvp0k}
\ee
It is straightforward to check that
\be
\|\varphi_{0,\gamma}\| \leq \|\varphi_0\|,\quad  \|\nabla \varphi_{0,\gamma}\| \leq \|\nabla \varphi_0\|,\quad
\gamma \|\Delta \varphi_{0,\gamma}\|\leq 2 \|\varphi_0\|.
\label{Lvp2k}
\ee
Then by the elliptic estimate, we get
\be
\|\varphi_{0,\gamma}\|_{H^2(\Omega)}\leq C(\|\Delta \varphi_{0,\gamma}\|+\|\varphi_{0,\gamma}\|)\leq C\left(\frac{1}{\gamma}+1\right)\|\varphi_{0}\|,
\label{Lvp3k}
\ee
where the constant $C>0$ depends only on $\Omega$. From \eqref{Lvp2k} and Young's inequality, we infer that
\begin{align*}
\|\varphi_{0,\gamma} -\varphi_{0}\|^2
& =  \gamma \int_\Omega  ( \Delta \varphi_{0,\gamma}-  \varphi_{0})(\varphi_{0,\gamma} -\varphi_{0})\,\mathrm{d}x \\
& \leq \gamma \|\nabla  \varphi_{0,\gamma}\| \|\nabla (\varphi_{0,\gamma} -\varphi_{0})\| + \gamma \|\varphi_{0}\|\|\varphi_{0,\gamma} -\varphi_{0}\|
\\
&\leq  2 \gamma \|\nabla  \varphi_{0}\|^2+ \frac12 \|\varphi_{0,\gamma} -\varphi_{0}\|^2+ \frac{\gamma^2}{2}\|\varphi_{0}\|^2,
\end{align*}
which implies $\varphi_{0,\gamma}$ strongly converges to $\varphi_0$ in $L^2(\Omega)$ as $\gamma\to 0$. This fact combined with \eqref{Lvp2k} further yields that $\varphi_{0,\gamma}$ weakly converges to $\varphi_0$ in $H^1(\Omega)$ as $\gamma\to 0$. Since $H^1(\Omega)$ is uniformly convex, we apply \cite[Proposition 3.32]{B2010} and \eqref{Lvp2k} to conclude
\be
\lim_{\gamma\to 0}\|\varphi_{0,\gamma}-\varphi_{0}\|_{H^1(\Omega)} =0.
\notag
\ee

Next, for any $\sigma_0$ satisfying $\sigma_0\ln \sigma_0\in L^1(\Omega)$ and $\sigma_0\geq 0$ almost everywhere in $\Omega$, we consider a family of approximations $\{\sigma_{0,n}\}_{n\in \mathbb{Z}^+}$, with the following properties (see, e.g., \cite[Section 2.2]{W2016})
$$
\sigma_{0,n} \in C^{\infty}_0(\Omega),  \quad \sigma_{0,n}\ge 0 \ \text { in } \Omega,\quad \sigma_{0,n}\not\equiv 0,\quad
\sigma_{0,n} \rightarrow \sigma_0 \text { in }\ L \log L(\Omega) \ \text { as } n\to +\infty.
$$
Here, $L \log L(\Omega)$ denotes the standard Orlicz space associated with the Young function $(0, \infty) \ni z \mapsto$ $z \ln (1+z)$. Without loss of generality, we assume that
$$ \int_\Omega \sigma_{0,n}\ln \sigma_{0,n}\,\mathrm{d}x \leq \int_\Omega \sigma_{0}\ln \sigma_{0}\,\mathrm{d}x +1,\quad
\forall\,n\in \mathbb{Z}^+.
$$

\textbf{The regularized problem.} For any given $\gamma\in \left( 0,\frac12\right]$,  $\epsilon\in(0,\epsilon_1]$, $n\in \mathbb{Z}^+$, we introduce the following regularized problem $(\bm{S}_{\gamma,\epsilon,n})$:
\begin{subequations}
\label{reg.sys}
\begin{alignat}{3}
&\partial_t  \bm{ v}+\bm{ v} \cdot \nabla  \bm {v}-\mathrm{div} \big(  2\eta(\varphi) D\bm{v} \big)+\nabla p=(\mu+\chi \sigma)\nabla \varphi,
\label{reg.1} \\
&\mathrm{div}\, \bm{v}=0,
\label{reg.11}\\
&\partial_t\varphi+\bm{v}\cdot\nabla\varphi =\mathrm{div}\big(m(\varphi)\nabla \mu\big) -\alpha \varphi+ h(\varphi, \sigma),
\label{reg.2} \\
&\mu= -\gamma^{8} \mathrm{div}(|\nabla \varphi|^2\nabla \varphi)
- \Delta \varphi + \Psi_{\epsilon}'(\varphi) - \chi \sigma,
\label{reg.3}\\
&\partial_t\sigma+\bm{v}\cdot\nabla\sigma -\Delta \sigma + \chi \mathrm{div}
(\sigma\nabla  \varphi)
=\beta(\varphi) \sigma-\kappa \sigma^2,
\label{reg.4}
\end{alignat}
\end{subequations}
in $\Omega \times (0,T)$, subject to the boundary conditions
\begin{alignat}{3}
&\bm{v}=\mathbf{0},\quad\partial_{\boldsymbol{n}} \varphi = \partial_{\boldsymbol{n}} \mu = \partial_{\boldsymbol{n}}\sigma =0,
\quad \textrm{on}\ \partial\Omega\times(0,T),
\label{reg.boundary}
\end{alignat}
and the initial conditions
\begin{alignat}{3}
&\bm{v}|_{t=0}=\bm{v}_{0},
\quad \varphi|_{t=0}=\varphi_{0,\gamma},
\quad \sigma|_{t=0}=\sigma_{0,n},
\qquad &\textrm{in}&\ \Omega.
\label{reg.ini}
\end{alignat}

\begin{remark}
The (exponential) growth of $\Psi_{\epsilon}'$ is critical in two dimensions and super critical in three dimensions in view of the Sobolev embedding for $H^1(\Omega)$.
In order to overcome this difficulty, we introduce an additional regularizing term $-\gamma^{8} \mathrm{div}(|\nabla \varphi|^2\nabla \varphi)$ in \eqref{reg.3}. Together with the regularized initial datum $\varphi_{0,\gamma}$, the $p$-Laplace term yields a stronger estimate for the approximate solution $\varphi^k$ in the Faedo--Galerkin scheme presented  below such that $\varphi^k \in L^\infty(0,T; W^{1,4}(\Omega))$. By the Sobolev embedding theorem, we have that $\varphi^k \in L^\infty(0,T; L^\infty(\Omega))$ ($d=2,3$), which is crucial to control the nonlinear term $\Psi_{\epsilon}'(\varphi^k)$.
\end{remark}

\subsection{The semi-Galerkin scheme}

We shall solve the regularized system \eqref{reg.sys}--\eqref{reg.ini} by a suitable semi-Galerkin scheme. Consider the family of eigenvalues $\{\lambda_i\}_{i=1}^{\infty}$ and corresponding eigenfunctions $\{\bm{y}_{i}(x)\}_{i=1}^{\infty}$ of the Stokes problem
\be
(\nabla \bm{y}_{i},\nabla \bm{\zeta})=\lambda_{i}(\bm{y}_{i},\bm{\zeta}),\quad  \forall\, \bm{\zeta} \in {\bm{H}^1_{0,\sigma}(\Omega)},\quad \textrm{with}\ \|\bm{y}_{i}\|=1.
\notag
\ee
It is well-known that $0<\lambda_1\leq \lambda_2 \leq \cdots \to +\infty$, $\{\bm{y}_{i}\}_{i=1}^{\infty}$ forms a complete orthonormal basis in $\bm{L}^2_{0,\sigma}(\Omega)$ and it is orthogonal in $\bm{H}^1_{0,\sigma}(\Omega)$.
Next, we consider the family of eigenvalues $\{\ell_i\}_{i=1}^{\infty}$ and corresponding eigenfunctions $\{z_{i}(x)\}_{i=1}^{\infty}$ of the Laplace operator with homogeneous Neumann boundary conditions
\be
(\nabla z_{i},\nabla w)=\ell_{i}(z_{i},w),\quad  \forall\, w \in {H^1(\Omega)}, \quad \textrm{ with }  \|z_{i}\|=1.
\notag
\ee
Then $0=\ell_1<\ell_2 \leq \cdots \to +\infty$, $z_1=1$,  $\{z_{i}\}_{i=1}^{\infty}$ forms a complete orthonormal basis in $L^2(\Omega)$ and it is orthogonal in $H^1(\Omega)$.
For every integer $k\geq 1$, we denote the finite-dimensional subspace of $\bm{L}^2_{0,\sigma}(\Omega)$ by
$$ \bm{Y}_{k}:=\textrm{span} \{\bm{y}_{1}(x) ,\cdots,\bm{y}_{k}(x)\}.$$
The orthogonal projection on $\bm{Y}_{k}$ with respect to the inner product in $\bm{L}^2_{0,\sigma}(\Omega)$ is denoted by $\bm{P}_{\bm{Y}_{k}}$. Similarly, we denote the finite-dimensional subspace of $L^2(\Omega)$ by
$$ Z_{k}:=\textrm{span} \{z_{1}(x) ,\cdots,z_{k}(x)\}.$$
The orthogonal projection on $Z_{k}$ with respect to the inner product in $L^2(\Omega)$ is denoted by $\bm{P}_{Z_{k}}$.
We note that $\bigcup_{k=1}^\infty \bm{Y}_{k}$ is dense in both $\bm{L}^2_{0,\sigma}(\Omega)$, $\bm{H}^1_{0,\sigma}(\Omega)$ and $D(\bm{S})$, while  $\bigcup_{k=1}^\infty Z_{k}$ is dense in both $L^2(\Omega)$, $H^1(\Omega)$ and $H_N^2(\Omega)$. Recalling that $\Omega$ is assumed to be sufficiently smooth, for instance, a $C^4$-domain, we have $\bm{y}_{i}\in \bm{H}^1_{0,\sigma}(\Omega)\cap \bm{H}^4(\Omega)$ as well as $z_i\in H^2_N(\Omega)\cap H^4(\Omega)$ for all $i\in \mathbb{Z}^+$. Moreover, for any fixed $k\in \mathbb{Z}^+$, the following inverse inequalities hold
\begin{align*}
&\|\bm{v}\|_{\bm{H}^j(\Omega)}\leq C_k\|\bm{v}\|,\quad \forall\, \bm{v}\in \bm{Y}_k,\quad j=1,2,3,4,\\
& \|\varphi\|_{H^j(\Omega)}\leq C_k\|\varphi\|,\quad \forall\,\varphi\in Z_k,\quad j=1,2,3,4.
\end{align*}
The symbol $C_k$ denotes a generic positive constant that depends on the approximating parameter $k$.
\medskip

\textbf{Finite dimensional approximation of the initial data.}
Let $k$ be an arbitrary positive integer.
For the initial velocity field $\bm{v}_{0}\in \bm{L}^2_{0,\sigma}(\Omega)$, we consider its finite dimensional approximation $\bm{P}_{\bm{Y}_{k}} \bm{v}_{0}$, which satisfies
\begin{align*}
\lim_{k\to +\infty} \|\bm{P}_{\bm{Y}_{k}} \bm{v}_{0}-  \bm{v}_{0}\|=0.
\end{align*}
Next, we consider the finite dimensional approximation for the regularized initial datum $\varphi_{0,\gamma}$ of the phase variable, that is,  $\bm{P}_{Z_{k}}\varphi_{0,\gamma}$. We have $\bm{P}_{Z_{k}}\varphi_{0,\gamma}\in H^2_N(\Omega)$ and
$$
\lim_{k\to+\infty}\|\bm{P}_{Z_{k}}\varphi_{0,\gamma} -\varphi_{0,\gamma}\|_{H^2(\Omega)}=0.
$$
From the Sobolev embedding theorem $H^2(\Omega) \hookrightarrow C(\overline{\Omega})$ that is valid for $d=2,3$, for every given $\gamma\in
\left( 0,\frac12\right]$, there exists an (sufficiently large) integer $\widehat{k}$ such that
\be
\|\bm{P}_{Z_{\widehat{k}}}\varphi_{0,\gamma} -\varphi_{0,\gamma}\|_{C(\overline{\Omega})}\le C\|\bm{P}_{Z_{\widehat{k}}}\varphi_{0,\gamma}-\varphi_{0,\gamma}\|_{H^2(\Omega)} \le \frac{\gamma}{2},
\notag
\ee
where the constant $C>0$ only depends on $\Omega$. Thus, for all integers $k\geq \widehat{k}$, it holds
\be
\label{PZK}
\|\bm{P}_{Z_{k}}\varphi_{0,\gamma}\|_{L^\infty(\Omega)}\le 1-\frac{\gamma}{2},
\quad   \|\bm{P}_{Z_{k}}\varphi_{0,\gamma}\|_{H^1(\Omega)}
\leq   \|  \varphi_{0,\gamma}\|_{H^1(\Omega)}\leq \|  \varphi_0\|_{H^1(\Omega)}.
\ee

\textbf{The semi-Galerkin scheme.}
Let $T>0$ be an arbitrarily fixed final time.
For every integer $k\geq \widehat{k}$, we consider the approximate solution $(\bm{v}^k,\varphi^k,\mu^k,\sigma^k)$ to the regularized system \eqref{reg.sys}--\eqref{reg.ini} such that the Galerkin ansatz
\be
\bm{v}^{k}(x,t):=\sum_{i=1}^{k}a_{i}^{k}(t)\bm{y}_{i}(x),
\quad \varphi^{k}(x,t):=\sum_{i=1}^{k}b_{i}^{k}(t)z_{i}(x),
\quad
\mu^{k}(x,t):=\sum_{i=1}^{k}c_{i}^{k}(t)z_{i}(x),
\notag
\ee
satisfy
\begin{align}
&(\partial_t  \bm{ v}^{k},\bm{\zeta})+\big(( \bm{ v}^{k} \cdot \nabla)  \bm {v}^{k},\bm{ \zeta}\big)
+\big(  2\eta(\varphi^{k}) D\bm{v}^{k},D\bm{\zeta}\big)=\big((\mu^{k}+\chi 
\sigma^{k})\nabla \varphi^{k},\bm {\zeta}\big),
\label{atest.1}\\
&(\partial_t\varphi^k,\xi)+(\bm{v}^k\cdot\nabla\varphi^k,\xi)
=-\big(m(\varphi^k)\nabla \mu^k,\nabla \xi\big)+\big(-\alpha \varphi^k+h(\varphi^k,\sigma^k),\xi\big),
\label{atest.2}
\end{align}
and
\begin{equation}
\label{atest.3}
 (\mu^k,\xi)=\big(-\gamma^{8} \mathrm{div}( |\nabla \varphi^k|^2 \nabla \varphi^k) - \Delta \varphi^k + \Psi'_{\epsilon}(\varphi^k)- \chi \sigma^k, \xi\big),
\end{equation}
in $(0,T)$ and for all $\bm{\zeta} \in \bm{Y}_{k}$, $\xi\in Z_{k}$, while the unknown function  $\sigma^k$ satisfies
\begin{align}
&\partial_t\sigma^k+\bm{v}^k\cdot\nabla\sigma^k -\Delta\sigma^k+\chi \mathrm{div}(\sigma^k\nabla\varphi^k)
=\beta(\varphi^k) \sigma^k-\kappa(\sigma^k)^2,
\label{atest.4}
		\end{align}
in $\Omega\times(0,T)$. The approximate solution $(\bm{v}^k,\varphi^k,\mu^k,\sigma^k)$
fulfills the boundary and initial conditions
\begin{align}
 & \bm{v}^k=\bm{0},\quad \partial_{\bm{n}}\varphi^k=\partial_{\boldsymbol{n}} \mu^k=\partial_{\boldsymbol{n}} \sigma^k=0,&&\quad \textrm{on}\ \partial\Omega\times(0,T),
\label{boundary1}
\\
& \bm{v}^{k}|_{t=0}=\bm{P}_{\bm{Y}_{k}} \bm{v}_{0},\quad \varphi^{k}|_{t=0}=\bm{P}_{Z_{k}}\varphi_{0,\gamma},\quad  \sigma^{k}|_{t=0}=\sigma_{0,n},&&\quad \text{in}\ \Omega.
\label{atest.ini0}
\end{align}
%

The following proposition states that, given $\gamma\in \left( 0,\frac12\right]$,  $\epsilon\in (0,\epsilon_1]$, $n\in\mathbb{Z}^+$, for every integer $k\geq \widehat{k}$, the semi-Galerkin scheme \eqref{atest.1}--\eqref{atest.ini0} admits a unique local solution.
\bp[Local solvability of the semi-Galerkin scheme]\label{p1}
Let $d =2,3$. Suppose that the assumptions (H1)--(H6) are satisfied, $T \in (0,+\infty)$, and the initial data $(\bm{v}_0,\varphi_0,\sigma_0)$ are given as in Theorem \ref{main}. Given $\gamma\in \left( 0,\frac12\right]$, $\epsilon\in (0,\epsilon_1]$, $n\in \mathbb{Z}^+$, then for every integer $k\geq \widehat{k}$, the semi-Galerkin scheme \eqref{atest.1}--\eqref{atest.ini0} admits a unique local solution $(\bm{v}^{k},\varphi^{k},\mu^{k},\sigma^{k})$ on certain interval $[0,T_{k}]\subset[0,T]$ satisfying
\begin{align}
&\bm{v}^{k} \in H^1(0,T_{k};\bm{Y}_k),
\notag \\
&\varphi^{k} \in H^1(0,T_{k};Z_{k}),\quad \mu^k\in C([0,T_k];Z_{k}),
\notag \\
&\sigma^{k}\in C^{2,1}(\overline{\Omega}\times [0,T_k]),
\quad \sigma^k(x, t) > 0  \ \ \text{in } \ \Omega\times (0,T_k].
\notag
\end{align}
The existence time $T_{k}\in (0,T]$ depend on the initial data, $\Omega$, $k$, $\gamma$, $\epsilon$, $n$  and coefficients of the system.
\ep
%
\begin{proof}
The proof of Proposition \ref{p1} is based on a fixed point argument through Schauder's theorem. \medskip

\textbf{Step 1}. Let $\widetilde{M}$ be a sufficiently large constant satisfying $\widetilde{M}\geq 2(\|\bm{v}_0\|^2+  \|\varphi_0\|^2  +1)$. The exact value of $\widetilde{M}$ will be determined later. Consider two given functions
\be
\bm{u}^{k}=\sum_{i=1}^{k}\widetilde{a}_{i}^{k}(t)\bm{y}_{i}(x)\in C^\delta([0,T];\bm{Y}_{k}),\quad \psi^{k}=\sum_{i=1}^{k}\widetilde{b}_{i}^{k}(t)z_{i}(x)\in C^\delta([0,T];Z_{k})
\nonumber
\ee
for some $\delta\in \left( 0,\frac12\right)$, which fulfill
\begin{align}
&\widetilde{a}_{i}^{k}(0)=(\bm{v}_{0},\bm{y}_{i}),\quad i=1,\cdots,k,\quad
\text{and}\quad \sup_{t\in [0,T]}  \sum_{i=1}^{k}|\widetilde{a}_{i}^{k}(t)|^{2}\le \widetilde{M},
\nonumber\\
&\widetilde{b}_{i}^{k}(0)=(\varphi_{0,\gamma},z_{i}),\quad i=1,\cdots,k,\quad
\text{and}\quad \sup_{t\in [0,T]}  \sum_{i=1}^{k}|\widetilde{b}_{i}^{k}(t)|^{2}\le \widetilde{M}.
\nonumber
\end{align}
Then it holds $\bm{u}^{k}(0)=\bm{P}_{\bm{Y}_{k}} \bm{v}_{0}$, $\psi^{k}(0)=\bm{P}_{Z_{k}}\varphi_{0,\gamma}$ and
\begin{align}
 \sup_{t\in [0,T]}\|\bm{u}^{k}(t)\|^2\leq \widetilde{M},
 \quad
 \sup_{t\in [0,T]}\|\psi^{k}(t)\|^2\leq \widetilde{M}.
 \label{aLL-up}
\end{align}
Let us first consider the following auxiliary equation for the chemical concentration
\begin{alignat}{3}
&\partial_t\sigma^k+\bm{u}^k\cdot\nabla\sigma^k -\Delta\sigma^k+\chi \mathrm{div}(\sigma^k\nabla \psi^k)=\beta(\psi^k)  \sigma^k-\kappa (\sigma^k)^2,
\label{1atest.4}
\end{alignat}
in $\Omega\times (0,T)$,
equipped with the boundary and initial conditions
\begin{alignat}{3}
&\partial_{\boldsymbol{n}} \sigma^k =0  \quad \text{on}\ \partial\Omega\times(0,T),
\qquad \sigma^{k}|_{t=0}=\sigma_{0,n} \quad \text{in}\ \Omega.
\label{boundary2}
\end{alignat}
Then we have
\bl\label{fp}
Given $(\bm{u}^k, \psi^k)$, problem \eqref{1atest.4}--\eqref{boundary2} admits a unique classical solution $\sigma^{k}$ on $[0,T]$ such that
$$
\sigma^{k}\in C(\overline{\Omega}\times [0,T])\cap C^{2,1}(\overline{\Omega}\times (0,T)),
\quad \sigma^k(x, t) > 0  \ \ \text{in } \ \overline{\Omega}\times (0,T].
$$
\el

The proof of Lemma \ref{fp} can be found in the Appendix \ref{app-1}.
According to Lemma \ref{fp}, we find that the following mapping is well defined
\begin{align}
\Phi^k_{1}:\  C^\delta([0,T];\bm{Y}_{k})\times C^\delta([0,T];Z_{k}) \ &\to\ \  \widehat{X}, \notag\\
(\bm{u}^{k},\psi^{k})\ &\mapsto\ \ \sigma^{k},\notag
\end{align}
where
\begin{align}
\widehat{X}&=L^{\infty}(0,T;H^1(\Omega)) \cap  L^{2}(0,T;H^2_N(\Omega))\cap H^1(0,T;L^2(\Omega)).
\notag
\end{align}
As seen in Appendix \ref{app-1}, $\Phi^k_{1}$ is a bounded operator from $C([0,T];\bm{Y}_{k})\times C([0,T];Z_{k})$ to $\widehat{X}$.

Next, we show that $\Phi^k_{1}$ is continuous with respect to the given data $(\bm{u}^{k},\psi^{k})$ in the topology of $X$, with
\begin{align}
X&= L^{\infty}(0,T;L^2(\Omega)) \cap  L^{2}(0,T;H^1(\Omega)).
\notag
\end{align}
To this end, let $\bm{u}^{k}_{1}$, $\bm{u}^{k}_{2}$ be two given vectorial functions with the same initial value $\bm{P}_{\bm{Y}_{k}} \bm{v}_{0}$, while $\psi^{k}_{1}$, $\psi^{k}_{2}$ be two given scalar functions with the same initial value $\bm{P}_{Z_{k}}\varphi_{0,\gamma}$. Both $(\bm{u}^{k}_{i},\psi^{k}_{i})$, $i=1,2$, satisfy the condition \eqref{aLL-up}.
Let $\sigma^{k}_{i} =\Phi^k_{1}(\bm{u}^{k}_{i},\psi^{k}_i)$, $i=1,2$, be the two corresponding solutions to problem \eqref{1atest.4}--\eqref{boundary2} given by Lemma \ref{fp} (with the same initial datum $\sigma_{0,n}$).
We denote their differences by
$$
\bm{u}^k= \bm{u}^k_1-\bm{u}^k_2,\quad \psi^{k}=\psi_{1}^{k}-\psi_{2}^{k},\quad \sigma^{k}=\sigma_{1}^{k}-\sigma_{2}^{k},
$$
which fulfill
\begin{align}
	&\partial_t  \sigma^{k}+{\bm{u}^k_1 \cdot \nabla \sigma^k}+{\bm{u}^k \cdot \nabla \sigma^k_2}-\Delta \sigma^k
	\notag\\
	&\quad =-\chi  \mathrm{div}\big(\sigma^k \nabla \psi_{1}^k+\sigma_{2}^k \nabla \psi^k\big)
+\big[\beta(\psi^k_1) -\kappa  (\sigma_1^k+\sigma_2^k)\big]\sigma^k
+\big(\beta(\psi^k_1)-\beta(\psi^k_2)\big) \sigma^k_2,
\label{2atest.3}
	\end{align}
in $\Omega\times (0,T)$, subject to the boundary and initial conditions
\begin{alignat}{3}
&\partial_{\boldsymbol{n}} \sigma^k =0  \quad \text{on}\ \partial\Omega\times(0,T),
\qquad \sigma^{k}|_{t=0}=0 \quad \text{in}\ \Omega.
\label{boundary3}
\end{alignat}
Testing \eqref{2atest.3} by $\sigma^k $, integrating over $\Omega$, using the fact $\mathrm{div}\,\bm{u}^k_1=0$ and integration by parts, we get
\begin{align}
&\frac{1}{2}\frac{\mathrm{d}}{\mathrm{d}t} \|\sigma^k\|^2+\|\nabla\sigma^k\|^2
\notag \\
&\quad =
(\bm{u}^{k} \sigma^{k}_2, \nabla\sigma^k)
+\chi \int_{\Omega} \sigma^k \nabla \psi_{1}^k \cdot \nabla \sigma^k \,\mathrm{d}x
+\chi \int_{\Omega} \sigma_{2}^k  \nabla \psi^k  \cdot \nabla \sigma^k \,\mathrm{d}x
\notag \\
&\qquad
+\int_{\Omega} \big[ \beta(\psi^k_1) -\kappa  (\sigma_1^k+\sigma_2^k)\big](\sigma^k)^2\,\mathrm{d}x
 + \int_{\Omega} \big(\beta(\psi^k_1)-\beta(\psi^k_2)\big) \sigma^k_2\sigma^k\,\mathrm{d}x
\notag \\
&\quad =: \sum_{j=1}^5 I_j.
\label{diffsig}
\end{align}
The right-hand side of \eqref{diffsig} can be estimated as follows. First,  it follows from H\"{o}lder's inequality, Young's inequality and the Sobolev emebdding theorem that
\begin{align}
I_1
&\le C\|\bm{u}^{k}\|_{\bm{L}^\infty(\Omega)}\|\sigma^{k}_{2}\|\|\nabla \sigma^{k}\|\notag\\
&\le C \|\bm{u}^{k}\|_{\bm{H}^2(\Omega)}\|\sigma^{k}_{2}\|\|\nabla\sigma^k\|\notag\\
& \le  C_k\|\sigma^{k}_{2}\|^2\|\bm{u}^{k}\|^2+\frac16 \|\nabla\sigma^k\|^2,
\notag
\end{align}
\be
\begin{aligned}
I_2
&\leq |\chi|\|\sigma^k\|  \|\nabla\psi_1^k\|_{\bm{L}^{\infty}(\Omega)} \|\nabla\sigma^k\|
\\
&\leq C\chi^2\|\psi_1^k\|_{H^3(\Omega)}^2\|\sigma^k\|^2
+\frac16\|\nabla\sigma^k\|^2
\\
&\leq C_k\chi^2 \|\psi_1^k\|^2\|\sigma^k\|^2 +\frac16\|\nabla\sigma^k\|^2,
\end{aligned}
\notag
\ee
and in a similar manner,
\be
\begin{aligned}
I_3&\leq C\chi^2 \|\sigma^k_2\|^2 \| \psi^k\|_{H^3(\Omega)}^2
+\frac16\|\nabla\sigma^k\|^2
\\
&\leq C_k\chi^2\|\sigma^k_2\|^2\|\psi^k\|^2 +\frac16\|\nabla\sigma^k\|^2.
\end{aligned}
\notag
\ee
Here, we have essentially used the facts that $\bm{u}_i^k$ and $\psi_i^k$, $i=1,2$, are finite dimensional. Next, we infer from (H5) and the nonnegativity of $\sigma_i^k$, $i=1,2$, that
\be
\begin{aligned}
	I_4&\leq b^*\|\sigma^k\|^2,
\end{aligned}
\notag
\ee
as well as
\be
\begin{aligned}
I_5  &\leq C\|\psi^k\|_{L^{\infty}(\Omega)}  \| \sigma^k_2 \| \|\sigma^k\|
\\
&\leq \|\psi^k\|_{H^2(\Omega)}^2 +C \|\sigma^k_2\|^2  \|\sigma^k\|^2
\\
&\leq C_k\|\psi^k\|^2+C \|\sigma^k_2\|^2  \|\sigma^k\|^2.
\end{aligned}
\notag
\ee
Combining the above estimates and integrating \eqref{diffsig} on $[0, t] \subset [0, T] $, we obtain
\be
\begin{aligned}
&  \|\sigma^k(t) \|^{2}
+\int_0^t  \|\nabla\sigma^k(s)\|^{2}\,\mathrm{d}s
  \leq C_k \int_0^t \|\sigma^k(s)\|^{2}\,\mathrm{d}s
+C_k \int_0^t  \left( \|\bm{u}^{k}(s)\|^2 + \| \psi^k(s)\|^{2} \right) \,\mathrm{d}s,
\label{sig}
\end{aligned}
\ee
where the estimate $\|\sigma_{2}^{k}\|_{L^{\infty}(0,T;L^2(\Omega)) } \leq C$ has been used.
An application of Gronwall's lemma to \eqref{sig} yields that
\begin{align}
& \|\sigma^k(t)\|^{2} +\int_0^t \|\nabla\sigma^k(s)\|^2\, \mathrm{d}s
\le C_T\Big(\sup_{s\in [0,t]}\|\bm{u}^{k}(s)\|^2+\sup_{s\in [0,t]}\|\psi^{k}(s)\|^2\Big),
\quad \forall\,t\in[0,T],
\label{auphih1}
\end{align}
where the constant $C_T>0$ depends on $T$ and $k$.
As a consequence, the solution operator $\Phi^k_{1}$ is continuous with respect to $(\bm{u}^{k},\psi^k)$ as a mapping from $C([0,T];\bm{Y}_{k})\times C([0,T];Z_{k})$ to $X$.
\medskip

\textbf{Step 2}. Given the function $\sigma^{k}=\Phi^k_{1}(\bm{u}^k,\psi^k)$ obtained in Step 1, we now look for the ansatz
\be
\bm{v}^k=\sum_{i=1}^{k}a_{i}^{k}(t)\bm{y}_{i}(x),\quad \varphi^k=\sum_{i=1}^{k}b_{i}^{k}(t)z_{i}(x),\quad  \mu^{k}(x,t):=\sum_{i=1}^{k}c_{i}^{k}(t)z_{i}(x),
\nonumber
\ee
that satisfy the following auxiliary system for the fluid velocity and phase variable:
\begin{align}
	& (\partial_t  \bm{v}^k,\bm{\zeta})
	+\big(( \bm{v}^k \cdot \nabla)  \bm{v}^k,\bm{ \zeta}\big)
    +\big(  2\eta(\varphi^{k}) D\bm{v}^k, D\bm{\zeta}\big)
=\big((\mu^{k}+\chi
\sigma^{k})\nabla \varphi^{k},\bm {\zeta}\big),
	\label{aatest3.c}\\
	&  (\partial_t  \varphi^{k},\xi)+( \bm{v}^{k} \cdot \nabla  \varphi^{k},\xi) +\big(m(\varphi^k)\nabla \mu^{k},\nabla \xi\big)
=\big(-\alpha \varphi^k + h(\varphi^k,\sigma^k),\xi\big),
\label{g1.a}\\
& (\mu^k,\xi)
=\big(-\gamma^{8} \mathrm{div}( |\nabla \varphi^k|^2 \nabla \varphi^k)
- \Delta \varphi^k + \Psi'_{\epsilon}(\varphi^k)- \chi \sigma^k, \xi\big),
\label{g4.d}
\end{align}
 in $(0,T)$ and for all $\bm{\zeta} \in \bm{Y}_{k}$, $\xi\in Z_{k}$.
In addition, $(\bm{v}^k,\varphi^k,\mu^k)$
satisfies the boundary and initial conditions
\begin{align}
 & \bm{v}^k=\bm{0},\quad \partial_{\bm{n}}\varphi^k=\partial_{\boldsymbol{n}} \mu^k=0,&&\quad \textrm{on}\ \partial\Omega\times(0,T),
\label{aatest.boundary1}
\\
& \bm{v}^{k}|_{t=0}=\bm{P}_{\bm{Y}_{k}} \bm{v}_{0},\quad \varphi^{k}|_{t=0}=\bm{P}_{Z_{k}}\varphi_{0,\gamma},&&\quad \text{in}\ \Omega.
\label{aatest3.cini}
\end{align}
Then we have
\bl\label{NSSa}
Given $\sigma^k=\Phi^k_{1}(\bm{u}^k,\psi^k)$, the Faedo--Galerkin scheme \eqref{aatest3.c}--\eqref{aatest3.cini} admits a unique solution   $(\bm{v}^k,\varphi^k,\mu^k)$ on $[0,T]$ such that
$$
\bm{v}^k\in C^1([0,T];\bm{Y}_{k}),\quad
\varphi^k \in C^1([0,T];Z_{k}),\quad
\mu^k\in C([0,T];Z_{k}).
$$
\el

The proof of Lemma \ref{NSSa} can be found in the Appendix \ref{app-2}.
Thanks to Lemma \ref{NSSa}, we define the following mapping, which is  determined by the unique solution to problem  \eqref{aatest3.c}--\eqref{aatest3.cini}:
\begin{align*}
\Phi^k_{2}:\quad  \widehat{X} &\to\  H^1(0,T;\bm{Y}_{k})\times H^1(0,T;Z_{k}),
\\
\sigma^{k} &\mapsto\ (\bm{v}^k,\varphi^k).
\end{align*}
As seen  in Appendix \ref{app-2},
$\Phi^k_2$ is a bounded operator from $\widehat{X}$ to $H^1(0,T;\bm{Y}_{k})\times H^1(0,T;Z_{k})$.

Next, we establish the continuity of $\Phi^k_2$. Given $\sigma^{k}_{i}\in \widehat{X}$, $i=1,2$, we define the corresponding solutions $(\bm{v}^{k}_{i},\varphi^{k}_{i})=\Phi^k_2(\sigma^{k}_{i})$, $i=1, 2$, as in Lemma \ref{NSSa} (with the same initial data $(\bm{P}_{\bm{Y}_{k}} \bm{v}_{0}, \bm{P}_{Z_{k}}\varphi_{0,\gamma})$).
The corresponding chemical potentials are denoted by $\mu_i^k$, $i=1,2$, respectively. As before, we denote
 $$
 \sigma^{k}=\sigma_{1}^{k}-\sigma_{2}^{k},\quad
 \bm{v}^{k}=\bm{v}^{k}_{1}-\bm{v}^{k}_{2},\quad
 \varphi^{k}=  \varphi_{1}^{k}-\varphi_{2}^{k},\quad
 \mu^{k}=  \mu_{1}^{k}-\mu_{2}^{k}.
 $$
Taking the difference of \eqref{aatest3.c} for $\bm{v}^{k}_{i}$, testing the resultant by $\bm{\zeta}=\bm{v}^k$, using integration by parts, we obtain
\begin{align}
\frac12 \frac{\mathrm{d}}{\mathrm{d}t}  \|\bm{v}^{k}\|^2
	+\big(  2\eta(\varphi^{k}_{1}) D\bm{v}^{k} ,D\bm{v}^k\big)
&= 	-\big(( \bm{ v}^{k} \cdot \nabla)  \bm {v}^{k}_{1},\bm{v}^k\big)
-\big(2(\eta(\varphi^{k}_{1})-\eta(\varphi^{k}_{2})) D\bm{v}^{k}_{2} ,D\bm{v}^k\big)
\notag\\
&\quad +\big((\mu^{k}_{1}+\chi
\sigma^k_1) \nabla\varphi^{k}_{1}-(\mu^{k}_{2}+\chi
\sigma^k_2) \nabla\varphi^{k}_{2},\bm{v}^k\big)
\notag\\
&  =:\sum_{j=6}^{8} I_{j}.
\label{diffvm}
\end{align}
The first two terms on the right-hand side of \eqref{diffvm} can be estimated as follows:
\begin{align}
I_{6}
&\le  \|\nabla \bm{v}^{k}_{1}\|_{\bm{L}^\infty(\Omega)}
\|\bm{v}^k\|^2
 \leq C\|\bm{v}^{k}_{1}\|_{\bm{H}^3(\Omega)}\|\bm{v}^k\|^2
 \leq C_k\|\bm{v}^{k}_{1}\| \|\bm{v}^k\|^2,
\notag
\end{align}
and
\begin{align}
I_{7}
&\le C\|\varphi^k\|_{L^\infty(\Omega)}\|D \bm{v}^k_2\|\|D \bm{v}^k\|
\notag\\
& \le C\|\varphi^k\|_{H^2(\Omega)}  \|D\bm{v}^k_2\|\|D \bm{v}^k\|
\notag\\
& \le \frac{\eta_*}{3}\|D\bm{v}^k\|^2 + C_k\|\bm{v}^k_2\|^2\|\varphi^{k}\|^2.
\notag
\end{align}
Taking $\xi=\mu_1^k$ in \eqref{g4.d} for the solution $\varphi_1^k$, we find
\begin{align*}
\| \mu_1^k\|^2
&\leq \left( \gamma^{8} \left\| |\nabla \varphi_1^k|^2 \nabla \varphi_1^k \right\| + \| \nabla \varphi_1^k \| \right) \| \nabla \mu_1^k\| + \left( \| \Psi'_\epsilon(\varphi_1^k) \| + |\chi| \| \sigma_1^k\| \right) \| \mu_1^k\|
\\
&\leq C_k \left( \gamma^{8} \left\| |\nabla \varphi_1^k|^2 \nabla \varphi_1^k \right\| + \| \nabla \varphi_1^k \| \right) \| \mu_1^k\| +   \left( \| \Psi'_\epsilon(\varphi_1^k) \| + |\chi|\| \sigma_1^k\| \right) \| \mu_1^k\|.
\end{align*}
Since for any given $k$ there exists some $\widetilde{C}_k>0$ such that $\| \varphi_i^k \|_{C([0,T]; H^4(\Omega))}\leq \widetilde{C}_k$, $i=1,2$, it follows from the Sobolev embedding theorem and the construction of $\Psi_\epsilon'$ that
\begin{equation}
\label{muk-1}
\| \mu_1^k\|\leq  \widetilde{C}_k+ C \| \sigma_1^k\|.
\end{equation}
Similar result also holds for $\mu_2^k$. Next, taking $\xi=\mu^k$ in \eqref{g4.d} for the difference $\mu^k$, we easily obtain
\begin{align}
\| \mu^k\|^2
&\leq C_k \left( \gamma^{8} \left\| |\nabla \varphi_1^k|^2 \nabla \varphi_1^k- |\nabla \varphi_2^k|^2 \nabla \varphi_2^k \right\| + \| \nabla \varphi^k \| \right) \| \mu^k\|
\notag
\\
&\quad + \big( \| \Psi'_\epsilon(\varphi_1^k) - \Psi'_\epsilon(\varphi_2^k)\| + |\chi| \| \sigma^k\| \big) \| \mu^k\|,
\notag
\end{align}
which entails that
\begin{equation}
\label{muk-dif}
\| \mu^k\|\leq  \widetilde{C}_k \| \varphi^k\| + C \| \sigma^k\|.
\end{equation}
As a consequence, writing
$$
I_8= \big((\mu^{k}_{1}+\chi
\sigma^k_1) \nabla\varphi^{k}, \bm{v}^k\big)
+ (\mu^{k}
+\chi \sigma^k) \nabla\varphi^{k}_{2},\bm{v}^k\big),
$$
we deduce that
\begin{align*}
|I_8|
&\leq C \big( \| \mu_1^k\| +\| \sigma_1^k\| \big) \| \nabla \varphi^k\|_{\mathbf{L}^\infty(\Omega)} \| \bm{v}^k\|
+\big( \| \mu^k\| +\| \sigma^k\| \big) \| \nabla \varphi_2^k\|_{\mathbf{L}^\infty(\Omega)} \| \bm{v}^k\|
\\
&\leq C \big(\widetilde{C}_k+ C  \| \sigma_1^k\| \big)  \| \varphi^k\|_{H^3(\Omega)} \| \bm{v}^k\|
+ C\big( \widetilde{C}_k \| \varphi^k\| + C \| \sigma^k\| \big)\| \varphi_2^k\|_{H^3(\Omega)} \| \bm{v}^k\|
\\
&\leq C(\| \sigma_1^k\|, \| \varphi_2^k\|, k) \big( \| \varphi^k\|^2 + \| \bm{v}^k\|^2 \big)
+ C_k  \| \sigma^k\|^2.
\end{align*}

Next, taking the difference of \eqref{g1.a} for $\varphi^k_i$, we obtain
\begin{align}
(\partial_t  \varphi^{k},\xi)+\alpha(\varphi^k,\xi)
&=-(\varphi^{k}\bm{v}^{k}_{1},\nabla \xi)-(\varphi^{k}_{2}\bm{v}^{k},\nabla \xi)
- \big(m(\varphi^k_1)\nabla \mu^{k},\nabla \xi\big)
\notag \\
&\quad - \big((m(\varphi^k_1)-m(\varphi^k_2))\nabla \mu^{k}_2,\nabla \xi\big)
+ \big(h(\varphi^k_1,\sigma^k_1)-h(\varphi^k_2,\sigma^k_2), \xi\big),
\label{atest111.a}
\end{align}
for any $\xi\in  Z_{k}$. Choosing $\xi=1$ in \eqref{atest111.a} yields that
$$
\frac{\mathrm{d}}{\mathrm{d}t} \overline{\varphi^{k}} + \alpha  \overline{\varphi^{k}} = \frac{1}{|\Omega|}\int_\Omega \big[h(\varphi^k_1,\sigma^k_1)-h(\varphi^k_2,\sigma^k_2)\big]\,\mathrm{d}x.
$$
Multiplying the above identity by $\overline{\varphi^{k}}$ and using Young's inequality with the Lipschitz continuity property of $h$, we get
\be
\frac12 \frac{\mathrm{d}}{\mathrm{d}t} \big|\overline{\varphi^{k}}\big|^2 + \alpha \big|\overline{\varphi^{k}}\big|^2
\leq \big|\overline{\varphi^{k}}\big|^2 + C\big(\|\varphi^{k}\|^2+ \|\sigma^k\|^2\big).
\label{diff-ave}
\ee
Taking $\xi=\mathcal{N}\big(\varphi^{k}-\overline{\varphi^{k}}\big)$ in \eqref{atest111.a}, we obtain
\begin{align}
& \frac{1}{2}\frac{\mathrm{d}}{\mathrm{d}t} \|\varphi^{k}-\overline{\varphi^{k}}\|_{V_0^{-1}}^2 +\alpha\|\varphi^k-\overline{\varphi^{k}}\|_{V_0^{-1}}
\notag \\
&\quad= -\big(\varphi^{k}\bm{v}^{k}_{1},\nabla \mathcal{N}\big(\varphi^{k}-\overline{\varphi^{k}}\big)\big) -\big(\varphi^{k}_{2}\bm{v}^{k},\nabla\mathcal{N}\big(\varphi^{k} -\overline{\varphi^{k}}\big)\big)
\notag \\
&\qquad - \big(m(\varphi^k_1)\nabla \mu^{k},\nabla \mathcal{N}\big(\varphi^{k}-\overline{\varphi^{k}}\big)\big)
- \big((m(\varphi^k_1)-m(\varphi^k_2))\nabla \mu^{k}_2,\nabla \mathcal{N}\big(\varphi^{k}-\overline{\varphi^{k}}\big)\big)
\notag \\
&\qquad  + \big(h(\varphi^k_1,\sigma^k_1)-\overline{h(\varphi^k_1,\sigma^k_1)} -h(\varphi^k_2,\sigma^k_2) +\overline{h(\varphi^k_2,\sigma^k_2)}, \mathcal{N}\big(\varphi^{k}-\overline{\varphi^{k}}\big)\big)
\notag \\
&\quad   =: \sum_{j=9}^{13} I_j.
\label{vhmna}
\end{align}
The first two terms on the right-hand side of \eqref{vhmna} can be estimated as follows:
\begin{align}
I_{9}
&\le \|\varphi^k\|_{L^6(\Omega)}\|\bm{v}^k_1\|_{\bm{L}^3(\Omega)} \|\mathcal{N}\big(\varphi^{k}-\overline{\varphi^{k}}\big)\|
\nonumber\\
&\le C\|\varphi^k\|_{H^1(\Omega)}\|\bm{v}^k_{1}\|_{\bm{H}^1(\Omega)} \|\varphi^{k}-\overline{\varphi^{k}}\|_{V_0^{-1}}
\notag\\
&\le \|\bm{v}^k_{1}\|^2 \|\varphi^{k}-\overline{\varphi^{k}}\|_{V_0^{-1}}^2 +C_k\|\varphi^{k}\|^2,
\notag
\end{align}
and
\begin{align}
I_{10}
&\le \|\varphi^k_2\|_{L^3(\Omega)}\|\bm{v}^k\|_{\bm{L}^6(\Omega)} \|\nabla\mathcal{N}\big(\varphi^{k} -\overline{\varphi^{k}}\big)\|
\notag \\
& \le\frac{\eta_*}{3}\|D\bm{v}^{k}\|^2 +C_k\|\varphi^k_2\|^2\|\varphi^{k}-\overline{\varphi^{k}}\|_{V_0^{-1}}^2.
\notag
\end{align}
Using \eqref{muk-dif}, we deduce that
\begin{align}
 I_{11}
&=\big(m'(\varphi^k_1)\mu^k\nabla \varphi^k_1 ,\nabla\mathcal{N}\big(\varphi^{k} -\overline{\varphi^{k}}\big)\big)
- \big(m(\varphi^k_1)\mu_k,  \varphi^{k} -\overline{\varphi^{k}} \big)
\notag \\
&\leq C\|\mu^{k}\|\|\nabla \varphi^k_1\|_{\bm{L}^\infty(\Omega)} \|\nabla \mathcal{N}\big(\varphi^{k}-\overline{\varphi^{k}}\big)\|
+ C \|\mu^{k}\|\|\varphi^{k} -\overline{\varphi^{k}}\|
\notag\\
&\leq C \big(\|\varphi^k_1\|_{H^3(\Omega)} +1\big) \big(\widetilde{C}_k \| \varphi^k\| + C \| \sigma^k\|\big)\| \varphi^{k}-\overline{\varphi^{k}} \|
\notag \\
&\leq C_k(\|\varphi^k_1\|+1)\big(\|\varphi^{k}\|^2+ \|\sigma^k\|^2\big),
\notag
\end{align}
while with a similar estimate \eqref{muk-1} for $\mu_2^k$, we find
\begin{align}
I_{12}
&\leq \|m(\varphi^k_1)-m(\varphi^k_2)\|_{L^\infty(\Omega)}\|\nabla \mu^{k}_2\|\|\nabla \mathcal{N}\big(\varphi^{k}-\overline{\varphi^{k}}\big)\|
\notag \\
&\leq C_k
\|\varphi^k\|
 \big( \widetilde{C}_k+ C \| \sigma_2^k\|\big)\|\varphi^{k}-\overline{\varphi^{k}}\|_{V_0^{-1}}
\notag \\
&\leq C(\|\sigma^k_2\|,k)\|\varphi^k\|^2.
\notag
\end{align}
Finally, it holds
\begin{align}
I_{13}
&\leq C ( \|\varphi^k\|+\|\sigma^k\|) \|\varphi^{k}-\overline{\varphi^{k}}\|_{V_0^{-1}}
\notag \\
&\leq \|\varphi^k\|^2+\|\sigma^k\|^2+ C\|\varphi^{k}-\overline{\varphi^{k}}\|_{V_0^{-1}}^2. \notag
\end{align}
In the above estimates, we have again essentially used the fact that $\bm{v}^k$, $\varphi^k$ and $\mu^k$ are finite dimensional. Hence, all the related estimates for higher order norms depend on the parameter $k$ at this stage.

The Poincar\'{e}--Wirtinger inequality yields
\begin{align}
\|\varphi^k\|^2
& \leq 2\|\varphi^k-\overline{\varphi^{k}}\|^2+ 2|\Omega| |\overline{\varphi^{k}}|^2
\notag \\
&\leq C\|\nabla (\varphi^k-\overline{\varphi^{k}})\|  \|\varphi^{k}-\overline{\varphi^{k}}\|_{V_0^{-1}}+ 2|\Omega| |\overline{\varphi^{k}}|^2
\notag\\
&\leq C_k\| \varphi^k \|  \|\varphi^{k}-\overline{\varphi^{k}}\|_{V_0^{-1}}+ 2|\Omega| |\overline{\varphi^{k}}|^2
\notag\\
&\leq \frac12 \| \varphi^k \|^2 +C_k  \|\varphi^{k}-\overline{\varphi^{k}}\|_{V_0^{-1}}^2 + 2|\Omega| |\overline{\varphi^{k}}|^2.
\notag
\end{align}
Collecting the above estimates, we deduce from \eqref{diffvm}, \eqref{diff-ave} and \eqref{vhmna} that
\begin{align}
& \frac{1}{2}\frac{\mathrm{d}}{\mathrm{d}t}
\left(\|\bm{v}^k\|^2+ \|\varphi^{k}-\overline{\varphi^{k}}\|_{V_0^{-1}}^2+ \big|\overline{\varphi^{k}}\big|^2\right)
+\eta_*\|D\bm{v}^{k}\|^2
\notag \\
&\quad \le C_k\left(\|\bm{v}^k\|^2+ \|\varphi^{k}-\overline{\varphi^{k}}\|_{V_0^{-1}}^2+ \big|\overline{\varphi^{k}}\big|^2\right)
+ C_k \|\sigma^k\|^2.
\label{astarphi}
\end{align}
Applying Gronwall's lemma to \eqref{astarphi}, we obtain (recalling that  $\bm{v}^k|_{t=0}=\bm{0},\,\varphi^k|_{t=0}=0$)
\begin{align}
&\|\bm{v}^k(t)\|^2+ \|\varphi^{k}(t)-\overline{\varphi^{k}}(t)\|_{V_0^{-1}}^2+ \big|\overline{\varphi^{k}(t)}\big|^2
 \le C_k e^{C_k t} \int_0^t\| \sigma^k(s)\|^2\,\mathrm{d}s,
\label{aum}
\end{align}
for all $t\in [0,T]$.
Since $\bm{v}^k$ and $\varphi^k$ are finite dimensional, we conclude from \eqref{aum} that the solution operator $\Phi^k_2$ is continuous as a mapping from $X$ to $C([0,T];\bm{Y}_{k})\times C([0,T];Z_{k})$.
\medskip

\textbf{Step 3.} Let us now define the composite mapping $\Phi^k:=\Phi_{2}^k\circ\Phi_{1}^k $ as
\begin{align*}
\Phi^k:\quad C^\delta([0,T];\bm{Y}_{k}) \times C^\delta([0,T];Z_{k}) \ &\ \to H^1(0,T;\bm{Y}_{k})\times H^1(0,T;Z_{k}),
\\
(\bm{u}^k,\psi^k)\ &\ \mapsto(\bm{v}^{k},\varphi^k).
\end{align*}
From the compactness of $H^1(0,T;\bm{Y}_{k})$ into
$C^\delta([0,T];\bm{Y}_{k})$ and $H^1(0,T;Z_{k})$ into
$C^\delta([0,T];Z_{k})$ (recalling that $\bm{Y}_{k}$ and $Z_{k}$ are finite-dimensional), we find that $\Phi^k$ is a compact operator from $C^\delta([0,T];\bm{Y}_{k}) \times C^\delta([0,T];Z_{k})$ into itself. On the other hand, it follows from the continuous estimates \eqref{auphih1} and \eqref{aum} that
\be
\begin{aligned}
&\sup_{t\in [0,T]}\|\bm{v}_1^{k}(t)-\bm{v}^{k}_2(t)\| +
\sup_{t\in [0,T]}\|\varphi_1^{k}(t)-\varphi^{k}_2(t)\|
\\
&\quad \leq C_T\left(\sup_{t\in [0,T]}\|\bm{u}_1^{k}(t)-\bm{u}_2^{k}(t)\|
+\sup_{t\in [0,T]}\|\psi_1^{k}(t)-\psi_2^{k}(t)\|\right).
\end{aligned}
\notag
\ee
Thanks to the boundedness of $(\bm{v}^k_i,\varphi^k_i)$ in $H^1(0,T;\bm{Y}_{k})\times H^1(0,T;Z_{k})$, $i=1,2$, we conclude by interpolation that $\Phi^k$ is a continuous operator from $C^\delta([0,T];\bm{Y}_{k}) \times C^\delta([0,T];Z_{k})$ into itself.

Take
$$
\widetilde{M}= 2\mathrm{e}( C_0+ 1) +2C_{0,k} + 2(\|\bm{v}_0\|^2+  \|\varphi_0\|^2  +1),
$$
where the positive constants $C_0$, $C_{0,k}$ are given in  \eqref{iniC0} and \eqref{me2}, respectively.
According to the estimates \eqref{3eeenerg-a} and \eqref{auvm1}, there exists a sufficiently small time $T_*\in (0,T]$ depending on $\widetilde{M}$ such that
\begin{equation}
\|\bm{v}^{k}(t)\|^2+ \|\varphi^{k}(t)\|^2
\leq
2 \mathrm{e} ( C_0 + 1 ) + 2 C_{0,k}<\widetilde{M},\qquad \forall\, t\in [0,T_*].
\notag
\end{equation}
Define
\be
\begin{aligned}
\bm{B}_k&=\Big\{(\bm{u}^k,\psi^k)\in C^\delta([0,T_*];\bm{Y}_{k}) \times C^\delta([0,T_*];Z_{k})\ :\sup_{t\in[0,T_*]}\|\bm{u}^k(t)\|^2\leq \widetilde{M},\\
& \qquad  \sup_{t\in[0,T_*]}\|\psi^k(t)\|^2\leq \widetilde{M},\ \bm{u}^{k}(0)= \bm{P}_{\bm{Y}_{k}} \bm{v}_{0},\ \psi^{k}(0)= \bm{P}_{Z_{k}} \varphi_{0,\gamma}\Big\},
\end{aligned}
\notag
\ee
which is a closed convex set in $C^\delta([0,T_*];\bm{Y}_{k}) \times C^\delta([0,T_*];Z_{k})$. Then for any $(\bm{u}^k,\psi^k)\in \bm{B}_k$, we find that $(\bm{v}^k,\varphi^k)=\Phi^k(\bm{u}^k,\psi^k)\in H^1(0,T_*;\bm{Y}_{k})\times H^1(0,T_*;Z_{k})\subset\subset C^\delta([0,T_*];\bm{Y}_{k}) \times C^\delta([0,T_*];Z_{k})$ and the pair $(\bm{v}^k,\varphi^k)$  satisfies
\begin{align}
\sup_{t\in[0,T_{*}]}\|\bm{v}^k(t)\|^2 \le\widetilde{M},\quad \sup_{t\in[0,T_{*}]}\|\varphi^k(t)\|^2 \le\widetilde{M}.\notag
\end{align}
As a result, it holds $(\bm{v}^k,\varphi^k)\in \bm{B}_k$.

Let us recall the classical Schauder's fixed point theorem:
\bl\label{SH}
Assume that $\bm{K}$ is a closed convex set in a Banach space
$\mathcal{B}$. Let $\mathcal{T}$ be a continuous mapping of $\bm{K}$ into itself satisfying the image $\mathcal{T}\bm{K}$ is precompact. Then there exists a fixed point in $\bm{K}$ for $\mathcal{T}$.
\el
Then, we conclude from Lemma \ref{SH} that for the above chosen small time $T_*$, there exists a fixed point $(\bm{v}^k,\varphi^k)$ for the mapping $\Phi^k$ in the set $\bm{B}_k$ defined above. Furthermore, $\sigma^k$ can be determined by $(\bm{v}^k,\varphi^k)$ as in Lemma \ref{fp}. This yields a local solution to the semi-Galerkin scheme \eqref{atest.1}--\eqref{atest.ini0} on the interval $[0,T_*]$. Uniqueness of the approximate solution $(\bm{v}^k,\varphi^k,\mu^k,\sigma^k)$ can be proved by the standard energy method and we omit the details here.
The proof of Proposition \ref{p1} is complete.
\end{proof}

\section{Existence of Global Weak Solutions}
\setcounter{equation}{0}
\label{proof-main}
In this section, we prove Theorem \ref{main} on the existence of global weak solutions to the original problem \eqref{sysNew}--\eqref{ini0new}. The essential step is to derive uniform estimates for the approximate solutions with respect to the approximating parameters $k$, $\gamma$, $\epsilon$ and $n$. Then we pass to the limit first as $k \to +\infty$ and then as $\gamma, \epsilon \to 0$, $n\to +\infty$ in suitable topologies, using compactness methods.

\subsection{Uniform estimates}
\label{unie}

We begin with the estimates for approximate solutions $\left\{ (\bm{v}^k,\varphi^k, \mu^k\,\sigma^k) \right\}$ to the semi-Galerkin scheme \eqref{atest.1}--\eqref{atest.ini0}, which are uniform with respect to $k\geq \widehat{k}$ and $t\in [0,T]$.
\medskip

\textbf{Step 1. Mass dynamics.}
As in \cite{RSS}, testing \eqref{atest.2} with $\xi=1$ yields
\be
\frac{\mathrm{d}}{\mathrm{d}t} \overline{\varphi^{k}} +\alpha \overline{\varphi^{k}} =\frac{1}{|\Omega|}\int_\Omega h(\varphi^k,\sigma^k)\,\mathrm{d}x.
\notag
\ee
Owing to (H4), we have
$$
-h^*\leq \frac{\mathrm{d}}{\mathrm{d}t} \overline{\varphi^{k}} +\alpha \overline{\varphi^{k} }
\leq h^*\!,
$$
so that
\be
\overline{\varphi_{0,\gamma}}e^{-\alpha t}-(1-e^{-\alpha t})\frac{h^*}{\alpha}
\leq
\overline{\varphi^{k}}(t)
\leq
\overline{\varphi_{0,\gamma}}e^{-\alpha t}+(1-e^{-\alpha t})\frac{h^*}{\alpha},
\quad \forall\, t\in[0,T_k].
\notag
\ee
From the above inequality and (H4) we find that for any $k\geq \widehat{k}$, it holds
\begin{align*}
&\overline{\varphi^{k}}(t)\in \left[-\frac{h^*}{\alpha},\,\overline{\varphi_{0,\gamma}}\right],
\qquad \text{if}\ \ \overline{\varphi_{0,\gamma}}\geq \frac{h^*}{\alpha},
\\
&\overline{\varphi^{k}}(t)\in \left[-\frac{h^*}{\alpha},\,\frac{h^*}{\alpha}\right],
\qquad \ \  \text{if}\ \ \overline{\varphi_{0,\gamma}}\in \left(-\frac{h^*}{\alpha},\,\frac{h^*}{\alpha}\right),
\\
&\overline{\varphi^{k}}(t)\in \left[\overline{\varphi_{0,\gamma}},\,\frac{h^*}{\alpha}\right],
\qquad \ \ \ \,\text{if}\ \ \overline{\varphi_{0,\gamma}}\leq -\frac{h^*}{\alpha},
\end{align*}
In view of \eqref{Lvp1k}, we obtain
\begin{align}
|\overline{\varphi^{k}}(t)|\leq 1-\rho_*,\quad \forall\, t\in [0,T_k].
\label{averphi1}
\end{align}
 for some $\rho_*\in (0,1)$ depending only on $\overline{\varphi_{0}}$, $h^*$ and $\alpha$, but not on the approximating parameters $k$, $\gamma$, $\epsilon$, $n$.
When the mass source vanishes, i.e., $S(\varphi, \sigma)\equiv 0$,
the same property holds with $\rho_*=1-|\overline{\varphi_0}|$, because the spatial mean of $\varphi^k$ is indeed conserved.
\begin{remark}\label{rem-ks}
Since we are only interested in the case as $\gamma, \epsilon\to 0$,
in what follows, we always assume
$$
\gamma \in \left(0,\,\frac{\rho_*}{2}\right],\quad
\epsilon\in \left(0,\,\min\Big\{\epsilon_1,\,\frac{\rho_*}{2}\Big\}\right].
$$
\end{remark}

\textbf{Step 2. Basic energy law.}  Thanks to Proposition \ref{p1}, we are allowed to test \eqref{atest.1} by $\boldsymbol{v}^k$, \eqref{atest.2} by $\mu^k$, \eqref{atest.3} by $\partial_t\varphi^k$, and \eqref{atest.4} by $\ln \sigma^k-\chi\varphi^k$, respectively. Adding the resultants together, we obtain the following energy identity:
\begin{align}
& \frac{\mathrm{d}}{\mathrm{d}t}\mathcal{E}^k(t)
+ \mathcal{D}^k(t)
+ \kappa \int_{\Omega}  (\sigma^k)^2 \ln \sigma^k\,\mathrm{d}x  =   \mathcal{R}^k(t),\quad \forall\,t\in(0,T_k),
\label{menergy}
\end{align}
where
\begin{align*}
\mathcal{E}^k(t)
&=\int_{\Omega}\left( \frac{1}{2}|\boldsymbol{v}^k|^{2}
+ \frac{\gamma^{8}}{4} |\nabla \varphi^k|^4
+\frac{1}{2}|\nabla \varphi^k|^{2}
+ \Psi_k(\varphi^k) +\sigma^k(\ln \sigma^k-1)
-\chi \sigma^k\varphi^k\right) \, \mathrm{d} x,
\\
\mathcal{D}^k(t)
& =\int_{\Omega} \left(2 \eta(\varphi^k)|D \boldsymbol{v}^k|^{2}
+m(\varphi^k)|\nabla \mu^k|^{2}
+ \sigma^k|\nabla(\ln \sigma^k-\chi \varphi^k)|^{2}\right)\, \mathrm{d} x,
\\
\mathcal{R}^k(t)
&=\int_{\Omega} \big[-\alpha\varphi^k+h(\varphi^k,\sigma^k)\big]\mu^k\,\mathrm{d}x
+ \int_\Omega \beta(\varphi^k)\sigma^k\ln\sigma^k\,\mathrm{d}x
\\
&\quad
-\chi \int_{\Omega} \big[\beta(\varphi^k) \sigma^k-\kappa   (\sigma^k)^2\big]\varphi^k\,\mathrm{d}x.
\end{align*}
Noticing that
$$
\|\bm{P}_{Z_{k}} \varphi_{0,\gamma}\|_{W^{1,4}(\Omega)}\leq C \|\bm{P}_{Z_{k}} \varphi_{0,\gamma}\|_{H^2(\Omega)}\leq C \|\varphi_{0,\gamma}\|_{H^2(\Omega)},
$$
where $C>0$ only depends on $\Omega$, thanks to the construction of the initial data (see \eqref{Lvp3k}, \eqref{PZK}) and Lemma \ref{You}, we control the initial approximate energy as follows
\begin{align}
\mathcal{E}^k(0) &=
\int_{\Omega}\left(\frac{1}{2}|\bm{P}_{\bm{Y}_{k}} \bm{v}_{0}|^{2}
+ \frac{\gamma^{8}}{4} |\nabla \bm{P}_{Z_{k}} \varphi_{0,\gamma}|^4
+\frac{1}{2}|\nabla \bm{P}_{Z_{k}}\varphi_{0,\gamma}|^{2}
 \right) \mathrm{d}x
\notag \\
&\quad +
\int_{\Omega}\big(
\Psi_\epsilon(\bm{P}_{Z_{k}}\varphi_{0,\gamma})
+\sigma_{0,n}(\ln \sigma_{0,n}-1)
-\chi \sigma_{0,n}\bm{P}_{Z_{k}}\varphi_{0,\gamma}\big) \mathrm{d}x
\notag \\
&\leq \frac12 \|\bm{v}_{0}\|^{2}
 + C \big(\gamma^2 \|\varphi_{0,\gamma}\|_{H^2(\Omega)}\big)^4
+ \frac12\|\varphi_{0,\gamma}\|_{H^1(\Omega)}^2 + |\Omega|\max_{r\in[-1,1]}|\Psi_0(r)|
\notag \\
&\quad + 2(1+|\chi|) \int_\Omega \sigma_0\ln \sigma_0\, \mathrm{d}x + C (1+|\chi|) |\Omega|
\notag \\
&\le C\Big(\|  \bm{v}_{0}\|, \|\varphi_{0}\|_{H^1(\Omega)}, \max_{r\in[-1,1]}|\Psi_0(r)|,
\int_\Omega \sigma_0\ln \sigma_0\, \mathrm{d}x, \chi, \Omega\Big)
\notag\\
&=: \mathcal{E}_0,
\label{iniE0}
\end{align}
where the upper bound $\mathcal{E}_0$ is independent of the parameters $k$, $\gamma$, $\epsilon$, $n$.
\medskip

\textbf{Step 3. Lower bound of the approximate total energy $\mathcal{E}^k$.}
To obtain uniform estimates for the approximate solution $(\bm{v}^k,\varphi^k, \mu^k\,\sigma^k)$, it is crucial to show that the total energy $\mathcal{E}^k(t)$ is uniformly semi-coercive from below with respect to $k$, $\gamma$, $\epsilon$ and $n$. The key point is to deal with the crossing term $-\chi \sigma^k \varphi^k$.

Recalling the construction of $\Psi_{0,\epsilon}$, we find
\begin{align}
\Psi_{0,\epsilon}(r)
& = \int_0^{1-\epsilon} \Psi'_{0,\epsilon}(s)\,\mathrm{d}s
+ \int^2_{1-\epsilon} \Psi'_{0,\epsilon}(s)\,\mathrm{d}s
+ \int_2^r \Psi'_{0,\epsilon}(s)\,\mathrm{d}s
\notag \\
&= \Psi_{0}(1-\epsilon)
+  (r-1+\epsilon) \Psi_{0}'(1-\epsilon)
\notag \\
&\quad
+ \left[\frac12 (1+\epsilon)^2 + \left(\frac{4|\chi|+3}{4(|\chi|+1)}+\epsilon\right)(r-2)\right] \Psi_{0}''(1-\epsilon)
\notag \\
&\quad
+ \frac{1}{16(|\chi|+1)^2}\Psi_{0}''(1-\epsilon) \left(\mathrm{e}^{4(|\chi|+1)r-8(|\chi|+1)}-1\right)
\notag \\
&\geq
\frac{\theta}{16(|\chi|+1)^2} \left(\mathrm{e}^{4(|\chi|+1)r-8(|\chi|+1)}-1\right),
\qquad \forall\, r\geq 2.
\label{PPsi}
\end{align}
Hence, there exists $r^*>2$ depending on $\theta$, $\theta_0$ and $\chi$, but not on $k$, such that for any $k\geq \widehat{k}$, it holds
\be
\begin{aligned}
\Psi_{0,\epsilon}(r)
\ge   \mathrm{e}^{(2|\chi|+1)r}
+ \theta_0  r^2 + 2|\chi|r+1,
\qquad \forall\, r\in[r^*,+\infty).
\label{low-2}
\end{aligned}
\ee	
We claim that
\be
\begin{aligned}
	&\frac{1}{2} \Psi_{0,\epsilon}(\varphi^k) -\frac{\theta_0}{2} (\varphi^k)^2 +\frac12\sigma^k(\ln \sigma^k-1)-|\chi| \sigma^k\varphi^k\geq 0,
\end{aligned}
\label{neg-1}
\ee	
for all $\sigma^k\geq 0$ and $\varphi^k\ge r^*$.
Indeed, if $\sigma^k\in[0,1]$ and $\varphi^k\ge r^*$, \eqref{neg-1} is a direct consequence of \eqref{low-2}. If $\sigma^k\ge 1$ and $\varphi^k\ge r^*$, it follows from \eqref{low-2} that
\begin{align}
	&\frac{1}{2} \Psi_{0,\epsilon}(\varphi^k) -\frac{\theta_0}{2} (\varphi^k)^2 +\frac12\sigma^k(\ln \sigma^k-1)-|\chi| \sigma^k\varphi^k
\notag \\
	&\quad\ge \frac12 \mathrm{e}^{(2|\chi|+1)\varphi^k} +\frac12 \sigma^k(\ln\sigma^k-1)-|\chi| (\sigma^k-1)\varphi^k.
\label{low-3}
\end{align}
Taking $a=2|\chi|\varphi^k$ and $b=\sigma^k-1$ in the generalized Young inequality \eqref{general young}, we get
\be
|\chi| (\sigma^k-1)\varphi^k \leq \frac12  \mathrm{e}^{2|\chi|\varphi^k} +\frac12 \sigma^k(\ln\sigma^k-1),
\label{low-1}
\ee
which together with \eqref{low-3} yields the claimed result \eqref{neg-1}.
By a similar procedure, we can find some constant $r_*<-2$ depending on $\theta$, $\theta_0$ and $\chi$, but not on $k$, such that for any $k\geq \widehat{k}$, it holds
\be
\begin{aligned}
	&\frac{1}{2} \Psi_{0,k}(\varphi^k) -\frac{\theta_0}{2} (\varphi^k)^2 +\frac12\sigma^k(\ln \sigma^k-1)+|\chi| \sigma^k\varphi^k\geq 0,
\end{aligned}
\label{neg-1a}
\ee	
for all $\sigma^k\geq 0$ and $\varphi^k\le r_*$.
Finally, for $\sigma^k\geq 0$ and $\varphi^k\in [r_*, r^*]$, it follows from \eqref{general young} that
\begin{align}
|\chi\sigma^k\varphi^k|
&\leq (1+|\chi|)\max\{-r_*,r^*\}\sigma^k
\notag \\
&\leq (1+|\chi|)\max\{-r_*,r^*\}
+ \frac12\mathrm{e}^{2(1+|\chi|)\max\{-r_*,r^*\}}
+ \frac12\sigma(\ln\sigma^k-1).
\label{low-4}
\end{align}
In conclusion, there exists a sufficiently large constant $C_*>0$ depending on $\theta$, $\theta_0$, $\chi$ and $\Omega$, but not on $k$, $\gamma$, $\epsilon$ and $n$, such that
\begin{align}
&|\chi|\int_{\Omega}  |\sigma^k \varphi^k |\, \mathrm{d}x
 \leq \frac{1}{2} \int_{\Omega}\left[ \Psi_{\epsilon}(\varphi^k)+ \sigma^k(\ln \sigma^k-1)\right]\,\mathrm{d}x+C_*,
\label{Lowerbd-1}
\end{align}
which yields
\begin{align}
&\mathcal{E}^k(t)\geq \frac12 \int_{\Omega}\left[ |\boldsymbol{v}^k|^{2}
+ |\nabla \varphi^k|^{2}
+ \Psi_\epsilon(\varphi^k) +\sigma^k(\ln \sigma^k-1)
\right] \mathrm{d} x +   \frac{\gamma^{8}}{4} \int_{\Omega} |\nabla \varphi^k|^4\, \mathrm{d} x  -C_*,
\label{Lowerbd-2}
\end{align}
where the constant $C_*>0$ is independent of $k$,  $\gamma$, $\epsilon$, $n$.
\medskip

\textbf{Step 4. Estimate of the remainder term $\mathcal{R}^k$.}
Similar to \cite[(3.8)]{RSS}, we deduce from (H4) and the Poincar\'{e}--Wirtinger inequality that
\begin{align}
\int_{\Omega} \big(-\alpha\varphi^k +h(\varphi^k,\sigma^k)\big)\mu^k\,\mathrm{d}x
&\leq C(\alpha+ h^*)\|\nabla \mu^k\|
+ (\alpha+h^*) \left|\int_\Omega \mu^k\,\mathrm{d}x\right|
\notag\\
&\leq C(\alpha+ h^*)\|\nabla \mu^k\| + (\alpha+h^*)\|\Psi_{0,\epsilon}'(\varphi^k)\|_{L^1(\Omega)}
\notag\\
&\quad +(\alpha+h^*)\theta_0|\Omega| |\overline{\varphi^k }| + |\chi|\int_\Omega \sigma^k\,\mathrm{d}x,
\label{es-mass1}
\end{align}
where $C>0$ only depends on $\Omega$. The last term on the right-hand side of \eqref{es-mass1} can be estimated as in \eqref{low-4}.
Recalling the construction of $\Psi_{0,\epsilon}'$, we have
$\Psi_{0,\epsilon}'(r)=\Psi_0'(r)$ for $r\in [-1+\epsilon,1-\epsilon]$. Hence, thanks to \cite[Proposition A.1]{MZ04} and the mass relation \eqref{averphi1}, we deduce that
$$
|\Psi_{0,\epsilon}'(r)|\leq c_1 \Psi_{0,\epsilon}'(r)\big( r - \overline{\varphi^k}\big)+c_2,
\qquad \forall\, r\in [-1+\epsilon,\,1-\epsilon],
$$
where the positive constants $c_1$, $c_2$ depend on $\rho_*$, but not on $k$, $\gamma$, $\epsilon$ and $n$.
Besides, for all $k\geq \widehat{k}$, we have (cf. Remark \ref{rem-ks})
$$
\Psi_{0,\epsilon}'(r)\big(r-\overline{\varphi^k}\big)\geq \left(1-\frac12\rho_*-\overline{\varphi^k}\right)\Psi_{0,\epsilon}'(r) \geq \frac12\rho_*\Psi_{0,\epsilon}'(r),\quad \forall\, r\geq 1-\epsilon.
$$
A similar reasoning yields
$$
\Psi_{0,\epsilon}'(r)\big(r-\overline{\varphi^k}\big)\geq \left(-1+\frac12\rho_*-\overline{\varphi^k}\right)\Psi_{0,\epsilon}'(r)\geq -\frac12\rho_*\Psi_{0,\epsilon}'(r),\quad \forall\, r\leq -1+\epsilon.
$$
Thus, there exist two positive constants $\widetilde{c}_1$, $\widetilde{c}_2$, both depending on $\rho_*$, but not on $k$ $\gamma$, $\epsilon$ and $n$ such that
\begin{align}
& \|\Psi_{0,\epsilon}'(\varphi^k)\|_{L^1(\Omega)}\leq \widetilde{c}_1 \int_\Omega \Psi_{0,\epsilon}'(\varphi^k)\big(\varphi^k-\overline{\varphi^k}\big)\,\mathrm{d}x +\widetilde{c}_2.
\label{Psi-L1}
\end{align}
Testing \eqref{atest.3} by $\xi=\varphi^k-\overline{\varphi^k}$, we get
\begin{align}
&  \gamma^{8}  \|\nabla \varphi^k\|_{\bm{L}^4(\Omega)}^4
+ \|\nabla \varphi^k\|^2+ \int_\Omega \Psi_{0,\epsilon}'(\varphi^k) (\varphi^k-\overline{\varphi^k})\,\mathrm{d}x
\notag \\
&\quad = (\mu^k-\overline{\mu^k},\varphi^k-\overline{\varphi^k})
+\theta_0\|\varphi^k-\overline{\varphi^k}\|^2
+ \chi\int_\Omega \sigma^k(\varphi^k-\overline{\varphi^k})\,\mathrm{d}x
\notag \\
&\quad \leq C\|\nabla \mu^k\|\|\nabla \varphi^k\|+ C\theta_0\|\nabla \varphi^k\|^2
+ \frac{1}{2} \int_\Omega \Psi_{\epsilon}(\varphi^k)\,\mathrm{d}x
 +\int_\Omega \sigma^k(\ln \sigma^k-1)\,\mathrm{d}x +C,
\label{Psi-L2}
\end{align}
where in the last line we have used estimates similar to \eqref{neg-1}, \eqref{neg-1a}, \eqref{low-4} and the Poincar\'{e}--Wirtinger inequality. By Young's inequality, we conclude from \eqref{es-mass1}, \eqref{Psi-L1} and \eqref{Psi-L2} that
\begin{align}
&\int_{\Omega} \big(-\alpha\varphi^k +h(\varphi^k,\sigma^k)\big)\mu^k\,\mathrm{d}x
\notag \\
& \quad \leq \frac{m_*}{2}\|\nabla \mu^k\|^2 + C\int_{\Omega}\left[\frac{1}{2}|\nabla \varphi^k|^{2}
+ \Psi_{\epsilon}(\varphi^k) +\sigma^k(\ln \sigma^k-1)
\right]\, \mathrm{d} x +C,
\label{es-mass2}
\end{align}
where the constant $C>0$ is independent of $k$,  $\gamma$, $\epsilon$ and $n$.

Next, thanks to (H5), we can apply the argument used in Step 3 to conclude that
\begin{align}
&\int_{\Omega} \beta(\varphi^k) \sigma^k \ln \sigma^k\,\mathrm{d}x
 \leq  2b^* \int_\Omega \sigma^k(\ln \sigma^k-1)\,\mathrm{d}x +C,
\label{es-react1}
\end{align}
and
\begin{align}
-\chi \int_{\Omega}  \beta(\varphi^k) \sigma^k \varphi^k\,\mathrm{d}x
&\le b^*|\chi|\int_{\Omega} | \sigma^k\varphi^k|\,\mathrm{d}x
\notag \\
& \leq \frac{b^*}{2} \int_{\Omega}\left[ \Psi_{\epsilon}(\varphi^k)+ \sigma^k(\ln \sigma^k-1)\right]\,\mathrm{d}x+C,
\label{es-react2}
\end{align}
where the constant $C>0$ is independence of $k$, $\gamma$, $\epsilon$ and $n$.

It remains to estimate the last term
$\chi \kappa \int_{\Omega}     (\sigma^k)^2  \varphi^k\,\mathrm{d}x$.
We apply the generalized Young's inequality again. Indeed, it follows from \eqref{PPsi} that there exists some $\widetilde{r^*}>2$ depending on $\theta$, $\theta_0$ and $\chi$, but not on $k$, $\gamma$, $\epsilon$ and $n$ such that for any $k\geq \widehat{k}$, it holds
\be
\begin{aligned}
\Psi_{0,\epsilon}(r)
\ge   \mathrm{e}^{(4|\chi|+1)r}
+ \theta_0  r^2 + 4|\chi|r,
\qquad \forall\, r\in \left[\widetilde{r^*},+\infty \right).
\label{low-2b}
\end{aligned}
\ee	
Then for $\sigma^k\geq 1$ and $\varphi^k\geq \widetilde{r^*}$, we find
\be
\begin{aligned}
	&\frac{1}{4} \Psi_{0,\epsilon}(\varphi^k) -\frac{\theta_0}{4} (\varphi^k)^2 +\frac12(\sigma^k)^2 \ln \sigma^k -|\chi| (\sigma^k)^2\varphi^k
\\
	&\quad\ge \frac14 \mathrm{e}^{(4|\chi|+1)\varphi^k} +\frac12 (\sigma^k)^2 \ln\sigma^k -|\chi| \big[(\sigma^k)^2-1\big]\varphi^k.
\end{aligned}
\label{low-3b}
\ee	
Taking now $a=4|\chi|\varphi^k$ and $b=(\sigma^k)^2-1$ in \eqref{general young}, we get
\be
|\chi| \big[(\sigma^k)^2-1\big]\varphi^k \leq \frac14  \mathrm{e}^{4|\chi|\varphi^k} +\frac12 (\sigma^k)^2\ln\sigma^k,
\label{low-1b}
\ee
which together with \eqref{low-3b} yields
\be
\begin{aligned}
	&\frac{1}{4} \Psi_{0,\epsilon}(\varphi^k) -\frac{\theta_0}{4} (\varphi^k)^2 +\frac12(\sigma^k)^2\ln \sigma^k -|\chi| (\sigma^k)^2\varphi^k\geq 0.
\end{aligned}
\label{neg-1b}
\ee	
Applying a similar argument  for \eqref{Lowerbd-1}, we obtain the following estimate:
\begin{align}
	& \chi\kappa \int_{\Omega}  (\sigma^k)^2  \varphi^k\,\mathrm{d}x
 \le  \frac{\kappa}{2}\int_{\Omega}  (\sigma^k)^2 \ln \sigma^k \,\mathrm{d}x + \frac{\kappa}{4} \int_\Omega \Psi_{\epsilon}(\varphi^k)
 \,\mathrm{d}x +\widetilde{C_*},
 \label{es-react3}
\end{align}
where the constant $\widetilde{C_*}>0$ is independent of $k$, $\gamma$, $\epsilon$ and $n$.
\medskip

\textbf{Step 5. Basic energy  estimates.}
Combining the estimates \eqref{es-mass2}, \eqref{es-react1}, \eqref{es-react2} and \eqref{es-react3}, we deduce from \eqref{menergy} and (H2), (H3) that
\begin{align}
& \frac{\mathrm{d}}{\mathrm{d}t}\mathcal{E}^k(t)
+ \int_{\Omega} \left(2 \eta_*|D \boldsymbol{v}^k|^{2}
+\frac{m_*}{2}|\nabla \mu^k|^{2}
+ \sigma^k|\nabla(\ln \sigma^k-\chi \varphi^k)|^{2} + \frac{\kappa}{2} (\sigma^k)^2\ln \sigma^k\right)\, \mathrm{d} x
\notag \\
&\quad \le C\mathcal{E}^k(t) +C',\quad \forall\,t\in(0,T_k),
\notag
\end{align}
where the constants $C, C'>0$ are independent of $k\geq \widehat{k}$, $\gamma$, $\epsilon$ and $n$. Applying Gronwall's lemma, we get
\begin{align}
&\mathcal{E}^k(t) \leq \left(\mathcal{E}^k(0)+\frac{C'}{C}\right) \mathrm{e}^{Ct}
\leq \left(\mathcal{E}_0+\frac{C'}{C}\right) \mathrm{e}^{Ct},
 \quad \forall\,t\in[0,T_k],
 \notag
\end{align}
and
\begin{align}
& \int_0^t \int_{\Omega} \left(2 \eta_*|D \boldsymbol{v}^k|^{2}
+\frac{m_*}{2}|\nabla \mu^k|^{2}
+ \sigma^k|\nabla(\ln \sigma^k-\chi \varphi^k)|^{2}+ \frac{\kappa}{2}(\sigma^k)^2\ln \sigma^k \right)\, \mathrm{d} x \mathrm{d}s
\notag \\
&\quad
  \leq \mathcal{E}_0+ \left(\mathcal{E}_0+\frac{C'}{C}\right) \mathrm{e}^{Ct}+C't, \quad \forall\,t\in[0,T_k].
  \notag
\end{align}
From the semi-coercivity of the energy $\mathcal{E}^k$ (see  \eqref{Lowerbd-2}), we  deduce the following uniform bounds
\begin{align}
&\|\bm{v}^k\|_{L^{\infty}(0, T_k ; \bm{L}^2(\Omega))}
 +\gamma^2\|\varphi^k\|_{L^{\infty}(0, T_k ; W^{1,4}(\Omega))}
+\|\varphi^k\|_{L^{\infty}(0, T_k ; H^1(\Omega))}
\notag \\
&\quad +\|\sigma^k \ln \sigma^k\|_{L^{\infty}(0, T_k ; L^{1}(\Omega))} +\|\Psi_{0,\epsilon}(\varphi^k)\|_{L^{\infty}(0, T_k; L^{1}(\Omega))}
\leq C,
\label{energy-a}
\end{align}
and
\begin{align}
&\|\bm{v}^k\|_{L^{2}(0, T_k ; \bm{H}^1(\Omega))}
+\|\nabla \mu^k\|_{L^{2}(0, T_k ; \bm{L}^{2}(\Omega))}
+\kappa \|(\sigma^k)^{2} \ln \sigma^k\|_{L^{1}(0, T_k; L^{1}(\Omega))}
\notag \\
&\quad+\| (\sigma^k)^\frac{1}{2}\nabla(\ln \sigma^k-\chi\varphi^k)\|_{L^{2}(0, T_k; \bm{L}^{2}(\Omega))}
\leq C,
\label{energy-b}
\end{align}
where the constant $C>0$ depends on $\mathcal{E}_0$, $\Omega$ and coefficients of the system,  but it is independent of $k$, $\gamma$, $\epsilon$ and $n$.

Testing \eqref{atest.3} by $\xi=1$ yields
\begin{align}
|\overline{\mu^{k}}|
&=|\Omega|^{-1}|(\Psi'_\epsilon(\varphi^{k}),1)-\chi(\sigma^{k},1)|
\notag \\
 & \le C\| \Psi'_\epsilon(\varphi^{k})\|_{L^1(\Omega)}
 +C|\chi|\int_\Omega \sigma^{k}(\ln \sigma^k-1)\,\mathrm{d}x+C|\chi|.
\label{average valuea}
\end{align}
Then from \eqref{Psi-L1}, \eqref{Psi-L2}, \eqref{energy-a}, \eqref{energy-b}, \eqref{average valuea} and the Poincar\'{e}--Wirtinger inequality, we conclude
that %
\begin{align}
\|\Psi'_\epsilon(\varphi^{k})\|_{L^2(0,T_k;L^1(\Omega))}\leq C,
\label{Psid-1}
\end{align}
as well as
\begin{equation}
\|  \mu^{k} \|_{L^{2}(0,T_k;H^1(\Omega))}\le C,
\label{mu}
\end{equation}
where the constant $C>0$ is independent of $k$, $\gamma$, $\epsilon$ and $n$.

\begin{remark}\label{rem:gloapp}
For any $k\geq \widehat{k}$, the estimates \eqref{energy-a}, \eqref{energy-b} permit to extend the local solution 
$(\bm{v}^k,\varphi^k,\mu^k,\sigma^k)$ 
from $[0,T_k]$ to the whole interval $[0,T]$ (recall that $T\in (0,+\infty)$ is an arbitrary final time). This yields a unique global solution at the approximate level. Moreover, the estimates \eqref{energy-a}, \eqref{energy-b} and \eqref{mu} are valid on $[0,T]$, where the constant $C>0$ may depend on $T$, but not on $k$, $\gamma$, $\epsilon$ and $n$.
\end{remark}

\textbf{Step 6. Further estimates on $\varphi^k$ and $\sigma^k$.}
In light of Remark \ref{rem:gloapp}, hereafter we just work with the global approximate solution $(\bm{v}^k,\varphi^k,\mu^k,\sigma^k)$ on $[0,T]$.

We observe that the estimates \eqref{energy-a}, \eqref{energy-b} and \eqref{mu} are independent of the parameter $\kappa$ for the logistic degradation. Below we derive some further estimates that rely on $\kappa>0$.

First, it follows from \eqref{energy-b} and the generalized Young's inequality \eqref{general young} that
\be
\|\sigma^k\|_{L^{2}(0, T; L^{2}(\Omega))} \leq C,
\label{energy-c}
\ee
where $C>0$ depends on the size of $\kappa^{-1}$, but is independent of $k$, $\gamma$, $\epsilon$ and $n$.

Next, testing \eqref{atest.3} by $-\Delta \varphi^k$, we obtain
\begin{align}
& \|\Delta \varphi^k \|^2 + \gamma^{8}  \int_\Omega \mathrm{div}( |\nabla \varphi^k|^2 \nabla \varphi^k)\Delta \varphi^k\,\mathrm{d}x
- \int_\Omega  \Psi'_{\epsilon}(\varphi^k)  \Delta \varphi^k\,\mathrm{d}x
\notag \\
&\quad =\int_\Omega (\mu^k+\chi \sigma^k)\Delta \varphi^k\,\mathrm{d}x
\notag\\
&\quad \leq \frac12 \|\Delta \varphi^k \|^2 + \|\mu^k\|^2+ \chi^2 \|\sigma^k\|^2.
\label{phikH2-1}
\end{align}
Using integration by parts and the fact that $\Psi_{0,\epsilon}''\geq \theta$, we find
\begin{align}
- \int_\Omega  \Psi'_{\epsilon}(\varphi^k)  \Delta \varphi^k\,\mathrm{d}x
& = \int_\Omega  \Psi''_{0,\epsilon}(\varphi^k) |\nabla \varphi^k|^2 \,\mathrm{d}x - \theta_0\|\nabla \varphi^k\|^2
\notag \\
&\geq (\theta-\theta_0)\|\nabla \varphi^k\|^2.
\label{phikH2-2}
\end{align}
A direct calculation yields
\begin{align}
& \int_\Omega \mathrm{div}( |\nabla \varphi^k|^2 \nabla \varphi^k)\Delta \varphi^k\,\mathrm{d}x
\notag \\
&\quad = - \int_\Omega \partial_j \partial_i( |\nabla \varphi^k|^2 \partial_i \varphi^k) \partial_j \varphi^k\,\mathrm{d}x
+\underbrace{\int_{\partial\Omega} \partial_i( |\nabla \varphi^k|^2 \partial_i \varphi^k) \partial_{\bm{n}}\varphi^k \,\mathrm{d}S}_{=0}
\notag\\
&\quad = \int_\Omega \partial_j ( |\nabla \varphi^k|^2 \partial_i \varphi^k) \partial_i \partial_j \varphi^k\,\mathrm{d}x
- \int_{\partial\Omega} \partial_j( |\nabla \varphi^k|^2 \partial_i \varphi^k) \bm{n}_i \partial_j \varphi^k\,\mathrm{d}x
\notag \\
&\quad = \underbrace{\int_\Omega |\nabla \varphi^k|^2 (\partial_j \partial_i \varphi^k)^2 \,\mathrm{d}x}_{\geq 0}
+  \frac{1}{2}
 \big\|\nabla |\nabla \varphi^k|^2\big\|^2
\notag \\
&\qquad - \underbrace{\int_{\partial\Omega} \partial_j  |\nabla \varphi^k|^2 (\partial_i \varphi^k \bm{n}_i) \partial_j\varphi^k\,\mathrm{d}x}_{=0}
- \int_{\partial\Omega}  |\nabla \varphi^k|^2 \partial_j \partial_i \varphi^k  \bm{n}_i \partial_j \varphi^k\,\mathrm{d}x,
\label{phikH2-3}
\end{align}
where $\partial_i$ denotes the partial derivative $\partial_{x_i}$ and $\bm{n}_i$ denotes the $i$-th component of the outer normal $\bm{n}$. Hereafter, we use Einstein's summation convention for repeated indices.
It remains to control the last term on the right-hand side of \eqref{phikH2-3}. According to \cite[(4.17)]{CM14} (which is a consequence of \cite[Equation (3, 1, 1, 8)]{Gr85}), we have
\begin{align}
& \Delta \varphi^k\partial_{\bm{n}}\varphi^k -\partial_j\partial_i\varphi^k\partial_j\varphi^k\bm{n}_i
\notag \\
&\quad = \mathrm{div}_\Gamma(\partial_{\bm{n}}\varphi^k \nabla_\Gamma \varphi^k)-\mathrm{tr}\Pi(\partial_{\bm{n}}\varphi^k)^2
-\Pi(\nabla_\Gamma \varphi^k,\nabla_\Gamma \varphi^k)
-2\nabla_\Gamma \varphi^k\cdot \nabla_\Gamma (\partial_{\bm{n}}\varphi^k)
\label{CM}
\end{align}
on $\partial\Omega$,  where $\Pi$ is the second fundamental form on $\partial \Omega$ and $\mathrm{tr}\Pi$ denotes its trace.
Since $\partial_{\bm{n}}\varphi^k=0$ on $\partial\Omega$, the identity \eqref{CM} yields
\begin{align}
\int_{\partial\Omega}  |\nabla \varphi^k|^2 \partial_j \partial_i \varphi^k  \bm{n}_i \partial_j \varphi^k\,\mathrm{d}x
=\int_{\partial\Omega} |\nabla_\Gamma \varphi^k|^2 \Pi(\nabla_\Gamma \varphi^k,\nabla_\Gamma \varphi^k)\,\mathrm{d}x.
\label{phikH2-4}
\end{align}
Here we recall that for a scalar differentiable function $f: \partial\Omega \to \mathbb{R}$, its surface (tangential) gradient is defined as
$\nabla_\Gamma f = \nabla \widetilde{f} - (\nabla \widetilde{f}\cdot \bm{n})\bm{n}$, where $\nabla$ denotes the usual gradient operator in $\mathbb{R}^d$, $\widetilde{f}$ is a smooth extension of $f$ to a $d$-dimensional neighbourhood $U$ of the surface $\partial\Omega$ such that $\widetilde{f}|_\Gamma =f$. The surface gradient $\nabla_\Gamma f$ depends only on
the values of $f$ on $\partial\Omega$ and its components are denoted by $\nabla_\Gamma f= (\underline{D}_if)_{i=1}^d$.
Thanks to the following trace inequality (see e.g., \cite[Chapter 2, (2.27)]{La}):
\begin{align}
\|f\|_{L^2(\partial\Omega)}\leq C_{\mathrm{tr}}\|f\|_{H^1(\Omega)}^\frac12\|f\|_{L^2(\Omega)}^\frac12,\qquad \forall\, f\in H^1(\Omega),
\notag
\end{align}
where $C_{\mathrm{tr}}>0$ only depends on $\Omega$ and $\partial\Omega$, we can estimate \eqref{phikH2-4}  by
\begin{align}
& \left|\int_{\partial\Omega} |\nabla_\Gamma\varphi^k|^2 \Pi(\nabla_\Gamma \varphi^k,\nabla_\Gamma \varphi^k)\,\mathrm{d}x\right|
\notag \\
&\quad \leq \|\Pi\|_{L^\infty(\partial\Omega)}
\int_{\partial\Omega} |\nabla_\Gamma \varphi^k|^4\,\mathrm{d}x
\notag \\
&\quad \leq \|\Pi\|_{L^\infty(\partial\Omega)}
\big\||\nabla\varphi^k|^2\big\|_{L^2(\partial\Omega)}^2
\notag \\
&\quad \leq
C\|\Pi\|_{L^\infty(\partial\Omega)}\big\||\nabla\varphi^k|^2 \big\|_{H^1(\Omega)} \big\||\nabla\varphi^k|^2\big\|
\notag \\
&\quad \leq
C\|\Pi\|_{L^\infty(\partial\Omega)}\big\|\nabla |\nabla\varphi^k|^2\big\|  \big\||\nabla\varphi^k|^2\big\| + C\|\Pi\|_{L^\infty(\partial\Omega)}\big\||\nabla\varphi^k|^2\big\|^2
\notag \\
&\quad \leq \frac14 \big\|\nabla |\nabla\varphi^k|^2\big\|^2
+ C (\|\Pi\|_{L^\infty(\partial\Omega)}^2+1) \|\nabla\varphi^k \|_{\bm{L}^4(\Omega)}^4.
\notag
\end{align}
Combining the above estimates, we deduce from \eqref{phikH2-1} that
\begin{align}
& \frac12 \|\Delta \varphi^k \|^2 + \frac{\gamma^{8}}{4} \big\|\nabla |\nabla\varphi^k|^2\big\|^2
\notag \\
&\quad \leq  \|\mu^k\|^2+ \chi^2 \|\sigma^k\|^2 + \theta_0\|\nabla \varphi^k\|^2
+ C (\|\Pi\|_{L^\infty(\partial\Omega)}^2+1) \big(\gamma^2\|\nabla\varphi^k \|_{\bm{L}^4(\Omega)}\big)^4.
\label{phikH2-6}
\end{align}
Since $\partial\Omega$ is sufficiently smooth, e.g., $\partial \Omega$ is of class $C^2$,  $\|\Pi\|_{L^\infty(\partial\Omega)}$ is bounded.
Then it follows from \eqref{energy-a}, \eqref{mu}, \eqref{energy-c}, \eqref{phikH2-6} and the elliptic estimate that
\begin{align}
& \|\varphi^k\|_{L^{2}(0, T; H^2(\Omega))}\le C,
\label{phikH2-7}
\end{align}
and
\begin{align}
& \gamma^4\big\||\nabla \varphi^k|^2\big\|_{L^2(0,T;H^1(\Omega))}\le C,
\label{pL-L2H1}
\end{align}
where $C>0$ is independent of $k$, $\gamma$, $\epsilon$ and $n$. On the other hand, \eqref{energy-a} implies
 \be
\|\varphi^k\|_{L^{\infty}(0, T; W^{1,4}(\Omega))}\le \frac{C}{\gamma^2},
\label{tm}
\ee
where $C>0$ is independent of $k$, $\gamma$, $\epsilon$ and $n$.
The $\gamma$-dependent estimate \eqref{tm} allows us to derive an $L^2_tL^2_x$-estimate for $\Psi'_{\epsilon}(\varphi^k)$.
From the definition of $\Psi'_{\epsilon}$ and the Sobolev embedding $W^{1,4}(\Omega)\hookrightarrow L^\infty(\Omega)$ in two and three dimensions, we deduce from \eqref{tm} that
\begin{align}
	\|\Psi'_{\epsilon}(\varphi^k)\|_{L^2(0,T;L^2(\Omega))}\leq C,
	\label{Psip0}
\end{align}
where $C>0$ depends on $\gamma$, $\epsilon$, $\mathcal{E}_0$, $\Omega$ and coefficients of the system, but it is independent of $k$ and $n$.

Next, testing \eqref{atest.4} by $\ln \sigma^k$, using integration by parts, H\"{o}lder's inequality and the generalized Young's inequality \eqref{general young}, we obtain
\begin{align}
&\frac{\mathrm{d}}{\mathrm{d} t} \int_{\Omega} \sigma^k(\ln \sigma^k-1)\,\mathrm{d}x
+  \int_{\Omega} \left(4 |\nabla (\sigma^k)^{\frac12}|^{2}+ \kappa  (\sigma^k)^{2} \ln \sigma^k\right) \mathrm{d}x
\notag \\
&\quad = \int_{\Omega} \left(-  \chi  \sigma^k  \Delta \varphi^k +  \beta(\varphi^k) \sigma^k \ln \sigma^k\right) \mathrm{d}x
\notag \\
&\quad \leq   \frac{\kappa}{2} \int_{\Omega} (\sigma^k)^{2} \ln \sigma^k\,\mathrm{d}x +\chi^2 \|\Delta\varphi^k\|^2
+ 2b^*\int_\Omega \sigma^k(\ln \sigma^k-1)\,\mathrm{d}x
  +C,
\label{sig-es-a}
\end{align}
where $C>0$ depends on $\kappa$, but is independent on $k$, $\gamma$, $\epsilon$, $n$.
Integrating \eqref{sig-es-a} in time, using  \eqref{energy-a} and  \eqref{phikH2-7}, we find
\begin{align}
 \|(\sigma^k)^{\frac12}\|_{L^{2}(0, T ; H^1(\Omega))} \leq C,
\label{sig-es-aa}
\end{align}
where $C>0$ depends on $\kappa$, but it is independent on $k$, $\gamma$, $\epsilon$, $n$.
Besides, by an argument similar to that for \cite[Lemma 2.11]{Lan2016}, we get
\begin{equation}
 \int_{\Omega} |\nabla \sigma^k |^{\frac{4}{3}}\,\mathrm{d}x= \int_{\Omega} \frac{ |\nabla \sigma^k|^{\frac{4}{3}}}{(\sigma^k)^{\frac{2}{3}}} (\sigma^k)^{\frac{2}{3}}\,\mathrm{d}x
 \leq \frac{2}{3}  \int_{\Omega} \frac{|\nabla \sigma^k|^2}{\sigma^k}\,\mathrm{d}x
 +\frac{1}{3} \ \int_{\Omega} |\sigma^k|^2\,\mathrm{d}x,
 \notag
\end{equation}
which together with \eqref{energy-c}, \eqref{sig-es-aa} leads to the following estimate
\be
\|\sigma^k\|_{L^{\frac{4}{3}}(0, T ; W^{1,\frac{4}{3}}(\Omega))} \leq C,
\label{sigma}
\ee
where $C>0$ depends on $\kappa$, but is independent on $k$, $\gamma$, $\epsilon$, $n$.

\begin{remark}\label{rem-TM}
We observe that the $p$-Laplace regularization in \eqref{reg.3} yields the $k$-independent estimate \eqref{Psip0}, while it brings extra difficulty for the derivation of \eqref{phikH2-7} in the case of a general bounded smooth domain. When $d=2$, this $p$-Laplace regularization is actually unnecessary. Setting $\gamma=0$, we can easily obtain \eqref{phikH2-7} from \eqref{phikH2-1}, without handling \eqref{phikH2-3} that does not have a definite sign. In order to recover \eqref{Psip0}, we recall the following Trudinger--Moser inequality in two dimensions, which is a direct consequence of \cite[Theorem 2.2]{NSY}: there exist positive constants $c_1, c_2$ such that
\begin{align}
\int_\Omega e^{|u|}\,\mathrm{d}x\leq c_1\exp\left\{
\left(\frac{1}{8\pi}+1\right)\|\nabla u\|^2+ c_2\|u\|_{L^1(\Omega)}^2\right\},
\quad \forall\, u\in H^1(\Omega).
\label{tm-2d}
\end{align}
Then by the definition of $\Psi'_{\epsilon}$ and
\eqref{energy-a}, we can apply \eqref{tm-2d} to conclude
\begin{align}
	\|\Psi'_{\epsilon}(\varphi^k)\|_{L^2(0,T;L^2(\Omega))}\leq C,
	\label{Psip0b}
\end{align}
where the constant $C>0$ depends on $\epsilon$, $\mathcal{E}_0$, $\Omega$ and coefficients of the system, but it is independent of $k$ and $n$.
Finally, let us emphasize that both estimates \eqref{Psip0} and \eqref{Psip0b} depend on $\epsilon$.
\end{remark}

\textbf{Step 7. Estimates on time derivatives.}
In the following, we derive estimates for the time derivatives of $\bm{v}^k$, $\varphi^k$ and $\sigma^k$.
By the Sobolev embedding $H^2(\Omega)\hookrightarrow L^\infty(\Omega)$ and \eqref{energy-a}, \eqref{phikH2-7}, we find
\begin{equation}
\int_0^{T} \|  \varphi^{k}(t)\bm{v}^{k}(t) \|^2\,\mathrm{d}t
\leq \|\bm{v}^{k}\|_{L^{\infty}(0,T;\bm{L}^2(\Omega))}^2 \int_0^{T} \|  \varphi^{k}(t)\|_{H^2(\Omega)}^2 \,\mathrm{d}t  \le C.
\nonumber
\end{equation}
This combined with \eqref{atest.2}, \eqref{energy-a} and \eqref{energy-b} yields
\begin{equation}
\|  \partial_{t}\varphi^{k}\|_{L^{2}(0,T;(H^1(\Omega))')} \le C.
\label{phimt}
\end{equation}

Next, for any test function $\bm{\zeta}\in L^{4}(0,T; \bm{H}^1_{0,\sigma}(\Omega))$, $d=2,3$, we infer from the Gagliardo--Nirenberg inequality  that
\be
\begin{aligned}
&\left|\int_0^{T}\big((\mu^k+\chi
\sigma^k)\nabla\varphi^k, \bm{\zeta}\big) \,\mathrm{d}t\right|
\\
&\quad \le  \int_0^{T} \|\mu^k\|_{L^6(\Omega)}\|\nabla \varphi^k\|_{\bm{L}^\frac32(\Omega)}\|\bm{\zeta}  \|_{\bm{L}^6(\Omega)}\,\mathrm{d}t
+ |\chi|  \int_0^{T} \|\sigma^k\|\|\nabla \varphi^k\|_{\bm{L}^3(\Omega)}\|\bm{\zeta}  \|_{\bm{L}^6(\Omega)}\,\mathrm{d}t
\\
&\quad \le C \|  \mu^{k} \|_{L^2(0,T;H^1(\Omega))}\|\varphi^k\|_{L^{\infty}(0,T;H^1(\Omega))}
 \|\bm{\zeta}  \|_{L^2(0,T;\bm{H}^{1}(\Omega))}
 \\
&\qquad + C \|\sigma^k\|_{L^2(0,T;L^2(\Omega))} \|\varphi^k\|_{L^{\infty}(0,T;H^1(\Omega))}^\frac12
 \|\varphi^k\|_{L^{2}(0,T;H^2(\Omega))}^\frac12
 \|\bm{\zeta}  \|_{L^4(0,T;\bm{H}^{1}(\Omega))}
\\
&\quad \leq C\|\bm{\zeta}  \|_{L^4(0,T;\bm{H}^{1}(\Omega))}.
\end{aligned}
\notag
\ee
Recalling the well-known estimate on the convection term in the Navier--Stokes system
\begin{align}
\| (\bm{v}^{k} \cdot \nabla)\bm{v}^{k}\|_{ L^{\frac{4}{d}}(0,T;(\bm{H}^1_{0,\sigma}(\Omega))')} \le C,
\quad d=2,3,
\notag
\end{align}
then by comparison in \eqref{atest.1}, we obtain
\be
 \| \partial_{t}\bm{v}^{k}\|_{ L^{\frac{4}{3}}(0,T;(\bm{H}^1_{0,\sigma}(\Omega))')}
 \le C.
 \label{vmt3d}
\ee

Since $W^{1,4}(\Omega)\subset L^\infty(\Omega)$ is valid in both two and three dimensions, we deduce from \eqref{atest.4}, H\"{o}lder's inequality and the Sobolev embedding theorem that (cf. \cite[(3.54)]{RSS} without fluid interaction)
\begin{align}
 \|\partial_{t}\sigma^{k}\|_{(W^{1,4}(\Omega))'}
& \leq C \|(\sigma^k)^{\frac12} \|_{L^4(\Omega)} \|(\sigma^k)^{\frac12} \nabla(\ln \sigma^k-\chi\varphi^k)\|
\notag \\
&\quad+C \big(\|\bm{v}^k\|_{\bm{L}^4(\Omega)} \|\sigma^k\|+\|\sigma^k\|+\|\sigma^k\|^2\big)
\notag \\
& \leq C \|\sigma^k \|_{L^1(\Omega)}^\frac18 \|(\sigma^k)^{\frac12} \|_{H^1(\Omega)}^\frac34\|(\sigma^k)^{\frac12} \nabla(\ln \sigma^k-\chi\varphi^k)\|
\notag \\
&\quad+C \big(\|\bm{v}^k\|^\frac14 \|\bm{v}^k\|_{\bm{H}^1(\Omega)}^\frac34 \|\sigma^k\|\big) + C\big(\|\sigma^k\|^2+1\big).
\label{sigmat-0}
\end{align}
Using \eqref{energy-a}, \eqref{energy-b}, \eqref{energy-c}, \eqref{sig-es-aa} and the generalized Young's inequality,
we find that the first two terms on the right-hand side of \eqref{sigmat-0}
are uniformly bounded in $L^{\frac87}(0,T)$, while the last term is uniformly bounded in $L^1(0,T)$. Hence, it holds
\begin{align}
\|\partial_{t}\sigma^{k}\|_{L^1(0,T;(W^{1,4}(\Omega))')}\leq C.
\label{sigmat-a}
\end{align}
We note that in \eqref{phimt}, \eqref{vmt3d} and \eqref{sigmat-a}, the positive constant $C$ depends on $\kappa$, but it is independent on $k$, $\gamma$, $\epsilon$ and $n$.
\medskip

\textbf{Step 8. Improved estimates with more regular initial datum  $\sigma_0$.}
Testing \eqref{atest.4} by $\sigma^k$, using integration by parts, we find
\begin{align}
&\frac12 \frac{\mathrm{d}}{\mathrm{d} t} \|\sigma^k\|^2 + \|\nabla \sigma^k\|^2
+\kappa \int_{\Omega} (\sigma^k)^{3} \mathrm{d}x
  = \chi  \int_{\Omega}   \sigma^k \nabla \varphi^k\cdot\nabla \sigma^k\,  \mathrm{d}x
  + \int_{\Omega} \beta(\varphi^k) (\sigma^k)^2\, \mathrm{d}x.
\label{es-sigmaL2}
\end{align}
It remains to control the first term on the right-hand side.
The argument in \cite{RSS} fails to apply here due to the projection operator $\bm{P}_{Z_k}$. Instead, we take advantage of the $p$-Laplace regularization. Using H\"{o}lder's inequality and the Gagliardo--Nirenberg inequality, we obtain
\begin{align}
\left|\int_{\Omega}  \chi \sigma^k \nabla \varphi^k\cdot\nabla \sigma^k \,\mathrm{d}x\right|
&\leq |\chi|\|\sigma^k\|_{L^4(\Omega)}\|\nabla \varphi^k\|_{\bm{L}^4(\Omega)}\|\nabla \sigma^k\|
\notag \\
&\leq C\|\sigma^k\|^\frac14 \|\nabla \varphi^k\|_{\bm{L}^4(\Omega)}\|\nabla \sigma^k\|^\frac{7}{4}+
C\|\sigma^k\|\|\nabla \varphi^k\|_{\bm{L}^4(\Omega)}
\notag \\
&\leq \frac12 \|\nabla \sigma^k\|^2+ C\big(1+\|\nabla \varphi^k\|_{\bm{L}^4(\Omega)}^8\big) \|\sigma^k\|^2
+ C\|\nabla \varphi^k\|_{\bm{L}^4(\Omega)}^2.
\notag
\end{align}
The above estimate together with \eqref{es-sigmaL2} leads to
\begin{align}
&\frac12 \frac{\mathrm{d}}{\mathrm{d} t} \|\sigma^k\|^2
+ \frac12 \|\nabla \sigma^k\|^2
+ \kappa \|\sigma^k\|_{L^3(\Omega)}^{3}
  \leq  C \big(1+\|\nabla \varphi^k\|_{\bm{L}^4(\Omega)}^8\big) \|\sigma^k\|^2
  + C\|\nabla \varphi^k\|_{\bm{L}^4(\Omega)}^2.
  \notag
\end{align}
Applying Gronwall's lemma and the estimate \eqref{energy-a}, we infer that
\begin{align}
&\|\sigma^k\|_{L^\infty(0,T;L^2(\Omega))\cap L^2(0,T;H^1(\Omega))\cap L^3(0,T;L^3(\Omega))} \leq C.
\label{hi-sig1b}
\end{align}
On the other hand, it follows from \eqref{phikH2-1} and \eqref{phikH2-2} that
\begin{align}
 \|\Delta \varphi^k \|^2
&  \leq \theta_0 \|\nabla \varphi^k\|^2  -\int_\Omega \nabla \mu^k\cdot \nabla \varphi^k\,\mathrm{d}x
+ \chi \int_\Omega   \sigma^k \Delta \varphi^k\,\mathrm{d}x
\notag\\
&  \leq \frac12 \|\Delta \varphi^k \|^2 + \theta_0 \|\nabla \varphi^k\|^2
+ \|\nabla \varphi^k\|\|\nabla \mu^k\| + \frac{\chi^2}{2} \|\sigma^k\|^2.
\label{phikH2-1a}
\end{align}
From \eqref{hi-sig1b}, \eqref{phikH2-1a}, we have
\begin{align}
&\|\varphi^k\|_{L^4(0,T;H^2(\Omega))}\leq C,
\label{hi-vhi1b}
\end{align}
and a comparison in \eqref{atest.4} yields
\begin{align}
\|\partial_t \sigma^k\|_{L^2(0,T;(H^1(\Omega))')}\leq C.
\label{hi-sig2b}
\end{align}
Here, the constant $C>0$ depends on $\kappa$, $\gamma$ and $\|\sigma_{0,n}\|$, but it is independent on $k$ and $\epsilon$.

\begin{remark}
When $d=2$, the $p$-Laplace regularization is again unnecessary (cf. Remark \ref{rem-TM}). Taking $\gamma=0$, we can apply an argument similar to that for \cite[Theorem 2.2]{RSS}. Indeed, the estimate \eqref{phikH2-1a} together with \eqref{energy-a}, \eqref{es-sigmaL2} enables us to conclude (see \cite[Section 4]{RSS} for details)
\begin{align}
&\frac12 \frac{\mathrm{d}}{\mathrm{d} t} \|\sigma^k\|^2 + \frac12 \|\nabla \sigma^k\|^2
+ \kappa \|\sigma^k\|_{L^3(\Omega)}^{3}
+ \frac12 \|\Delta \varphi^k\|^4
  \leq  C\big(1+\|\sigma^k\|^4+\|\nabla \mu^k\|^2\big).
  \notag
\end{align}
Applying Gronwall's lemma, we can recover the estimates \eqref{hi-sig1b}, \eqref{hi-vhi1b}, \eqref{hi-sig2b} by using \eqref{energy-a}, \eqref{energy-b}, \eqref{energy-c}, where the positive constant $C$ depends on $\kappa$ and $\|\sigma_{0,n}\|$, but is independent of $k$, $\epsilon$ (note that $\gamma=0$).
\end{remark}

\subsection{Proof of Theorem \ref{main}}
We first pass to the limit as $k\to +\infty$ and establish the existence of a global weak solution to the regularized problem  $(\bm{S}_{\gamma,\epsilon,n})$.

\bp\label{p2}
Let $d =2,3$, $\gamma\in (0,\frac{1}{2}]$, $\epsilon\in (0,\epsilon_1]$, $n\in \mathbb{Z}^+$ and $T \in (0,+\infty)$. Suppose that the assumptions (H1)--(H6) are satisfied, and the initial data $(\bm{v}_0,\varphi_0,\sigma_0)$ are given as in Theorem \ref{main}. The regularized problem \eqref{reg.sys}--\eqref{reg.ini} admits a global weak solution $(\bm{v}^{\gamma,\epsilon,n},\varphi^{\gamma,\epsilon,n}, \mu^{\gamma,\epsilon,n}, \sigma^{\gamma,\epsilon,n})$ on $[0,T]$ such that
\begin{align}
&\bm{v}^{\gamma,\epsilon,n} \in L^{\infty}(0,T;\bm{L}^2_{0,\sigma}(\Omega)) \cap L^{2}(0,T;\bm{H}^1_{0,\sigma}(\Omega))\cap  W^{1,\frac{4}{3}}(0,T;(\bm{H}^1_{0,\sigma}(\Omega))'),
\notag \\
&\varphi^{\gamma,\epsilon,n} \in L^{\infty}(0,T;W^{1,4}(\Omega))\cap L^{4}(0,T;H^2_{N}(\Omega)) \cap H^{1}(0,T;(H^1(\Omega))'),
\notag \\
&\mu^{\gamma,\epsilon,n} \in   L^{2}(0,T;H^1(\Omega)),\quad |\nabla \varphi^{\gamma,\epsilon,n}|^2\nabla \varphi^{\gamma,\epsilon,n}\in L^2(0,T;\bm{H}^1(\Omega)),
\notag\\
& \Psi_{\epsilon}(\varphi^{\gamma,\epsilon,n})\in L^\infty(0,T;L^1(\Omega)),\quad \Psi_{\epsilon}'(\varphi^{\gamma,\epsilon,n})\in L^2(0,T;L^2(\Omega)),
\notag \\
&\sigma^{\gamma,\epsilon,n} \in L^{\infty}(0, T ; L^{2}(\Omega))\cap L^2(0,T;H^1(\Omega))\cap H^1(0,T;(H^1(\Omega))'),
\notag \\
&\sigma^{\gamma,\epsilon,n}(x, t) \geq 0 \quad \text {a.e. in}\ \Omega\times(0,T),
\notag
\notag
\end{align}
and the following identities hold
\begin{subequations}
\begin{alignat}{3}
& \langle\partial_t  \bm{ v}^{\gamma,\epsilon,n} ,\bm{\zeta} \rangle_{(\bm{H}^1_{0,\sigma})',\bm{H}^1_{0,\sigma}}
{\color{black}{- (\bm{v}^{\gamma,\epsilon,n}  \otimes\bm{ v}^{\gamma,\epsilon,n} ,D\bm{ \zeta})}}
+\big(2\eta(\varphi^{\gamma,\epsilon,n} ) D\bm{v}^{\gamma,\epsilon,n} ,D\bm{ \zeta}\big)
&& \notag \\
& \quad =\big((\mu^{\gamma,\epsilon,n} +\chi\sigma^{\gamma,\epsilon,n} ) \nabla \varphi^{\gamma,\epsilon,n} , \bm {\zeta}\big), \quad&& \textrm{a.e. in }(0,T), &\label{test3.c-reg}
\\
& \langle \partial_t \varphi^{\gamma,\epsilon,n} ,\xi \rangle_{(H^1(\Omega))',H^1(\Omega)}
+ ({\bm{v}^{\gamma,\epsilon,n}  \cdot \nabla \varphi^{\gamma,\epsilon,n} },\xi)+ \big(m(\varphi^{\gamma,\epsilon,n} )\nabla \mu^{\gamma,\epsilon,n} ,\nabla \xi\big)
&& \notag \\
& \quad =
\big(-\alpha \varphi^{\gamma,\epsilon,n} + h(\varphi^{\gamma,\epsilon,n} , \sigma^{\gamma,\epsilon,n} ),\xi\big),\quad && \textrm{a.e. in }(0,T),&\label{test1.a-reg}
\\
& \mu^{\gamma,\epsilon,n} = -\gamma^{8}\mathrm{div}(|\nabla \varphi^{\gamma,\epsilon,n} |^2 \nabla \varphi^{\gamma,\epsilon,n} )- \Delta \varphi^{\gamma,\epsilon,n}  +\Psi_{\epsilon}'(\varphi^{\gamma,\epsilon,n} ) -\chi \sigma^{\gamma,\epsilon,n},\quad && \textrm{a.e. in }\Omega\times(0,T),&\label{test4.d-reg}
\\
& \langle \partial_t\sigma^{\gamma,\epsilon,n},\xi \rangle_{(H^1(\Omega))',H^1(\Omega)}
- (\bm{v}^{\gamma,\epsilon,n}  \sigma^{\gamma,\epsilon,n} ,\nabla \xi) + (\nabla\sigma^{\gamma,\epsilon,n} , \nabla \xi)
&& \notag\\
&\quad = \chi (\sigma^{\gamma,\epsilon,n} \nabla \varphi^{\gamma,\epsilon,n},\nabla \xi)
+ \big(\beta(\varphi^{\gamma,\epsilon,n} )\sigma^{\gamma,\epsilon,n}  - \kappa\big(\sigma^{\gamma,\epsilon,n} \big)^2, \xi\big),
&& \textrm{a.e. in }(0,T)& \label{test2.b-reg}
\end{alignat}
\end{subequations}
 for all test functions $\bm {\zeta} \in \bm{H}^{1}_{0,\sigma}(\Omega)$, $\xi\in H^1(\Omega)$. Moreover, the initial conditions are fulfilled
\begin{align}
&\bm{v}^{\gamma,\epsilon,n} |_{t=0}=\bm{v}_{0},\quad \varphi^{\gamma,\epsilon,n} |_{t=0}=\varphi_{0,\gamma},\quad
\sigma^{\gamma,\epsilon,n} |_{t=0}=\sigma_{0,n}, \quad\textrm{a.e. in } \Omega.	
\notag
\end{align}
\ep
\begin{proof}
We have already shown that for every $k\geq \widehat{k}$, the semi-Galerkin scheme \eqref{atest.1}--\eqref{atest.ini0} admits a unique global solution $(\bm{v}^{k},\varphi^{k},\mu^{k},\sigma^{k})$ on $[0,T]$ with suitable estimates that are independent of $k$. The uniform estimates \eqref{energy-a}, \eqref{energy-b}, \eqref{mu},
 \eqref{Psip0}, \eqref{sigma}, \eqref{phimt}, \eqref{vmt3d},
\eqref{hi-sig1b}, \eqref{hi-vhi1b} and \eqref{hi-sig2b} are sufficient for us to apply theorems of weak compactness and the Aubin--Lions--Simon lemma \cite{si87} to extract a suitable subsequence that approach a limit in corresponding topologies as $k\to +\infty$. Here, we just mention that in order to show \eqref{test4.d-reg}, we have applied \cite[Theorem 2.1]{CM19} on the optimal second-order regularity for the $p$-Laplace problem to obtain the regularity property $|\nabla \varphi^{\gamma,\epsilon,n}|^2\nabla \varphi^{\gamma,\epsilon,n}\in L^2(0,T;\bm{H}^1(\Omega))$.
Since the remaining procedure is standard, we omit the details.
\end{proof}

Now we are in a position to prove our main result.
\medskip
\\
\noindent
\textbf{Proof of Theorem \ref{main}.} For simplicity, we shall take the limit $\gamma, \epsilon\to 0$ and $n\to +\infty$ simultaneously. Hence, we can take $\gamma=\epsilon=\frac{1}{n}$ (for sufficiently large $n$) and simply denote the approximate solutions $(\bm{v}^{\gamma,\epsilon,n},\varphi^{\gamma,\epsilon,n},\mu^{\gamma,\epsilon,n}, \sigma^{\gamma,\epsilon,n})$ obtained in Proposition \ref{p2} by $(\bm{v}^{n},\varphi^{n},\mu^{n}, \sigma^{n})$.
In view of the uniform estimates obtained in Section \ref{unie}, Lebesgue's dominated convergence theorem, Fatou's lemma, and the weak lower semicontinuity of norms,  we have the following estimates that are independent of $n$ (and thus for $\gamma=\epsilon=\frac{1}{n}$):
\begin{align}
&\|\bm{v}^n\|_{L^{\infty}(0, T; \bm{L}^2(\Omega))}
 +\gamma^2\|\varphi^n\|_{L^{\infty}(0, T; W^{1,4}(\Omega))}
+\|\varphi^n\|_{L^{\infty}(0, T; H^1(\Omega))}
\notag \\
&\quad +\|\sigma^n \ln \sigma^n\|_{L^{\infty}(0, T; L^{1}(\Omega))} +\|\Psi_{0,\epsilon}(\varphi^n)\|_{L^{\infty}(0, T; L^{1}(\Omega))}
\leq C,
\label{energy-an}
\end{align}
\begin{align}
&\|\bm{v}^n\|_{L^{2}(0, T; \bm{H}^1(\Omega))}
+\|\mu^n\|_{L^{2}(0, T ; H^{1}(\Omega))}
+ \|\varphi^n\|_{L^{2}(0, T; H^2(\Omega))}
\leq C,
\label{energy-bn}
\end{align}
\begin{align}
& \gamma^4\big\||\nabla \varphi^n|^2\big\|_{L^2(0,T;H^1(\Omega))}\le C,
\label{pL-L2H1n}
\end{align}
\begin{align}
\| \sigma^n\|_{L^{2}(0, T; L^{2}(\Omega))}
+ \|\sigma^n\|_{L^{\frac{4}{3}}(0, T ; W^{1,\frac{4}{3}}(\Omega))}
+ \|(\sigma^n)^2 \ln \sigma^n\|_{L^{1}(0, T; L^{1}(\Omega))}
\leq C,
\label{energy-cn}
\end{align}
\begin{equation}
\| \partial_{t}\bm{v}^{k}\|_{ L^{\frac{4}{3}}(0,T;(\bm{H}^1_{0,\sigma})')}
+ \|  \partial_{t}\varphi^{k}\|_{L^{2}(0,T;(H^1(\Omega))')} \le C.
\label{phimtn}
\end{equation}
Using the fact $H^2(\Omega)\hookrightarrow W^{1,6}(\Omega)\hookrightarrow L^\infty(\Omega)$, we find
\begin{align}
 \|\partial_{t}\sigma^{n}\|_{(H^2(\Omega))'}
& \leq C \|\bm{v}^n\|_{\bm{L}^3(\Omega)} \|\sigma^n\|
+ C \| \sigma^n\| \| \nabla \varphi^n\|_{\bm{L}^3(\Omega)}
+ C(\|\sigma^n\|+\|\sigma^n\|^2\big)
\notag \\
& \leq C \big(\|\bm{v}^n\|^\frac12 \|\bm{v}^n\|_{\bm{H}^1(\Omega)}^\frac12 \|\sigma^n\|\big)
+ C \| \sigma^n\| \| \nabla \varphi^n\|_{\bm{L}^3(\Omega)}
+ C\big(\|\sigma^n\|^2+1\big),
\notag
\end{align}
which combined with \eqref{energy-bn}, \eqref{energy-cn} yields
\begin{align}
 \|\partial_{t}\sigma^{n}\|_{L^1(0,T;(H^2(\Omega))')}\leq C.
 \label{sigmatn}
\end{align}
In order to pass to the limit, we also need to handle the nonlinear term $\Psi'_{0,\epsilon}(\varphi^n)$. Since for every $n$, $\Psi_{0,\epsilon}(\varphi^n)\in L^2(0,T;H^1(\Omega))$ (though the corresponding estimate depends on $n$), then we are allowed to test \eqref{test4.d-reg} by $\Psi'_{0,\epsilon}(\varphi^n)$. After integration by parts, we find
\begin{align}
 &\underbrace{\int_\Omega \big(\gamma^{8}|\nabla \varphi^{n} |^4 + |\nabla \varphi^{n}|^2 \big) \Psi''_{0,\epsilon}(\varphi^n) \,\mathrm{d}x }_{\geq 0}
 + \|\Psi_{0,\epsilon}'(\varphi^{n} ) \|^2\notag \\
&\quad = \int_\Omega\big(\mu^{n}+\chi \sigma^{n}+\theta_0\varphi^n\big)\Psi_{0,\epsilon}'(\varphi^{n} ) \,\mathrm{d}x
\notag\\
&\quad \leq \frac12 \|\Psi_{0,\epsilon}'(\varphi^{n} ) \|^2+ C \big(\|\mu^{n}\|^2+\|\sigma^{n}\|^2+\|\varphi^n\|^2\big).
\notag
\end{align}
Integrating in time over $[0,T]$, we conclude from \eqref{energy-bn} and \eqref{energy-cn} that
\begin{align}
	\|\Psi'_{\epsilon}(\varphi^n)\|_{L^2(0,T;L^2(\Omega))}\leq C.
	\label{Psip0n}
\end{align}

Collecting the above uniform estimates that are independent of $n$, by classical compactness results, we are able to pass to the limit as $n\to +\infty$ (in the sense of a suitable non-relabelled subsequence). We denote the limit by $(\bm{v},\varphi,\mu,\sigma)$, which is a global weak solution to the original problem \eqref{sysNew}--\eqref{ini0new} on $[0,T]$. Comparing with \cite{RSS}, no truncation was used in our approximating scheme for the chemical concentration $\sigma$, thus the procedure to pass to the limit turns out to be simpler. The limiting procedure for the subsystem of $(\bm{v}^n,\varphi^n)$ is similar to that for the Navier--Stokes--Cahn--Hilliard system, see, e.g., \cite{B, MT}, while the limiting procedure for the equation of $\sigma^n$ can be done as in \cite{Lan2016}. In particular,   thanks to
\eqref{energy-bn} and \eqref{phimtn}, the Aubin--Lions--Simon theorem  entails that (up to a subsequence) $\varphi^n \to \varphi$ strongly in $L^2(0,T; W^{1,q}(\Omega))$, for any $q\in [2, 6)$. Owing to the weak convergence $\sigma^n  \rightharpoonup \sigma \in L^2(0,T; L^2(\Omega))$, the latter allows the  passage to the limit in the nonlinear terms involving the product $\sigma^n \nabla \varphi^n$.
Besides, due to the singularity of $\Psi$, we can apply the classical method for the Cahn--Hilliard equation to conclude that (see, e.g., \cite{GGW,Mi19})
$$
\varphi\in L^{\infty}(\Omega\times (0,T))\ \textrm{with}\ \ |\varphi(x,t)|<1\ \ \textrm{a.e.\ in}\ \Omega\times(0,T).
$$
If in addition, $\sigma_0\in L^2(\Omega)$, improved regularity of the weak solution $(\varphi, \sigma)$ can be derived from the uniform estimates obtained in Step 8 in Section \ref{unie} in two dimensions. In particular, we can apply the classical result for the Cahn--Hilliard equation to obtain (see, e.g., \cite{CG,GG,GGW})
\begin{align}
\|\varphi\|_{L^2(0,T;W^{2,q}(\Omega))} \textcolor{black}{+\|\Psi'(\varphi)\|_{L^2(0,T;L^q(\Omega))}} \leq C,\quad \forall\, q\in [2,+\infty).
\notag
\end{align}
This concludes the proof of Theorem \ref{main} and we omit the details.
\hfill $\square$

\section{Uniqueness of Weak Solution in Two Dimensions}
\setcounter{equation}{0}
\label{sec:uni}

In this section, we prove Theorem \ref{main-2} on the uniqueness of global weak solutions in two dimensions. The proof mainly follows the idea in \cite{GMT} for the two dimensional Navier--Stokes--Cahn--Hilliard system with non-constant viscosity. Additional efforts are required to handle the $\sigma$-equation, cf. \cite{H}.
\medskip 

\noindent \textbf{Proof of Theorem \ref{main-2}.} 
Let $(\bm{v}_{1},\varphi_{1},\mu_1,\sigma_{1})$ and $(\bm{v}_{2},\varphi_{2},\mu_2,\sigma_{2})$ be two weak solutions to problem \eqref{maineq}--\eqref{ini0} given by Theorem \ref{main}-(2) subject to the initial data $(\bm{v}_{01},\sigma_{01},\varphi_{01})$ and $(\bm{v}_{02},\sigma_{02},\varphi_{02})$, respectively. 
Denote the differences by 
$$
(\bm{v},\varphi,\mu,\sigma)=(\bm{v}_{1}-\bm{v}_{2}, \varphi_{1}-\varphi_{2},\mu_1-\mu_2,\sigma_{1}-\sigma_{2}).
$$
Following the detailed calculations performed in \cite[Section 4]{H}, we deduce that
\begin{align}
\frac{\mathrm{d}}{\mathrm{d}t}W(t)+ \frac{\eta_*}{2} \|\bm{v}\|^2 + \frac{\widehat{m}}{2}\|\nabla\varphi\|^2 + \frac12 \|\sigma\|^2
\leq \widehat{C}Z(t)W(t)\ln\left(\frac{\widetilde{C}}{W(t)}\right),
\label{uniA}
\end{align}
where  
\begin{align}
 W(t):=\|\nabla\bm{S}^{-1}\bm{v}(t)\|^2+\|\varphi(t)\|_{(H^1(\Omega))'}^2 +\|\sigma(t)\|_{(H^1(\Omega))'}^2+|\overline{\varphi}(t)|,
 \notag
 \end{align}
 and 
\begin{align}
Z(t)&=\|\nabla \bm{v}_{1}(t)\|^2+\|\nabla \bm{v}_{2}(t)\|^2+\|\varphi_{1}(t)\|_{W^{2,3}(\Omega)}^2 +\|\varphi_{2}(t)\|_{W^{2,3}(\Omega)}^2
+\|\varphi_{1}(t)\|_{H^{2}(\Omega)}^4
\notag\\
&\quad +\|\Psi'(\varphi_{1}(t))\|_{L^1(\Omega)} +\|\Psi'(\varphi_{2}(t))\|_{L^1(\Omega)}
+\|\sigma_2(t)\|_{H^1(\Omega)}^2 
\notag\\
&\quad 
 + \|\sigma_2(t)\|^4
+ \|\sigma_1(t)\|^2\|\sigma_1(t)\|_{H^1(\Omega)}^2 +1. 
\notag 
\end{align}
Here, $\widehat{C}$, $\widetilde{C}$ are positive constants depending on norms of the initial data, $\Omega$, and the coefficients of the system, $\widetilde{C}$ may also depend on $T$. In the derivation of \eqref{uniA}, we only need modifications due to the new cross diffusion term as well as the logistic term, that is,
\begin{align}
&|\chi(\sigma_1 \nabla \varphi_1,\nabla \mathcal{N}_1\sigma)-\chi(\sigma_2 \nabla \varphi_2,\nabla \mathcal{N}_1\sigma)|
\notag \\
&\quad \leq |\chi(\sigma_1 \nabla \varphi,\nabla \mathcal{N}_1\sigma)|+|\chi(\sigma \nabla \varphi_2,\nabla \mathcal{N}_1\sigma)|
\notag \\
&\quad \leq C\|\sigma_1\|_{L^4(\Omega)}\|\nabla \varphi\|\|\nabla \mathcal{N}_1\sigma\|_{\bm{L}^4(\Omega)}
+C\|\sigma\|\| \nabla \varphi_2\|_{\bm{L}^\infty(\Omega)}\|\nabla \mathcal{N}_1\sigma\|
\notag \\
&\quad \leq C\|\sigma_1\|^\frac12\|\sigma_1\|_{H^1(\Omega)}^\frac12 \|\nabla \varphi\| \|\nabla \mathcal{N}_1\sigma\|^\frac12\|\sigma\|^\frac12
+C \|\sigma\|\| \varphi_2\|_{W^{2,3}(\Omega)}\|\nabla \mathcal{N}_1\sigma\|
\notag\\
&\quad \leq \frac{\widehat{m}}{4}\|\nabla \varphi\|^2 + \frac{1}{8}\|\sigma\|^2 +C \|\sigma_1\|^2\|\sigma_1\|_{H^1(\Omega)}^2  \|\sigma\|_{(H^1(\Omega))'}^2
+ C\| \varphi_2\|_{W^{2,3}(\Omega)}^2\|\sigma\|_{(H^1(\Omega))'}^2,
\notag 
\end{align}
and by H\"{o}lder's inequality together with Agmon's inequality, 
\begin{align*}
 \left|\int_\Omega (\kappa \sigma_1^2-\kappa\sigma_2^2)\mathcal{N}_1\sigma\,\mathrm{d}x\right|
& \leq C\big(\|\sigma_1\|+\|\sigma_2\|)\|\sigma\| \|\mathcal{N}_1\sigma\|_{L^\infty(\Omega)}
\\
& \leq C\big(\|\sigma_1\|+\|\sigma_2\|)\|\sigma\| \|\mathcal{N}_1\sigma\|^\frac12\|\sigma\|^\frac12
\\
& \leq \frac{1}{8}\|\sigma\|^2+ C \big(\|\sigma_1\|^4+ \|\sigma_2\|^4)  \|\sigma\|_{(H^1(\Omega))'}^2.
\end{align*}
Integrating \eqref{uniA} on $[0,t]\subset [0,T]$, we get
\begin{align}
W(t)\leq W(0)+ C_1\int_0^t Z(s)W(s)\ln \left(\frac{\widetilde{C}}{W(s)}\right)ds,\quad \text{for a.e.}\, t\in [0,T].
\label{uniA1}
\end{align}
From the regularity of weak solutions, it holds
$$
Z(t)\in L^1(0,T).
$$
Then we are in a position to apply the Osgood lemma (see \cite[Lemma 3.4]{BCD}): if $W(0)=0$, we find $W(t) = 0 $ for all $t \in [0, T]$, which yields the uniqueness of global weak solutions to problem \eqref{maineq}--\eqref{ini0};
if $W(0)>0$, we can derive a continuous dependence estimate like in \cite{GMT,H}, that is,
\be
W(t)\le \widetilde{C}\left(\frac{W(0)}{\widetilde{C}}\right)^{\exp\left(-\widehat{C}\int_{0}^{t} Z(s)\, \mathrm{d}s\right)},\quad \ \forall\, t\in[0,T].
\label{contiuniq}
\ee
The proof of Theorem \ref{main-2} is complete. 
\hfill $\square$.

\appendix
\section{Appendix}
\setcounter{equation}{0}

\subsection{Proof of Lemma \ref{fp}}
\label{app-1}
The existence of a classical solution $\sigma^k$ to problem \eqref{1atest.4}--\eqref{boundary2} can be proved by using the semigroup approach (see e.g., \cite{W2010,W2012}).
Assume $\|\sigma_{0,n}\|_{L^\infty(\Omega)}\leq K$ for some $K>0$. Consider $T_K\in (0,1)$ to be specified below and the closed set
$$
\Sigma_K= \big\{\sigma^k\in C([0,T_K];C(\overline{\Omega}))\,|\, \|\sigma^k\|_{L^\infty(0,T_K;L^\infty(\Omega))}\leq K+1\big\}.
$$
For $\sigma^k\in \Sigma_K$ and $t\in [0,T_K]$, we define the nonlinear mapping $\mathcal{G}$ as follows
\begin{align}
\mathcal{G}(\sigma^k)(t)
& =e^{t\Delta}\sigma_{0,n}-\int_0^t e^{(t-s)\Delta}\big[\bm{u}^k(s)\cdot\nabla\sigma^k(s)+ \chi \mathrm{div}(\sigma^k(s)\nabla \psi^k(s))\big]\,\mathrm{d}s
\notag \\
&\quad + \int_0^t e^{(t-s)\Delta}\big[\beta(\psi^k(s))  \sigma^k(s)-\kappa(\sigma^k(s))^2\big]\,\mathrm{d}s.
\label{G}
\end{align}
Here, $(e^{t\Delta})_{t\geq 0}$ denotes the Neumann heat semigroup in $\Omega$. By the maximum principle, we find
$$
\|e^{t\Delta}\sigma_{0,n}\|_{L^\infty(\Omega)}\leq \|\sigma_{0,n}\|_{L^\infty(\Omega)}\leq K,\quad \forall\,t\in [0,T_K],
$$
and
\begin{align*}
& \int_0^t \big\|e^{(t-s)\Delta}\big[\beta(\psi^k(s)) \sigma^k(s)-\kappa(\psi^k(s)) (\sigma^k(s))^2 \big]\big\|_{L^\infty(\Omega)}\,\mathrm{d}s
\\
&\quad \leq ( b^*+\kappa) \int_0^t \big(\|\sigma^k(s)\|_{L^\infty(\Omega)}+
\|\sigma^k(s)\|_{L^\infty(\Omega)}^2\big)\,\mathrm{d}s\\
&\quad \leq ( b^*+\kappa)(K+2)^2T_K,
\quad
\forall\, t\in [0,T_K].
\end{align*}
We recall that $\mathcal{A}_1=-\Delta +I$ subject to the homogeneous Neumann boundary condition is a sectorial operator in $L^4(\Omega)$. In two and three dimensions, it holds $D(\mathcal{A}_1^{\frac{7}{16}})\hookrightarrow C(\overline{\Omega})$. Applying \cite[Lemma 1.3]{W2010jde}, the incompressibility condition and the fact $T_K\in (0,1)$, we can conclude that
\begin{align*}
&\int_0^t \|e^{(t-s)\Delta}\big[\bm{u}^k(s)\cdot\nabla\sigma^k(s)+ \chi \mathrm{div}(\sigma^k(s)\nabla \psi^k(s))\big]\|_{L^\infty(\Omega)}\,\mathrm{d}s
\\
&\quad \leq C \int_0^t \|\mathcal{A}_1^{\frac{7}{16}}e^{(t-s)\Delta} \mathrm{div}\big(\bm{u}^k(s) \sigma^k(s)+\chi\sigma^k(s)\nabla \psi^k(s)\big) \|_{L^4(\Omega)}\,\mathrm{d}s\\
&\quad \leq C \int_0^t (t-s)^{-\frac{15}{16}}\|\bm{u}^k(s) \sigma^k(s)+\chi\sigma^k(s)\nabla \psi^k(s)\|_{L^4(\Omega)}\,\mathrm{d}s\\
&\quad \leq C\max\{1,\,|\chi|\}T_K^{\frac{1}{16}} \sup_{s\in [0,T_K]}
\big(\|\bm{u}^k(s)\|_{\bm{L}^4(\Omega)} +\|\nabla \psi^k(s)\|_{\bm{L}^4(\Omega)}\big)\|\sigma^k(s)\|_{L^\infty(\Omega)}\\
&\quad \leq C_k \widetilde{M}^\frac12 (K+1) T_K^{\frac{1}{16}}.
\end{align*}
Combining the above estimates, we infer from \eqref{G} that
$$
 \|\mathcal{G}(\sigma^k)(t)\|_{L^\infty(\Omega)}
 \leq K + (b^* +\kappa)(K+2)^2T_K + C_k \widetilde{M}^\frac12 (K+1) T_K^{\frac{1}{16}},\quad \forall\, t\in [0,T_K].
$$
Using the same idea, given $\sigma^k_1,\,\sigma^k_2\in \Sigma_K$, we can estimate the difference
\begin{align*}
 & \|\mathcal{G}(\sigma^k_1)(t)-\mathcal{G}(\sigma^k_2)(t)\|_{L^\infty(\Omega)}\\
 &\quad \leq  \int_0^t \big\|e^{(t-s)\Delta}\big[\beta(\psi^k(s)) (\sigma_1^k(s)-\sigma_2^k(s)) \big]\big\|_{L^\infty(\Omega)}\,\mathrm{d}s
 \\
 &\qquad + \int_0^t \big\|e^{(t-s)\Delta}\big[ \kappa (\sigma_1^k(s)-\sigma_2^k(s))(\sigma_1^k(s)+\sigma_2^k(s))\big] \big\|_{L^\infty(\Omega)}\,\mathrm{d}s\\
 &\qquad +\int_0^t \big\|e^{(t-s)\Delta} \mathrm{div}[\bm{u}^k(s)(\sigma_1^k(s)-\sigma_2^k(s))] \big\|_{L^\infty(\Omega)}\,\mathrm{d}s
 \\
 &\qquad +\int_0^t \big\|e^{(t-s)\Delta} \chi \mathrm{div}\big[(\sigma_1^k(s)-\sigma_2^k(s))\nabla \psi^k(s)\big]\big\|_{L^\infty(\Omega)}\,\mathrm{d}s
 \\
 &\quad \leq (b^* +\kappa^*) (2K+3) T_K\sup_{s\in [0,T_K]}\|\sigma_1^k(s)-\sigma_2^k(s)\|_{L^\infty(\Omega)}
 \\
 &\qquad
  + C_k \widetilde{M}^\frac12 T_K^{\frac{1}{16}}\sup_{s\in [0,T_K]}\|\sigma_1^k(s)-\sigma_2^k(s)\|_{L^\infty(\Omega)},\quad \forall\, t\in [0,T_K].
\end{align*}
Hence, taking $T_K\in (0,1)$ sufficiently small, we find that $\mathcal{G}$ maps $\Sigma_K$ into itself and is indeed a contraction on $\Sigma_K$. Applying the contraction mapping principle, we obtain a unique $\sigma^k\in \Sigma_K$ such that
$$\mathcal{G}(\sigma^k)=\sigma^k.$$

Since $\bm{u}^k$, $\psi^k$ are sufficiently regular, according to standard bootstrap
arguments involving the regularity theories for second order parabolic equations, we can show that $\sigma^k$ is a classical solution to problem \eqref{1atest.4}--\eqref{boundary2} on $[0,T_K]$. A standard
extension argument yields that there exists a maximal existence time $T_{\max}\in (0,T]$ such that problem \eqref{1atest.4}--\eqref{boundary2} admits a classical solution $\sigma^k$ on $[0,T_{\max})$, moreover, if $T_{\max}<T$, then $\|\sigma^k(t)\|_{L^\infty(\Omega)}\to +\infty$ as $t\nearrow T_{\max}$.

Since $\sigma_{0,n}\geq 0$, the weak maximum principle for parabolic equations (or the Stampacchia truncation argument) easily yields that $\sigma^k(x,t)\geq 0$ in $\overline{\Omega}\times [0,T_{\max})$. Then, the strong maximum principle implies that $\sigma^k(x,t)> 0$ for $\overline{\Omega}\times (0,T_{\max})$. If this is not the case, we would have $\sigma^k\equiv 0$ on $\overline{\Omega}\times [0,t_0]$ for some $t_0\in (0,T_{\max})$, in particular, $\sigma^k(0)=\sigma_{0,n}\equiv 0$, which leads to a contradiction to the assumption that $\sigma_{0,n} \not\equiv 0$.

Next, we show that $T_{\max}=T$.
Testing \eqref{1atest.4} with $\sigma^{k}$ and using integration by parts, we find
\be
\begin{aligned}
	 &\frac{1}{2} \frac{\mathrm{d}}{\mathrm{d} t}\|\sigma^{k}\|^{2}+\|\nabla \sigma^{k}\|^{2}= \chi \int_{\Omega} \sigma^k\nabla \psi^k \cdot \nabla \sigma^k\,\mathrm{d}x
+\int_{\Omega}\big[\beta(\psi^k) (\sigma^k)^2-\kappa(\sigma^k)^3 \big]\,\mathrm{d}x.
\label{1tot1}
\end{aligned}
\ee
The first term on the right-hand side can be estimated by using H\"{o}lder's and Young's inequalities,
\begin{align}
 \chi \int_{\Omega} \sigma^k \nabla \psi^k \cdot \nabla \sigma^k \,\mathrm{d}x
& \le |\chi| \|\sigma^k\| \|\nabla \psi^k\|_{\bm{L}^{\infty}(\Omega)} \|\nabla \sigma^k\|
\notag \\
& \le  C_k \|\psi^k\|^2\|\sigma^k\|^2 + \frac12\|\nabla \sigma^k\|^2.
\notag
\end{align}
From \eqref{1tot1} and nonnegativity of $\sigma^k$, we find
\begin{align}
  \frac{\mathrm{d}}{\mathrm{d} t}\|\sigma^{k}\|^{2} + \|\nabla \sigma^{k}\|^{2}
  \leq C_k (\|\psi^k\|^2 +1)\|\sigma^k\|^2.
\label{1tot3}
\end{align}
It follows from Gronwall's lemma that
\begin{align}
\|\sigma^{k}(t)\|^{2}
 &\leq \|\sigma_{0,n}\|^2 \, \mathrm{e}^{Ct(M_1+1)},
 \quad \forall\, t\in [0,T_{\max}),
\label{3eeenerg-a}
\end{align}
where
$$
M_1=\sup_{t\in[0,T]}\|\psi^k(t)\|^2.
$$
Here, the positive constant $C$ depends on $k$, but it is independent of $M_1$. Next, testing \eqref{1atest.4} by $-\Delta\sigma^{k}$ yields
\begin{align}
\frac{1}{2} \frac{\mathrm{d}}{\mathrm{d} t}\|\nabla\sigma^{k}\|^{2}+\|\Delta\sigma^{k}\|^{2}
&=  \int_{\Omega} \bm{u}^k\cdot\nabla \sigma^k \Delta \sigma^k\,\mathrm{d}x
+\chi \int_{\Omega} \nabla\sigma^k\cdot\nabla \psi^k \Delta \sigma^k\,\mathrm{d}x
\notag \\
&\quad+\chi \int_{\Omega} \sigma^k\Delta \psi^k \Delta \sigma^k\,\mathrm{d}x
-\int_{\Omega}\big[\beta(\psi^k) \sigma^k -\kappa (\sigma^k)^2\big]\Delta\sigma^k\,\mathrm{d}x
\notag \\
&=: \sum_{j=1}^4 K_j.
\label{2tot1}
\end{align}
Applying H\"{o}lder's, Young's inequalities, we have
\begin{equation}
\begin{aligned}
K_1& \le \|\bm{u}^k\|_{\bm{L}^{\infty}(\Omega)} \|\nabla \sigma^k\| \|\Delta \sigma^k\|
\le  C_k \|\bm{u}^k\|^2\|\nabla\sigma^k\|^2 + \frac16\|\Delta \sigma^k\|^2,
\notag
\\
K_2& \le |\chi| \|\nabla\sigma^k\|  \|\nabla \psi^k\|_{\bm{L}^{\infty}(\Omega)} \|\Delta \sigma^k\|
\le  C_k \|\psi^k\|^2\|\nabla\sigma^k\|^2 + \frac16\|\Delta \sigma^k\|^2,
\notag
\\
K_3
& \le |\chi| \|\sigma^k\|_{L^6(\Omega)} \|\Delta \psi^k\|_{L^{3}(\Omega)} \|\Delta \sigma^k\|
\le  C_k \|\psi^k\|^2\|\sigma^k\|_{H^1(\Omega)}^2 + \frac16 \|\Delta \sigma^k\|^2,
\notag
\end{aligned}
\end{equation}
and
\begin{equation}
\begin{aligned}
K_4
& =\int_{\Omega}(\beta(\psi^k) -2\kappa  \sigma^k)|\nabla\sigma^k|^2\,\mathrm{d}x
+\int_{\Omega} \beta'(\psi^k) \sigma^k  \nabla\psi^k\cdot\nabla\sigma^k\,\mathrm{d}x
\\
&\le  b^*\|\nabla\sigma^k\|^2
+C \|\nabla\psi^k\|_{\bm{L}^{\infty}(\Omega)}\|\sigma^k\|\|\nabla\sigma^k\|
\\
&\le C_k( \|\psi^k\|^2+1)\|\sigma^k\|_{H^1(\Omega)}^2.
\notag
\end{aligned}
\end{equation}
Combining the above estimates, we infer from \eqref{1tot3} and \eqref{2tot1} that
\begin{equation}
\begin{aligned}
&\frac{\mathrm{d}}{\mathrm{d} t}\|\sigma^{k}\|^{2}_{H^1(\Omega)}+\| \sigma^{k}\|^{2}_{H^2(\Omega)}
\leq C_k \big( \|\bm{u}^k\|^2+\|\psi^k\|^2+1 \big)\|\sigma^k\|^2_{H^1(\Omega)}.
\label{2tot5}
\end{aligned}
\end{equation}
Then it follows from Gronwall's lemma and \eqref{2tot5} that
\begin{align}
\|\sigma^{k}(t)\|^{2}_{H^1(\Omega)}
& \leq \|\sigma_{0,n}\|_{H^1(\Omega)}^2 \mathrm{e}^{C t(M_2+1)},
\quad \forall\, t\in [0,T_{\max}),
\label{4eeenerg-a}
\end{align}
where
\begin{align*}
& M_2 = \sup_{t\in[0,T]}\big(\|\bm{u}^k(t)\|^2+  \|\psi^k(t)\|^2\big).
\end{align*}

Let us proceed with the estimate the $L^\infty$-norm of $\sigma^k$.
The Sobolev embedding theorem in three dimensions together with \eqref{4eeenerg-a} implies that $\|\sigma^{k}(t)\|_{L^6(\Omega)}\leq C$ for all $t\in [0,T_{\max})$.  For any $q\geq 6$, testing \eqref{1atest.4} with $(\sigma^k)^{q-1}$, integrating over $\Omega$, using the incompressibility condition and integration by parts, we find
\begin{align}
	&\frac{1}{q} \frac{\mathrm{d}}{\mathrm{d} t}\int_\Omega (\sigma^k)^q\,\mathrm{d}x
+\frac{4(q-1)}{q^2}\int_\Omega |\nabla(\sigma^k)^\frac{q}{2}|^2\,\mathrm{d}x
+  \int_\Omega \kappa(\sigma^k)^{q+1}\,\mathrm{d}x
\notag \\
	&\quad=  \frac{2\chi(q-1)}{q}\int_\Omega (\sigma^k)^\frac{q}{2}\nabla(\sigma^k)^\frac{q}{2}\cdot \nabla \psi^k\,\mathrm{d}x
+ \int_\Omega  \beta(\psi^k)(\sigma^k)^q\,\mathrm{d}x
\notag \\
&\quad \leq \frac{2(q-1)}{q^2}\int_\Omega |\nabla(\sigma^k)^\frac{q}{2}|^2\,\mathrm{d}x
+  \left( C(q-1)\| \nabla \psi^k\|_{\bm{L}^\infty(\Omega)} + b^*\right) \int_\Omega  (\sigma^k)^q\,\mathrm{d}x.
\label{qtot1}
\end{align}
Since $\| \nabla \psi^k\|_{\bm{L}^\infty(\Omega)}\leq C_k\|\psi^k\|$ is bounded, we can apply the Moser--Alikakos iteration argument to the differential inequality \eqref{qtot1}  to conclude that (see, e.g., the calculations in \cite[Theorem 3.1]{Xiang15})
$$
\|\sigma^k(t)\|_{L^\infty(\Omega)}\leq C,\qquad \forall\, t\in [0,T_{\max}),
$$
where $C>0$ depends on $\|\sigma_{0,n}\|_{L^\infty(\Omega)}$, $T$ and $k$, but not on $T_{\max}$. As a consequence, we can conclude that $T_{\max}=T$ and the estimates  \eqref{3eeenerg-a}, \eqref{4eeenerg-a} hold for all  $t\in [0,T]$.

The proof of Lemma \ref{fp} is complete.
\hfill $\square$

\subsection{Proof of Lemma \ref{NSSa}}
\label{app-2}
In view of equation \eqref{g4.d}, we find that $\mu^k$ (i.e., $\{c_i^{k}\}_{i=1}^{k}$) can be uniquely determined by $\varphi^k$ and $\sigma^k$. Hence, problem \eqref{aatest3.c}--\eqref{aatest3.cini} can be reduced to a system with $2k$ nonlinear ordinary differential equations for the time-dependent coefficients $\{a_{i}^{k}\}_{i=1}^{k}$ (by taking $\bm{\zeta}=\bm{y}_i$, $i=1,\cdots,k$) and $\{b_{i}^{k}\}_{i=1}^{k}$ (by taking $w=z_i$, $i=1,\cdots,k$). From the assumptions (H2)--(H4), the construction of $\Psi_\epsilon$  and the regularity property of $\sigma^k$, we can apply the Cauchy--Lipschitz theorem for nonlinear ODE systems to conclude the existence and uniqueness of local solutions $\{a_{i}^{k}\}_{i=1}^k$, $\{b_{i}^{k}\}_{i=1}^{k}\subset  C^1([0,T_k])$, which satisfy the above mentioned ODE system on a certain time interval $[0,T_k]\subset[0,T]$.
Therefore, we obtain a unique local solution $(\bm{v}^k,\varphi^k,\mu^k)$ to problem \eqref{aatest3.c}--\eqref{aatest3.cini} that satisfies
$$
\bm{v}^k\in C^1([0,T_{k}];\bm{Y}_{k}),
\quad
\varphi^k\in C^1([0,T_{k}];Z_{k}),
\quad
\mu^k\in C([0,T_{k}];Z_{k}).
$$

We now show that the existence time $T_k$ can be extended to $T$.
Testing \eqref{aatest3.c} by $\bm{v}^{k}$, using the incompressibility condition and integration by parts, we have
\begin{align}
&\frac12 \frac{\mathrm{d}}{\mathrm{d}t}\|\bm{v}^k\|^2 +\int_{\Omega}2\eta(\varphi^{k} )|D\bm{v}^{k}|^2 \, \mathrm{d}x
  = \big( \mu^{k}\nabla \varphi^{k},\bm{v}^{k}\big)
+\big(\chi \sigma^{k} \nabla \varphi^{k},\bm{v}^{k}\big),
\label{aenergy}
\end{align}
for any $t\in(0,T_k)$.
By H\"{o}lder's and Poincar\'{e}'s inequality, we get
\be
\begin{aligned}
\big(\chi \sigma^{k} \nabla \varphi^{k},\bm{v}^{k}\big)
&  \le  |\chi| \| \sigma^{k}\|  \|\nabla\varphi^{k}\|_{\bm{L}^{\infty}(\Omega)}\|\bm{v}^{k}\|
\\
&  \le  C  \| \sigma^{k}\|  \| \varphi^{k}\|_{H^3(\Omega)}  \|\bm{v}^{k}\|
\\
& \le C_k\| \sigma^{k}\| \left( \|\nabla \varphi^{k}\|^2 +\|\bm{v}^{k}\|^2 +1 \right).
\end{aligned}
\notag
\ee
 Testing \eqref{g1.a} by $\mu^{k}$, we find
\begin{align}
& \frac{\mathrm{d}}{\mathrm{d}t}
\left(\frac{1}{2}\|\nabla \varphi^{k}\|^2+ \int_\Omega  \Psi_\epsilon(\varphi^{k})\, \mathrm{d}x
 + \int_\Omega\frac{\gamma^{8}}{4} |\nabla \varphi^k|^4\, \mathrm{d}x\right)
 + \int_\Omega m(\varphi^k)|\nabla \mu^{k}|^2\,\mathrm{d}x
\notag \\
&\quad  = -\big(\bm{v}^{k} \cdot \nabla \varphi^{k},\mu^{k}\big)
+ \int_\Omega\big(-\alpha\varphi^k+h(\varphi^k,\sigma^k)\big) \mu^k\,\mathrm{d}x
+\chi \int_\Omega \partial_t\varphi^{k} \sigma^{k}\, \mathrm{d}x,
\label{menergy1}
\end{align}
for any $t\in(0,T_k)$. We note that the first terms on the right-hand side of \eqref{aenergy} and \eqref{menergy1} cancel each other.
Next, we infer from \eqref{es-mass2} that
\begin{align}
&\int_{\Omega} \big(-\alpha\varphi^k +h(\varphi^k,\sigma^k)\big)\mu^k\,\mathrm{d}x
\notag \\
& \quad
\leq \frac{m_*}{2}\|\nabla \mu^k\|^2 + C\int_{\Omega}\left(\frac{1}{2}|\nabla \varphi^k|^{2}
+ \Psi_\epsilon(\varphi^k)\right) \mathrm{d} x
+ C\|\sigma^k\|^2 +C.
\label{mueea1}
\end{align}
Concerning the last term on the right-hand side of \eqref{menergy1}, we take
$\xi=\partial_t\varphi^k$ in \eqref{g1.a} and apply (H3), (H4) obtaining
\begin{align}
\|\partial_t\varphi^k\|^2
&=  ( \bm{v}^{k}  \varphi^{k},\nabla \partial_t\varphi^k) -\big(m(\varphi^k)\nabla \mu^{k},\nabla \partial_t\varphi^k\big)
+\big(-\alpha \varphi^k + h(\varphi^k,\sigma^k),\partial_t\varphi^k\big)
\notag\\
&\leq \|\bm{v}^{k}\|_{\bm{L}^6(\Omega)}\|\varphi^{k}\|_{L^3(\Omega)}\|\nabla \partial_t\varphi^k\| + m^* \|\nabla \mu^{k}\|\|\nabla \partial_t\varphi^k \|
+ \big(\alpha\|\varphi^k\|+h^*|\Omega|^\frac12\big) \|\partial_t\varphi^k \|
\notag \\
&\leq C_k(\|\bm{v}^{k}\| \|\varphi^{k}\| + \|\nabla \mu^{k}\|+ \|\varphi^k\|+1)\|\partial_t\varphi^k \|,
\notag
\end{align}
which yields
\begin{align}
\|\partial_t\varphi^k\|
 \leq C_k(\|\bm{v}^{k}\| \|\varphi^{k}\| + \|\nabla \mu^{k}\|+ \|\varphi^k\|+1).
\notag
\end{align}
Then, by Young's inequality, we have
\begin{align}
\chi \int_\Omega \partial_t\varphi^{k}\sigma^{k} \, \mathrm{d}x
& \le |\chi|\|\partial_t\varphi^{k}\| \|\sigma^{k}\|
\notag \\
&\leq \frac{m_*}{4}\|\nabla \mu^k\|^2
+C_k (1+\|\sigma^k\|^2)+ C_k(1+\|\bm{v}^{k}\|) \|\varphi^{k}\|\|\sigma^k\|.
\label{me0}
\end{align}
Adding \eqref{aenergy} with \eqref{menergy1}, we infer from \eqref{mueea1}, \eqref{me0}, \eqref{averphi1} and the Poincar\'{e}--Wirtinger inequality  that
\begin{align}
& \frac{\mathrm{d}}{\mathrm{d}t}
\left(\frac12\|\bm{v}^{k}\|^2
+\frac{1}{2}\|\nabla \varphi^{k}\|^2
+\int_\Omega \Psi_\epsilon (\varphi^{k})\, \mathrm{d}x
+ \int_\Omega\frac{\gamma^{8}}{4} |\nabla \varphi^k|^4\, \mathrm{d}x \right)
\notag \\
&\qquad
+2\eta_*\int_{\Omega} |D\bm{v}^{k}|^2 \, \mathrm{d}x
+ \frac{m_*}{4}\|\nabla  \mu^{k}\|^2
\notag \\
&\quad  \le C_{1,k}(1+\|\sigma^{k}\|^2)\left(\frac12\|\bm{v}^{k}\|^2
+\frac{1}{2}\|\nabla \varphi^{k}\|^2
+\int_\Omega \Psi_\epsilon(\varphi^{k})\, \mathrm{d}x
\right)
	+C_{1,k}(1+\|\sigma^{k}\|^2).
\label{menergy2}
\end{align}
Moreover, in light of \eqref{PZK} (cf. \eqref{iniE0}), the initial value of the energy can be controlled as follows
\begin{align}
&\int_{\Omega}\left[\frac{1}{2}|\bm{P}_{\bm{Y}_{k}} \bm{v}_{0}|^{2}
+\frac{1}{2}|\nabla \bm{P}_{Z_{k}}\varphi_{0,\gamma}|^{2}
+  \Psi_\epsilon(\bm{P}_{Z_{k}}\varphi_{0,\gamma})
+ \frac{\gamma^{8}}{4} |\nabla \bm{P}_{Z_{k}}\varphi_{0,\gamma}|^4 \right] \mathrm{d}x
\notag\\
&\quad \le C\Big(\|  \bm{v}_{0}\|, \|\varphi_{0}\|_{H^1(\Omega)}, \max_{r\in[-1,1]}|\Psi_0(r)|, \Omega\Big)
 =: C_0.
\label{iniC0}
\end{align}
Exploiting the convexity of $\Psi_{0,\epsilon}$, we have
\be
\frac{1}{2}\|\nabla \varphi^{k}\|^2
+\int_\Omega \Psi_\epsilon(\varphi^{k})\, \mathrm{d}x
\ge \frac{1}{2}\|\varphi^{k}\|^2-C_{0,k},
\label{me2}
\ee
for some $C_{0,k}\geq 0$.
Hence, applying Gronwall's lemma to \eqref{menergy2}, we deduce from \eqref{me2} that
\be
\begin{aligned}
&\|\bm{v}^{k}(t)\|^2+ \|\varphi^{k}(t)\|^2
\le 2 \mathrm{e}^{N_1(t)}\big( C_{0} + N_1(t)\big)
+ 2C_{0,k}=:N_2(t),
\quad \forall\, t\in [0,T_k],
\label{auvm1}
\end{aligned}
\ee
 with
\be
\begin{aligned}
N_1(t)&=t C_{1,k}\Big(1+ \sup_{t\in[0,T]} \|\sigma^{k}(t)\|^2\Big).
\notag
\end{aligned}
\ee
Thanks to \eqref{auvm1}, a further integration in time of \eqref{menergy2} gives
\begin{equation}
\int_0^t \| \nabla  \mu^k(s)\|^2 \, \mathrm{d}s \leq
C_0 + C_{0,k} + N_1(t)N_2(t) + N_1(t)=:N_3(t), \quad \forall \, t \in [0,T_k].\notag
\end{equation}
Based on the above estimates and \eqref{3eeenerg-a},
we can extend the unique local solution $(\bm{v}^k,\varphi^k)$ to the whole interval $[0,T]$,
with the same estimate as \eqref{auvm1}.
The boundedness of $(\bm{v}^k,\varphi^k)$ in $H^1(0,T;\bm{Y}_{k})\times H^1(0,T;Z_{k})$
can be easily derived from \eqref{aatest3.c} and \eqref{g1.a} by using the energy method.
Indeed, it is easily seen that
\begin{align*}
\int_0^T \| \partial_t \bm{v}^k(s)\|^2 \, \mathrm{d}s
&\leq C_k \int_0^T \left( \| \bm{v}^k(s)\|^4+ \| \bm{v}^k(s)\|^2+
\| \nabla  \mu^k(s)\|^2 \| \varphi^k(s)\|^2
+\| \sigma^k(s)\|^2 \| \varphi^k(s)\|^2  \right) \mathrm{d}s
\\
&\leq C_k \left( T N^2_2(T) + TN_2(T) + N_3(T) N_2(T) + T N_2(T) \sup_{t\in [0,T]}\|\sigma^{k}(t)\|^2  \right)
\end{align*}
and
\begin{align*}
\int_0^T \| \partial_t \varphi^k(s)\|^2 \, \mathrm{d}s
&\leq C_k \int_0^T \left( \| \bm{v}^k(s)\|^2 \| \varphi^k(s)\|^2
+\| \nabla  \mu^k(s)\|^2 + \| \varphi^k(s)\|^2 + C \right) \mathrm{d}s
\\
&\leq C_k \big(  T N^2_2(T) + N_3(T) + TN_2(T) + C T \big).
\end{align*}
The proof of Lemma \ref{NSSa} is complete.
\hfill $\square$

\section*{Acknowledgments}
\noindent
Part of this work was conducted during A. Giorgini's visit to Fudan University in 2023 under the ``Fudan Fellow" program. A. Giorgini is supported by the MUR grant Dipartimento di Eccellenza 2023--2027.
{\color{black}{J.-N. He’s research is partially supported by Zhejiang Provincial Natural Science Foundation of China (Grant Number LQ24A010011). She also acknowledges the support by the RFS Grant from the Research Grants Council (Project P0047825, PI: Prof. Xianpeng Hu) for her research stay in Research Center for Nonlinear Analysis, The Hong Kong Polytechnic University.}} H. Wu is partially supported by NNSFC 12071084
and the Shanghai Center for Mathematical Sciences.
A. Giorgini is a member of Gruppo Nazionale per l'Analisi Matematica, la Probabilit\'{a} e le loro Applicazioni (GNAMPA), Istituto Nazionale di Alta Matematica (INdAM).
H. Wu is a member of the Key Laboratory of Mathematics for Nonlinear Sciences (Fudan University), Ministry of Education of China.


\end{document}